\documentclass[dvipdfmx,12pt]{article}
\usepackage{latexsym}
\usepackage{amsmath}
\usepackage{amssymb}

\usepackage{color}

\input cyracc.def

\newfont{\tencyr}{wncyr10}

\begin{document}

\begin{center}
 {\Large \bf  On arithmetic Dijkgraaf-Witten theory}
 \end{center}

\vspace{.4cm}

\begin{center}

Hikaru HIRANO, Junhyeong KIM and Masanori MORISHITA

\end{center}

\begin{center}
{\it Dedicated to the memory of Professor Toshie Takata}
\end{center}

\vspace{.8cm}

\footnote[0]{2010 Mathematics Subject Classification: 11R,  57M, 81T. \\
Key words: arithmetic Chern-Simons theory, arithmetic topology, Dijkgraaf-Witten theory, topological quantum field theory}
\\
{\small
Abstract:  We present  basic constructions and properties in arithmetic Chern-Simons theory with finite gauge group along the line of topological quantum field theory. For a finite set $S$ of finite primes of a number field $k$, we construct arithmetic analogues of  the Chern-Simons 1-cocycle, the prequantization bundle for a surface and the Chern-Simons functional for a $3$-manifold. We then construct arithmetic analogues for $k$ and $S$ of the quantum Hilbert space (space of conformal blocks) and the Dijkgraaf-Witten partition function in (2+1)-dimensional Chern-Simons TQFT.  We show some basic and functorial properties of those arithmetic analogues. Finally we show decomposition  and gluing formulas for arithmetic Chern-Simons invariants and arithmetic Dijkgraaf-Witten partition functions.}
\vspace{0.8cm}
\\
{\small \bf 
\begin{center}
Contents.
\end{center}
Introduction\\
1. Preliminaries on torsors and group cochains\\
2. Classical theory\\
2.1. Arithmetic prequantization bundles and arithmetic Chern-Simons 1-cocycles\\
2.2. Arithmetic Chern-Simons functionals\\
3. Quantum theory\\
3.1. Arithmetic quantum spaces\\
3.2. Arithmetic Dijkgraaf-Witten partition functions\\
4. Some basic and functorial properties\\
4.1 Change of the 3-cocycle $c$\\
4.2. Change of number fields\\
4.3. The case that $S$ is empty\\
4.4. Disjoint union of finite set of primes and reversing the orientation of $\partial V_S$\\
5. Decomposition and gluing  formulas\\
5.1. Arithmetic Chern-Simons functionals and arithmetic Dijkgraaf-Witten partition functions for $V_S$\\
5.2. Decomposition formulas for arithmetic Chern-Simons invariants and gluing formulas for arithmetic Dijkgraaf-Witten partition functions\\
References
}

\vspace{0.8cm}

\begin{center}
{\bf Introduction} 
\end{center}

In [Ki] Minhyong Kim initiated to study {\em arithmetic Chern-Simons theory} for number rings, which is based on the ideas of Dijkgraaf-Witten theory for $3$-manifolds ([DW]) and the analogies between 3-manifolds and number rings, knots and primes in arithmetic topology ([Mo2]). We note that Dijkgraaf-Witten theory may be seen as a $3$-dimensional Chern-Simons  gauge theory with finite gauge group (cf. [FQ], [G], [Wa], [Y] etc). Among other things, Kim constructed an arithmetic analog of the Chern-Simons functional, which is defined on a space of Galois representations over a totally imaginary number field. In the subsequent paper  [CKKPY]  Kim and his collaborators showed a decomposition formula for arithmetic Chern-Simons invariants and applied it to concrete computations for some examples.  Later, Kim's construction was extended over arbitrary number field which may have real primes ([H], [LP]). Computations of arithmetic Chern-Simons invariants  have also been carried out for some examples, by employing number-theoretic considerations  in [AC], [BCGKPT], [CKKPPY], [H] and  [LP]. In [CKKPPY], the arithmetic Chern-Simons correlation functions for  finite cyclic gauge groups were computed in terms of arithmetic linking numbers. It should be noted that Kim also considered arithmetic Chern-Simons functionals for the case where the gauge groups are $p$-adic Lie groups ([Ki; 3]). By {\em arithmetic Dijkgraaf-Witten theory} in the title, we mean arithmetic Chern-Simons theory with finite gauge group in the sense of Kim.\\

The purpose of this paper is to add some basic constructions and properties to Kim's theory and lay a foundation for arithmetic Dijkgraaf-Witten theory   along the line of topological quantum field theory, TQFT for short, in the sense of Atiyah ([A1]).  TQFT is a framework to produce topological invariants for manifolds. For example, the Jones polynomials of knots can be obtained in the context of $(2 + 1)$-dimensional Chern-Simons TQFT with compact connected gauge group (cf. [A2], [Ko], [Wi]).  For the TQFT structure of Dijkgraaf-Witten theory, we consult [DW], [FQ], [G], [Wa], [Y].  In this paper, following Gomi's treatment [G] and Kim's original ideas [Ki], we construct an arithmetic analogue of Dijkgraaf-Witten TQFT in a certain special situation, namely, we construct arithmetic analogues,  for a finite set $S$ of finite primes of a number field $k$,  of the prequantization bundles, the Chern-Simons $1$-cocycle,  the Chern-Simons functional,  the quantum Hilbert space (space of conformal blocks) and the Dijkgraaf-Witten partition function. Arithmetic Dijkgraaf-Witten invariants are new arithmetic invariants for a number field, which may be seen as variants of  (non-abelian) Gaussian sums.\\

We fix a finite group $G$ and a 3-cocycle $c \in Z^3(G,\mathbb{R}/\mathbb{Z})$. For an oriented compact manifold $X$ with a fixed triangulation, let ${\cal F}_X$ be the space of gauge fields  associated to $G$ and let ${\cal G}_X$ be the gauge group ${\rm Map}(X,G)$ acting on ${\cal F}_X$. Note that ${\cal F}_X$ and ${\cal G}_X$ are finite sets and that the quotient space ${\cal M}_X := {\cal F}_X/{\cal G}_X$ is identified with  ${\rm Hom}(\pi_1(X), G)/G$ if $X$ is connected, where  ${\rm Hom}(\pi_1(X), G)/G$ is the quotient of the set of homomorphisms from the fundamental group 
$\pi_1(X)$ of $X$ to $G$ by the conjugate action of $G$. 

As for the classical theory in the sense of physics, we construct, using the 3-cocycle $c$, the following correspondences 
{\small $$ \begin{array}{ccc}
\mbox{oriented closed surface}\; \Sigma  & \rightsquigarrow & \mbox{1-cocycle} \; \lambda_{\Sigma}  \in Z^1({\cal G}_{\Sigma}, {\rm Map}({\cal F}_{\Sigma}, \mathbb{R}/\mathbb{Z})), \\
\mbox{oriented compact 3-manifold}\; M  & \rightsquigarrow &  \mbox{ 0-chain} \; CS_{M}  \in C^0({\cal G}_M, {\rm Map}({\cal F}_M, \mathbb{R}/\mathbb{Z})),
\end{array} \leqno{(0.1)}$$ 
}
which satisfy
 $$d CS_M = {\rm res}^* \lambda_{\partial M}, \leqno{(0.2)}$$
where ${\rm res} : {\cal F}_M \; ({\rm resp.} \; {\cal G}_M) \rightarrow {\cal F}_{\partial M} \; ({\rm resp.}\; {\cal G}_{\partial M})$ is the restriction map and $d : C^0({\cal G}_M, {\rm Map}({\cal F}_M, \mathbb{R}/\mathbb{Z}) \rightarrow C^1({\cal G}_M, {\rm Map}({\cal F}_M, \mathbb{R}/\mathbb{Z}))$ is the coboundary homomorphism of group cochains. The key ingredient to construct $\lambda_{\Sigma}$ and $CS_M$ is the transgression homomorphism $C^i(G, \mathbb{R}/\mathbb{Z}) \rightarrow C^{i-d}({\cal G}_X, {\rm Map}({\cal F}_X, \mathbb{R}/\mathbb{Z}))$ with $d = \dim X$ and, in fact, $\lambda_{\Sigma}$ and $CS_M$ are given by the images of $c$ for $i=3$, $X = \Sigma$ and $M$, respectively ([G]). Then we can construct a ${\cal G}_{\Sigma}$-equivariant principal $\mathbb{R}/\mathbb{Z}$-bundle ${\cal L}_{\Sigma}$ and the associated complex line bundle $L_{\Sigma}$ over  ${\cal F}_{\Sigma}$, using $ \lambda_{\Sigma}$, and hence the complex line bundle $\overline{L}_{\Sigma}$  over ${\cal M}_X$. In fact, ${\cal L}_{\Sigma}$ is the product bundle ${\cal F}_{\Sigma} \times \mathbb{R}/\mathbb{Z}$ on which ${\cal G}_{\Sigma}$ acts by $(\rho_{\Sigma},m).g = (\rho_{\Sigma}.g, m + \lambda_{\Sigma}(g,\rho_{\Sigma}))$ for $\rho_{\Sigma} \in {\cal F}_{\Sigma}, m \in \mathbb{R}/\mathbb{Z}$ and $g \in {\cal G}_{\Sigma}$. We call $\lambda_{\Sigma}$ the {\em Chern-Simons $1$-cocycle}. The line bundle $L_{\Sigma}$ (or $\overline{L}_{\Sigma}$) is called the  {\em prequantization complex line bundle} for a surface $\Sigma$. The 0-chain $CS_M$  is called the  {\em Chern-Simons functional}  for a $3$-manifold $M$. We see that $CS_M$ is a ${\cal G}_M$-equivariant section of ${\rm res}^*{\cal L}_{\Sigma}$ over ${\cal F}_M$. 

As for the quantum theory,  the formalism of $(2+1)$-dimensional TQFT is given by the following correspondences (functor from the cobordism category of surfaces to the category of complex vector spaces) 
$$ \begin{array}{ccc}
\mbox{oriented closed surface}\; \Sigma \; & \rightsquigarrow & \; \mbox{quantum Hilbert space} \; {\cal H}_{\Sigma}, \\
\mbox{oriented 3-manifold}\; M \; & \rightsquigarrow & \; \mbox{partition function} \; Z_M \in {\cal H}_{\partial M},
\end{array} \leqno{(0.3)}$$
which satisfy several axioms (cf. [A1]). Here we notice the following two axioms:\\
$(0.4)$  {\em functoriality}:  An orientation preserving homeomorphism $f : \Sigma \stackrel{\approx}{\rightarrow} \Sigma'$ induces an isomorphism ${\cal H}_{\Sigma}  \stackrel{\sim}{\rightarrow} {\cal H}_{\Sigma'}$ of Hilbert quantum spaces. Moreover, if $f$ extends to an orientation preserving homeomorphism $M \stackrel{\approx}{\rightarrow} M'$, with $\partial M = \Sigma, \partial M' = \Sigma'$, then $Z_{M}$ is sent to $Z_{M'}$ under the induced isomorphism ${\cal H}_{\partial M} \stackrel{\sim}{\rightarrow} {\cal H}_{\partial M'}$.\\
$(0.5)$  {\em multiplicativity and involutority}:  For disjoint surfaces $\Sigma_1, \Sigma_2$ and the surface $\Sigma^* $ = $\Sigma$ with the opposite orientation, we require
$${\cal H}_{\Sigma_1 \sqcup \Sigma_2} = {\cal H}_{\Sigma_1} \otimes {\cal H}_{\Sigma_2}, \;\; {\cal H}_{\Sigma^*} = ({\cal H}_{\Sigma})^*,$$
 where $({\cal H}_{\Sigma})^*$ is the dual space of ${\cal H}_{\Sigma}$. Moreover, if $\partial M_1 = \Sigma_1 \sqcup \Sigma_2, \partial M_2 = \Sigma_2^* \sqcup \Sigma_3$ and $M$ is the $3$-manifold obtained by gluing $M_1$ and $M_2$ along $\Sigma_2$, then we require 
$$ < Z_{M_1}, Z_{M_2}> \, = \, Z_M, $$ 
where $<\cdot , \cdot > : {\cal H}_{\Sigma_1 \sqcup \Sigma_2} \times {\cal H}_{\Sigma_2^* \sqcup \Sigma_3} \rightarrow {\cal H}_{\Sigma_1 \sqcup \Sigma_3}$ is the natural  gluing pairing of  quantum Hilbert spaces. This multiplicative property is indicative of the ``quantum" feature of the theory (cf. [A1]). 

The construction of the Hilbert space ${\cal H}_{\Sigma}$ is phrased as the {\it geometric quantization}. We note that  ${\cal H}_{\Sigma}$ is known to be isomorphic to the space of conformal blocks for the surface $\Sigma$ when the gauge group is a compact connected group (cf. [Ko]). Elements of ${\cal H}_{\Sigma}$ are called (non-abelian) theta functions (cf. [BL]). For Dijkgraaf-Witten theory, ${\cal H}_{\Sigma}$ is constructed, in an analogous manner,  as the space of ${\cal G}_{\Sigma}$-equivariant sections of the prequantization line bundle $L_{\Sigma}$ over ${\cal F}_{\Sigma}$, in other words, the space of sections of $\overline{L}_\Sigma$ over ${\cal M}_\Sigma$:
$$\begin{array}{ll} 
 {\cal H}_{\Sigma} & = \{ \vartheta : {\cal F}_{\Sigma} \rightarrow \mathbb{C}\, | \, \vartheta(\varrho_{\Sigma}.g) = e^{2\pi\sqrt{-1}\lambda_{\Sigma}(g)(\vartheta)} \vartheta(\varrho_{\Sigma})\; \forall g \in {\cal G}_{\Sigma}, \, \varrho_{\Sigma} \in {\cal F}_{\Sigma} \}\\
                   & = \Gamma({\cal M}_{\Sigma},\overline{L}_{\Sigma}).\\
\end{array}  \leqno{(0.6)}$$
In quantum field theories, partition functions are given as path integrals. In Dijkgraaf-Witten theory,
the {\em Dijkgraaf-Witten partition function} $Z_M \in {\cal H}_{\partial M}$ is defined by the following finite sum fixing the boundary condition:
$$ \displaystyle{ Z_M(\varrho_{\partial M}) = \frac{1}{\# G} \sum_{{\scriptstyle \varrho \in {\cal F}_M}\atop{\scriptstyle {\rm res}(\varrho) = \varrho_{\partial M}}}   e^{2\pi \sqrt{-1} CS_M(\varrho)} \;\; \;\; (\varrho_{\partial M} \in {\cal F}_{\partial M}).} \leqno{(0.7)}$$
The value $Z_M(\varrho_{\partial M})$ is called the {\em Dijkgraaf-Witten invariant} of $\varrho_{\partial M} \in {\cal F}_{\partial M}$. We note that when $[c]$ is trivial and $S$ is empty, then ${\cal F}_{\Sigma} = \{ * \}$ and the Dijkgraaf-Witten invariant $Z_M(*)$, denoted by $Z(M)$,  coincides with the (averaged) number of homomorphism from $\pi_1(M)$ to $G$:
$$ Z(M) = \frac{ \# {\rm Hom}(\pi_1(M), G)}{\# G}, \leqno{(0.8)}$$
which is the classical invariant for the connected 3-manifold $M$.\\

 Now let us turn to the arithmetic. First, let us recall the basic analogies in arithmetic topology which bridges $3$-dimensional topology and number theory ([Mo2]. See also [Mh], [NU]). Let $k$ a number field of finite degree over the rationals $\mathbb{Q}$. Let ${\cal O}_k$ be the ring of integers of $k$ and set $X_k := {\rm Spec}({\cal O}_k)$. Let $X_k^{\infty}$ denote the set of infinite primes of $k$ and set $\overline{X}_k := X_k \sqcup X_k^{\infty}$.  We see  $X_k$, $X_k^{\infty}$ and $\overline{X}_k$  as analogues of a non-compact $3$-manifold $M$, the set of ends and the end-compactification $\overline{M}$, respectively. A maximal ideal $\frak{p}$ of ${\cal O}_k$ is identified with the residue field ${\rm Spec}({\cal O}_k/\frak{p}) = K(\widehat{\mathbb{Z}},1)$ ($\widehat{\mathbb{Z}}$ being the profinite completion of $\mathbb{Z}$), which is seen as an analogue of the circle $S^1 = K(\mathbb{Z},1)$. We see the mod $\frak{p}$  reduction map ${\rm Spec}(\mathbb{F}_{\frak{p}}) \hookrightarrow X_k$ as an analogue of a knot, an embedding $S^1  \hookrightarrow M$. Let ${\cal O}_{\frak{p}}$ be the ring of $\frak{p}$-adic integers and let $k_{\frak{p}}$ be the $\frak{p}$-adic field. We denote  ${\rm Spec}({\cal O}_\frak{p})$ and ${\rm Spec}(k_\frak{p})$ by  $V_\frak{p}$ and $\partial V_\frak{p}$, respectively . We see $V_\frak{p}$ and $\partial V_\frak{p}$ as analogue of a tubular neighborhood of a knot and its boundary torus, respectively. So  we see the \'{e}tale fundamental group $\Pi_{\frak{p}}$ of ${\rm Spec}(k_{\frak{p}})$, which is the absolute Galois group ${\rm Gal}(\overline{k}_{\frak{p}}/k_{\frak{p}})$ ($\overline{k}_{\frak{p}}$ being an algebraic closure of $k_{\frak{p}}$),  as an analogue of the peripheral group of a knot. (To be precise, the tame quotient of $\Pi_{\frak{p}}$ may be seen as a closer  analogue of the peripheral group. cf. [Mo2; Chapter 3]) 

 Let $S = \{ \frak{p}_1, \dots , \frak{p}_r \}$ be a finite set of maximal ideals of ${\cal O}_k$. Let $\overline{X}_S := \overline{X}_k \setminus S$. We see $S$   
 and $\overline{X}_S$  as an analogue of a link in a $3$-manifold and the link complement, respectively. We may also see $\overline{X}_S$ as an analogue of a compact $3$-manifold with boundary (union of tori), where  $\partial V_S := {\rm Spec}(k_{\frak{p}_1}) \sqcup \cdots \sqcup {\rm Spec}(k_{\frak{p}_r})$ plays an analogous role of the boundary tori, ``$\partial \overline{X}_S = \partial V_S$". The modified \'{e}tale fundamental group $\Pi_S$ of $\overline{X}_S$, introduced in [H; 2.1] taking real primes into account, is
 the Galois group of the maximal subextension $k_S$ of $k$ which is unramified at any (finite and infinite) prime outside $S$, as an analogue of the link group.

We list herewith  some analogies which will be used in this paper.
\vspace{0.8cm}\\
\begin{tabular}{|c|c|}
\hline
oriented, connected, closed  & compactified spectrum of  \\
3-manifold $\overline{M}$ & number ring $\overline{X}_k = \overline{{\rm Spec}({\cal O}_k)}$ \\
\hline
knot & prime \\
 $ {\cal K} : S^{1} \hookrightarrow M$ &   $\{ \frak{p} \} \ = {\rm Spec}({\cal O}_k/\frak{p}) \hookrightarrow \overline{X}_k$\\ 
 \hline 
 link & finite set of maximal ideals\\
 ${\cal L} = {\cal K}_1 \sqcup \cdots \sqcup {\cal K}_r$ & $S = \{ \frak{p}_1, \dots , \frak{p}_r\}$\\
 \hline 
 tubular n.b.d of a knot & $\frak{p}$-adic integer ring \\
 $V_{\cal K}$ & $V_\frak{p} = {\rm Spec}({\cal O}_{\frak{p}})$\\
 boundary torus & $\frak{p}$-adic field\\
 $\partial V_{\cal K} $ &
 $\partial V_\frak{p} = {\rm Spec}(k_{\frak{p}})$ \\
 peripheral group & local absolute Galois group\\
 $\pi_1(\partial V_{\cal K})$ &  $\Pi_{\frak{p}} = {\rm Gal}(\overline{k}_{\frak{p}}/k_{\frak{p}})$ \\
\hline
tubular n.b.d of a link & union of $\frak{p}_i$-adic integer rings \\
$V_{\cal L} = V_{{\cal K}_1}\sqcup \cdots \sqcup  V_{{\cal K}_r}$ & $V_S = {\rm Spec}({\cal O}_{\frak{p}_1}) \sqcup \cdots \sqcup {\rm Spec}({\cal O}_{\frak{p}_r})$\\
boundary tori &   union of $\frak{p}_i$-adic fields \\
$\partial V_{\cal L} = \partial V_{{\cal K}_1}\sqcup \cdots \sqcup \partial V_{{\cal K}_r}$ & $ \partial V_S = {\rm Spec}(k_{\frak{p}_1}) \sqcup \cdots \sqcup {\rm Spec}(k_{\frak{p}_r})$\\
\hline
link complement & complement of a finite set of primes\\
 $X_{\cal L} = \overline{M} \setminus {\rm Int}(V_{\cal L})$ & $\overline{X}_S = \overline{X}_k \setminus S$\\
 link group & maximal Galois group with \\
 $\Pi_{\cal L} = \pi_1(X_{\cal L})$ & given ramification $\Pi_S = {\rm Gal}(k_S/k)$\\
 \hline
\end{tabular}
\vspace{.8cm}\\
 Based on the above analogies, we construct an arithmetic analogue of Dijkgraaf-Witten TQFT in a special situation, which corresponds to the case that $M$ is a link complement and $\Sigma$ is the boundary tori of a tubular neighborhood of a link. Notations being as above, let $N$ be an integer $>1$ and assume that the number field $k$ contains a primitive $N$-th root $\zeta_N$ of unity. We fix a finite group  $G$ and  a 3-cocycle $c \in Z^3(G,\mathbb{Z}/N\mathbb{Z})$.  Let $F$ be a subfield of $\mathbb{C}$ such that $\zeta_N$ is contained in $F$ and $\overline{F} = F$ ($\overline{F}$ being the complex conjugate). Let $S$ be a finite set of finite primes $S = \{ \frak{p}_1, \dots , \frak{p}_r \}$ of $k$ such that any finite prime dividing $N$ is contained in $S$. Let $\overline{X}_S := \overline{X}_k \setminus S$ and let $\partial V_S := {\rm Spec}(k_{\frak{p}_1}) \sqcup \cdots \sqcup {\rm Spec}(k_{\frak{p}_r})$ as before so that $\partial V_S$ plays a role of the boundary of $\overline{X}_S$, ``$\partial \overline{X}_S = \partial V_S$". For arithmetic analogues of the spaces of gauge fields ${\cal F}_{\Sigma}$ and ${\cal F}_M$, we consider ${\cal F}_{S}:= \prod_{i=1}^r {\rm Hom}_{\rm cont}(\Pi_{\frak{p}_i}, G)$ and ${\cal F}_{\overline{X}_S} := {\rm Hom}_{\rm cont}(\Pi_S,G)$, respectively, where ${\rm Hom}_{\rm cont}(-, G)$ denotes the set of continuous homomorphisms to $G$. For an arithmetic analog of the gauge groups  ${\cal G}_{\Sigma}$ and ${\cal G}_M$, we simply take the group $G$ acting on ${\cal F}_S$ and ${\cal F}_{\overline{X}_S}$ by conjugation. Set ${\cal M}_S := {\cal F}_S/G$. 

As for the classical theory in the arithmetic side, we firstly develop a local theory at a finite prime $\frak{p}$, namely, we construct the {\em arithmetic prequantization principal $\mathbb{Z}/N\mathbb{Z}$-bundle} ${\cal L}_{\frak{p}}$ and the associated {\em arithmetic prequantization $F$-line bundle} $L_{\frak{p}}$ for $\partial V_{\frak{p}}$, which are $G$-equivariant  bundles over ${\cal F}_{\frak{p}} := {\rm Hom}_{\rm cont}(\Pi_{\frak{p}}, G)$. By choosing a section $x_{\frak{p}} \in \Gamma({\cal F}_{\frak{p}},{\cal L}_{\frak{p}})$, we construct  the {\em arithmetic Chern-Simons $1$-cocycle} $\lambda_{\frak{p}}^{x_{\frak{p}}} \in Z^1(G, {\rm Map}({\cal F}_{\frak{p}}, \mathbb{Z}/N\mathbb{Z}))$.  The key idea for the constructions is due to M. Kim ([Ki]), who used the conjugate $G$-action on $c$ and the canonical isomorphism
$$ {\rm inv}_{\frak{p}} : H^2(\Pi_{\frak{p}}, \mathbb{Z}/N\mathbb{Z}) \stackrel{\sim}{\longrightarrow} \mathbb{Z}/N\mathbb{Z}$$
in the theory of Brauer groups of  local fields. We note that this isomorphism tells us that $\partial V_{\frak{p}}$ is ``orientable" and we choose (implicitly) the ``orientation"  of $\partial V_{\frak{p}}$ corresponding to $1 \in \mathbb{Z}/N\mathbb{Z}$. 

Getting together the local theory over $S$, we construct  the {\em arithmetic prequantization principal $\mathbb{Z}/N\mathbb{Z}$-bundle} ${\cal L}_S$ and the associated {\em arithmetic prequantization $F$-line bundle} $L_S$ for $\partial V_S$, which are
$G$-equivariant  bundles over ${\cal F}_S$. By choosing a section $x_S$ 
of ${\cal L}_S$ over ${\cal F}_{S}$, we construct  the {\em arithmetic Chern-Simons $1$-cocycle} $\lambda_S^{x_S} \in Z^1(G, {\rm Map}({\cal F}_S, \mathbb{Z}/N\mathbb{Z}))$ and show that ${\cal L}_S$ (resp. $L_S$) is isomorphic to the product bundle ${\cal L}_S^{x_S} = {\cal F}_S \times \mathbb{Z}/N\mathbb{Z}$ (resp. $L_S^{x_S} = {\cal F}_S \times F$) on which $G$ acts by $(\rho_S, m).g = (\rho_S.g, m + \lambda_S^{x_S}(g, \rho_S))$ (resp. $(\rho_S, z).g = (\rho_S.g, z \zeta_N^{\lambda_S^{x_S}(g, \rho_S)})$) for $\rho_S \in {\cal F}_S$, $m \in \mathbb{Z}/N\mathbb{Z}, z \in F$ and $g \in G$. By employing $H^3(\Pi_S, \mathbb{Z}/N\mathbb{Z}) = 0$, the {\em arithmetic Chern-Simons functional} $CS_{\overline{X}_S}$ for $\overline{X}_S$ is defined as a $G$-equivariant section of ${\rm res}_S^{*}({\cal L}_S)$ over ${\cal F}_{\overline{X}_S}$, where ${\rm res}_S : {\cal F}_{\overline{X}_S} \rightarrow {\cal F}_{S}$ is the restriction map induced by the natural homomorphisms $\Pi_{\frak{p}} \rightarrow \Pi_S$ for $\frak{p} \in S$. Using the section $x_S$, it can be regarded as a $G$-equivariant functional $CS_{\overline{X}_S}^{x_S} : {\cal F}_{\overline{X}_S} \rightarrow \mathbb{Z}/N\mathbb{Z}$. Thus we construct the following correspondences 
$$ \begin{array}{ccc}
\partial V_S  \; & \rightsquigarrow & \; \mbox{1-cocycle} \; \lambda_{S}^{x_S}  \in Z^1(G, {\rm Map}({\cal F}_{S}, \mathbb{Z}/N\mathbb{Z})), \\
\overline{X}_S  \; & \rightsquigarrow & \; \mbox{0-chain} \; CS_{\overline{X}_S}^{x_S}  \in C^0(G, {\rm Map}({\cal F}_{\overline{X}_S}, \mathbb{Z}/N\mathbb{Z})),
\end{array} \leqno{(0.9)}$$ 
which satisfy 
 $$d CS_{\overline{X}_S}^{x_S} = {\rm res}_S^* \lambda_{S}^{x_S}. \leqno{(0.10)}$$
We may regard (0.9), (0.10) as arithmetic analogues of (0.1), (0.2)  in a special situation that corresponds to the case $\Sigma$ is a boundary tori of a link and $M$ is a link complement.  

As for the quantum theory in the arithmetic side, following the topological side, we define the {\em arithmetic quantum space} ${\cal H}_S$ for $\partial V_S$ to be the space of $G$-equivariant sections of the arithmetic prequantization $F$-line bundle $L_S$ over ${\cal F}_S$. Choosing a section $x_S \in \Gamma({\cal F}_S, {\cal L}_S)$, it is isomorphic to the space ${\cal H}_S^{x_S}$ given by
$$\begin{array}{ll} 
 {\cal H}_S^{x_S} & = \{ \theta : {\cal F}_S \rightarrow F \, | \, \theta(\rho_S.g) = \zeta_N^{\lambda_{S}(g)(\rho_S)} \theta(\rho_S)\; \forall g \in G, \, \rho_S \in {\cal F}_S \}\\
                   & = \Gamma({\cal M}_{S},\overline{L}_{S}^{x_S}),\\
\end{array} \leqno{(0.11)}$$
where $\overline{L}_S^{x_S}$ is the quotient of $L_S^{x_S}$ by the action of $G$. 
The {\it arithmetic Dijkgraaf-Witten invariant} $Z_{\overline{X}_S}^{x_S}(\rho_S)$ of $\rho_S \in {\cal F}_S$  with respect to $x_S$ is then defined by the following finite sum fixing the boundary condition:
$$ \displaystyle{ Z_{\overline{X}_S}^{x_S}(\rho_S) = \frac{1}{\# G} \sum_{{\scriptstyle \rho \in {\cal F}_{\overline{X}_S}}\atop{\scriptstyle {\rm res}_S(\rho) = \rho_S}}   \zeta_N^{CS_{\overline{X}_S}^{x_S}(\rho)}.} \leqno{(0.12)} $$
Then we can show that $Z_{\overline{X}_S}^{x_S} \in {\cal H}_S^{x_S}$. Since the spaces ${\cal H}_S^{x_S}$, when $x_S$ is varied, are naturally isomorphic each other, ${\cal H}_S$ is identified with $(\bigsqcup {\cal H}_S^{x_S})/\sim$, where the equivalence relation $\sim$ identifies elements via the isomorphisms between ${\cal H}_S^{x_S}$'s. Hence $Z_{\overline{X}_S}^{x_S}$ determine the element $Z_{\overline{X}_S} \in {\cal H}_S$, which we call the {\it arithmetic Dijkgraaf-Witten partition function} for $\overline{X}_S$. Thus we construct the following correspondences 
$$ \begin{array}{ccc}
 \partial V_S  \; & \rightsquigarrow & \; \mbox{ arithmetic quantum 
 space} \; {\cal H}_{S}, \\
\overline{X}_S   \; & \rightsquigarrow & \; \mbox{arithmetic Dijkgraaf-Witten partition function} \; Z_{\overline{X}_S} \in {\cal H}_S,
\end{array} \leqno{(0.13)}$$
which satisfy some properties similar to the axioms in $(2+1)$-dimensional TQFT. 
We note that when $[c]$ is trivial and $S$ is empty, then the arithmetic Dijkgraaf-Witten invariant $Z_{\overline{X}_S}$, denoted by $Z(\overline{X}_k)$, coincides with the (averaged) number of continuous homomorphism from the modified \'{e}tale fundamental group  $\pi_1(\overline{X}_k)$ of $\overline{X}_k$ ([H;2.1]), which is the Galois group of maximal extension of $k$ unramified at all finite and infinite primes, to $G$:
$$ Z(\overline{X}_k) = \frac{ \# {\rm Hom}_{\rm cont}(\pi(\overline{X}_k),G) }{\# G}, \leqno{(0.14)}$$
which is the classical invariant for the number field $k$.
We may regard (0.11), (0.12), (0.13) and (0.14)  as an arithmetic analogues of (0.6), (0.7), (0.3) and (0.8) respectively, in a special situation that corresponds to the case $\Sigma$ is a boundary tori of a link and $M$ is a link complement. 

We note that elements of ${\cal H}_S$ may be seen as arithmetic analogs of (non-abelian) theta functions. In this respect it may be interesting to observe that the arithmetic Dijkgraaf-Witten invariants $Z_{\overline{X}_S}^{x_S}(\rho_S)$ in (0.12) look like (non-abelian) Gaussian sums.\\

Next we show some basic and functorial properties of arithmetic Chern-Simons 1-cocycles, arithmetic prequantization bundles, arithmetic Chern-Simons invariants, arithmetic quantum spaces and arithmetic Dijkgraaf-Witten partition function \\
(i) when we change the 3-cocycle $c$ in the cohomology class $[c]$,\\
(ii) when we change the pair of $k$ and $S$ to the isomorphic one,  \\
(iii) when $S$ is an empty set, and \\
(iv) when $S$ is a disjoint union of finite sets of finite primes and when we reverse the orientation of $\partial V_S$.\\
As for (ii) and (iv), we show the following properties:\\
(0.15) {\em functoriality}: If there are isomorphisms $\xi_i : k_{\frak{p}_i} \stackrel{\sim}{\rightarrow} k'_{\frak{p}_i'}$ ($1 \leq i \leq r$), then they induce the isomorphism ${\cal H}_S  \stackrel{\sim}{\rightarrow} {\cal H}_{S'}$ for $S = \{ \frak{p}_1, \dots , \frak{p}_r \}, S' = \{ \frak{p}_1', \dots , \frak{p}_r' \}$. Moreover, if $\xi : k \stackrel{\sim}{\rightarrow} k'$ is an isomorphism of number fields such that $\xi(\frak{p}_i) = \frak{p}'_i$ and $\xi$ induces isomorphisms $k_{\frak{p}_i} \stackrel{\sim}{\rightarrow} k'_{\frak{p}'_i}$, then $\xi$ induces the isomorphism ${\cal H}_S \stackrel{\sim}{\rightarrow} {\cal H}_{S'}$ which sends $Z_{\overline{X}_S}$ to $Z_{\overline{X}_{S'}}$.\\
(0.16) {\em multiplicativity and involutority}: For disjoint sets $S_1, S_2$ of finite sets of finite primes and $\partial V_S^{*}$ = $\partial V_S$ with the opposite orientation for a finite set $S$ of finite primes (cf. 4.4 below for the meaning), we show
$$ {\cal H}_{S_1 \sqcup S_2} = {\cal H}_{S_1} \otimes {\cal H}_{S_2}, \;\; {\cal H}_{S^*} = ({\cal H}_S)^*,$$
where ${\cal H}_{S^*}$ denotes the arithmetic quantum space for $\partial V_S^*$ and $({\cal H}_S)^*$ is the dual space of ${\cal H}_S$.\\
These properties (0.15) and (0.16) may be regarded as arithmetic analogues of the axioms (0.4) and (0.5) in $(2+1)$-dimensional TQFT.\\

Finally we show decomposition formulas for arithmetic Chern-Simons invariants, which generalize, in our framework,  the ``decomposition formula" by Kim and his collaborators ([CKKPY]),  and show gluing formulas for arithmetic Dijkgraaf-Witten partition functions. Let $S_1$ and $S_2$ be disjoint sets of finite primes of $k$, where $S_1$ may be empty and $S_2$ is non-empty. We assume that any prime dividing $N$ is contained in $S_2$ if $S_1$ is empty and that any prime dividing of $N$ is contained in $S_1$ if $S_1$ is non-empty. We set $S := S_1 \sqcup S_2$.  When $S_1$ is empty, $\overline{X}_{S_1} = \overline{X}_k$ and we mean by $CS_{\overline{X}_{S_1}}$ the arithmetic Chern-Simons functional $CS_{\overline{X}_k}$ defined in [H] (see also [LP]).  We can also define the arithmetic Chern-Simons functional  $CS_{V_{S_2}}$ for $V_{S_2}$  as a section of $\tilde{\rm res}_{S_2}^{*}({\cal L}_{S_2})$ over ${\cal F}_{V_{S_2}} := \prod_{\frak{p} \in S_2} {\rm Hom}_{\rm cont}(\tilde{\Pi}_{\frak{p}},G)$, where $\tilde{\Pi}_{\frak{p}} := \pi_1^{\mbox{\'{e}t}}(V_{\frak{p}})$ and $\tilde{\rm res}_S : {\cal F}_{V_{S_2}} \rightarrow {\cal F}_{S_2}$ is the restriction map induced by the natural homomorphism $\Pi_{\frak{p}} \rightarrow \tilde{\Pi}_{\frak{p}}$. Then we have the following decomposition formula 
$$CS_{\overline{X}_{S_1}}(\rho) \boxplus CS_{V_{S_2}}((\rho \circ u_{\frak{p}})_{\frak{p} \in S_2} ) = CS_{\overline{X}_S}(\rho \circ \eta_S) , \leqno{(0.17)}$$
 where $\rho \in {\rm Hom}_{\rm cont}(\Pi_{S_1}, G)$, and $\eta_S : \Pi_S \rightarrow \Pi_{S_1}$, $u_{\frak{p}} :  \tilde{\Pi}_{\frak{p}} \rightarrow \Pi_{S_1}$ are natural homomorphisms induced by $\overline{X}_S \rightarrow \overline{X}_{S_1}$, $V_{\frak{p}} \rightarrow \overline{X}_{S_1}$ for $\frak{p} \in S_2$, respectively, and $\boxplus : {\cal L}_{S_1} \times {\cal L}_{S_2} \rightarrow {\cal L}_S$ is the natural ``sum" of arithmetic prequantization principal $\mathbb{Z}/N\mathbb{Z}$-bundles (cf. (4.3.1), (4.3.2)). When $S_1$ is empty, the formula (0.13) is a reformulation of the decomposition formula in [CKKPY]. As for arithmetic Dijkgraaf-Witten partition functions, we have the following gluing formula. Note that $\overline{X}_{S_1}$ may be obtained by gluing $\overline{X}_S$ and $V_{S_2}^*$ along $\partial V_{S_2}$, where $V_{S_2}^* = V_{S_2}$ with the opposite orientation.  Then we have 

$$ < Z_{\overline{X}_S}, Z_{V_{S_2}^{*}} >\;  = \; Z_{\overline{X}_{S_1}}, \leqno{(0.18)}$$ 
where $< \cdot , \cdot > : {\cal H}_S \times {\cal H}_{S_2^*} \rightarrow {\cal H}_{S_1}$ is the gluing pairing of arithmetic quantum spaces (cf. (5.2.3)). We may regard (0.16)  as an arithmetic analog of the gluing formula in the axiom (0.5) in $(2+1)$-dimensional TQFT.\\
The contents of this paper are organized as follows. In Section 1, we collect some basic facts on torsors and group cochains, which will be used in the subsequent sections. In Section 2, we construct arithmetic prequantization bundles, arithmetic Chern-Simons $1$-cocycles  and the arithmetic Chern-Simons functionals. These constructions correspond to the classical theory of topological Dijkgraaf-Witten TQFT. In Section 3, we construct arithmetic quantum spaces and the arithmetic Dijkgraaf-Witten partition functions. These constructions correspond to the quantum theory of topological Dijkgraaf-Witten TQFT. In Section 4, we show some basic and functorial properties of arithmetic prequantization bundles, arithmetic Chern-Simons $1$-cocycles, arithmetic Chern-Simons invariants and arithmetic Dijkgraaf-Witten invariants. In Section 5,  we show decomposition formulas for arithmetic Chern-Simons invariants and gluing formulas for arithmetic Dijkgraaf-Witten partition functions.
\\
\\
{\em Notation.} For a $G$-equivariant fiber bundle $\varpi : E \rightarrow B$ for a group $G$, we denote by $\Gamma(B, E)$ (resp. $\Gamma_G(B,E)$)  the set of  sections (resp. the set of $G$-equivariant sections) of $\varpi$. In this paper, we deal with the case where the base space $B$ is a finite (discrete) set. \\
\\
{\em Acknowledgment.} We would like to thank Kiyonori Gomi, Tomoki Mihara, Yuji Terashima, Masahito Yamazaki and Michihisa Wakui for useful communications. We are grateful to Gomi for answering our questions patiently. We would like to thank the referee for careful reading of the paper and useful comments. The first author is supported by Grant-in-Aid for JSPS Fellow (DC1) Grant Number 20J21684. The second author was supported by Grant-in-Aid for JSPS Fellow (DC1) Grant Number 17J02472. The third author is supported by Grant-in-Aid for Scientific Research (KAKENHI) (B) Grant Number JP17H02837.
\\ 
 \begin{center}
 {\bf 1. Preliminaries on torsors and group cochains}
\end{center}

In this section, we collect some baisc facts on torsors for an additive group and group cochains, which will be used in the subsequent sections.\\
\\
{\bf 1.1. Torsors for an additive group.} Let $A$ be an additive group, where the identity element of $A$ is denoted by $0$. An {\em $A$-torsor} is defined by a non-empty set $T$ equipped with action of $A$ from the right 
$$T \times A \longrightarrow T; \; (t,a) \mapsto t.a,$$
which is simply transitive. So, for any elements $s, t \in T$, there exists uniquely $a \in A$ such that $s = t.a$. We denote such an $a$ by $s - t$:
$$  a = s - t \stackrel{\rm def}{\Longleftrightarrow} s = t.a. \leqno{(1.1.1)}$$

For $A$-torsors $T$ and $T'$, a morphism $ f : T \rightarrow T'$ is defined by a map of sets, which  satisfies 
$$f(t.a) = f(t).a \leqno{(1.1.2)}$$
 for all $t \in T$ and $a \in A$. We easily see that any morphism of $A$-torsors is an isomorphism. 

Defining the action of $A$ on $A$ by $(t,a) \in A \times A \mapsto t+a \in A$, $A$ itself becomes an $A$-torsor. We call it a {\em trivial $A$-torsor}.  A morphism $f : A \rightarrow A$ of trivial $A$-torsors is given by $f(a) = a + \lambda$ for any $a \in A$ with $\lambda = f(0)$. Choosing an element $t \in T$, any $A$-torsor $T$ is isomorphic to the trivial $A$-torsor by the morphism
$$ \varphi_t : T \stackrel{\sim}{\longrightarrow} A; \; s \mapsto \varphi_t(s) := s - t.  \leqno{(1.1.3)}$$
We call $\varphi_t$ the {\it trivialization} at $t$.

Here are some properties concerning $A$-torsors, which will be used in the subsequent sections. \\
\\
{\bf Lemma 1.1.4.} (1) {\it Let $T$ be an $A$-torsor. For $s, t, u \in T$ and $a \in A$, we have the following equality in $A$}:
$$ s - s = 0, \;\; s - u = (s - t) + (t - u), \;\; s.a - t = (s - t) + a.$$
(2) {\em Let $T, T'$ be $A$-torsors and let $f : T \rightarrow T'$ be 
a morphism of $A$-torsors. Then, for $s, t \in T$, we have the following equality in $A$}:
$$ s - t = f(s) - f(t).$$
(3) {\em Let $T, T'$ be $A$-torsors and let $f : T \rightarrow T'$ be 
a morphism of $A$-torsors. Fix $t \in T$ and $t' \in T'$, and let $\lambda(f; t, t') := f(t) - t'$. Then we have the following commutative diagram}:
$$ \begin{array}{rcl}
 T & \stackrel{f}{\longrightarrow}   & T'\\
\varphi_t \downarrow & & \downarrow \varphi_{t'} \\
A &\stackrel{+ \lambda(f; t, t')}{\longrightarrow}& A.
\end{array}
$$
For other choices $s \in T$ and $s' \in T'$, we have 
$$\lambda(f; s, s') = \lambda(f;t,t') + (s - t) - (s' - t').$$
(4) {\it For an $A$-torsor $T$ and a subgroup $B$ of $A$, we note that the quotient set $T/B$ is an $A/B$-torsor by $(t \, \mbox{mod} \, B).(a \, \mbox{mod}\, B) := (t.a \, \mbox{mod}\, B)$ for $t \in T$ and $a \in A$.}\\
\\
{\it Proof.} (1) These equalities follow from the definition of group action and (1.1.1). \\
(2) This follows from (1.1.1) and (1.1.2).\\
(3) The former assertion follows from (1.1.3). For the latter assertion, we note the following commutative diagram.
$$ \begin{array}{ccccccc}
T & \stackrel{{\rm id}}{\longrightarrow} &  T & \stackrel{f}{\longrightarrow} &  T' & \stackrel{{\rm id}}{\longrightarrow}  & T'\\
\;\;\;\; \downarrow \varphi_{s}  & & \;\;\;\; \downarrow \varphi_{t} & & \;\;\;\; \downarrow \varphi_{t'} & & \;\;\;\; \downarrow \varphi_{s'}\\
A & \stackrel{ + (s - t)}{\longrightarrow} & A & \stackrel{+ \lambda(f; t, t')}{\longrightarrow} & A & \stackrel{- (s' - t') }{\longrightarrow} & A.
\end{array}
$$
Since the composite map in the lower row is $+ \lambda(f; s, s')$ by the former assertion, the latter assertion follows.\\
(4) This is easily seen.  $\;\;\;\;\;\;\;\; \Box$\\
\\
{\bf 1.2. Conjugate action on  group cochains.} Let $\Pi$ be a profinite group and let $M$ be an additive  discrete group on which $\Pi$ acts continuously from the left. Let $C^n(\Pi,M)$  $(n \geq 0)$ be the group of continuous $n$-cochains of $\Pi$ with coefficients in $M$ and let $d^{n+1} : C^n(\Pi,M) \rightarrow C^{n+1}(\Pi,M)$ be the coboundary homomorphisms defined by
{\small
$$\begin{array}{ll}
(d^{n+1}\alpha^n)(\gamma_1,\dots , \gamma_{n+1}) := &  \gamma_1\alpha^n(\gamma_2,\dots , \gamma_{n+1})\\
                                                                                  & + \displaystyle{ \sum_{i=1}^n (-1)^i \alpha^n(\gamma_1,\dots , \gamma_{i-1},\gamma_i\gamma_{i+1}, \gamma_{i+2},\dots, \gamma_{n+1}) }\\    
                                                                                  & + (-1)^{n+1}\alpha^n(\gamma_1,\dots , \gamma_n)
 \end{array}
 \leqno{(1.2.1)}
 $$
 }
 for $ \alpha^n \in C^n(\Pi,M)$ and $\gamma_1,\dots,\gamma_{n+1} \in \Pi$.  Let $Z^n(\Pi,M) := {\rm Ker}(d^{n+1})$ and $B^n(\Pi,M) := {\rm Im}(d^{n})$ be the subgroups of $C^n(\Pi,M)$ consisting of $n$-cocycles and $n$-coboundaries, respectively, and let $H^n(\Pi,M) := Z^n(\Pi,M)/B^n(\Pi,M)$, the $n$-th cohomology group of $\Pi$ with coefficients in $M$. By convention, we put $C^n(\Pi,M) =0$ for $n <0$. We sometimes write $d$ for $d^{n}$ simply if no misunderstanding is caused.

Note that $\Pi$ acts on $C^n(\Pi,M)$ from the left by 
$$ (\sigma.\alpha^n)(\gamma_1,\dots, \gamma_n) := \sigma \alpha^n( \sigma^{-1}\gamma_1\sigma, \dots , \sigma^{-1}\gamma_n \sigma) \leqno{(1.2.2)}$$
for $\alpha^n \in C^n(\Pi,M)$ and $\sigma, \gamma_1, \dots , \gamma_n \in \Pi$. By (1.2.1) and (1.2.2), we see that this action commutes with the coboundary homomorphisms:
$$ d^{n+1}(\sigma.\alpha^i) = \sigma.d^{n+1}(\alpha^i) \;\;\;\; (\alpha^i \in C^i(\Pi,M)). \leqno{(1.2.3)}$$

Now we shall describe the action of $\Pi$ on  $C^n(\Pi,M)$ in a concrete manner. For $\sigma$, $\sigma_1$, $\sigma_2 \in \Pi$,  $0 \leq i \leq j \leq n$ $(n \geq 1)$, and $1 \leq k \leq n-1$, we define  the maps $ s_i=s^n_i (\sigma) : \Pi^n \rightarrow \Pi^{n+1}$, 
$s_{i,j}=s^n_{i,j} (\sigma_1,\sigma_2) : \Pi^n \rightarrow \Pi^{n+2}$ and $t_k=t^n_k :\Pi^n \rightarrow \Pi^{n-1}$
by 
$$ \begin{array}{l}
s_i(g_1,g_2,\dots,g_n) :=  (g_1,\dots,g_i ,\sigma, \sigma^{-1} g_{i+1} \sigma, \dots, \sigma^{-1} g_{n} \sigma), \\
s_{i,j}(g_1,g_2,\dots,g_n) := (g_1,\dots,g_i ,\sigma_1, {\sigma_1}^{-1} g_{i+1} \sigma_1, \dots, {\sigma_1}^{-1} g_{j} \sigma_1,  \\
 \;\;\;\; \;\;\;\;\;\;\;\;\;\;\;\;\;\;\;\;\;\;\;\;\;\;\;\; \;\;\;\; \sigma_2, {(\sigma_1 \sigma_2)}^{-1} g_{j+1} \sigma_1 \sigma_2, \dots, {(\sigma_1 \sigma_2)}^{-1} g_{n} \sigma_1 \sigma_2),  \\
t_k(g_1,g_2,\dots,g_n)  :=(g_1,\dots,g_{k-1},g_k g_{k+1}, g_{k+2}, \dots, g_n)
\end{array}  \leqno{(1.2.4)}$$
for $(g_1,g_2,\dots,g_n) \in \Pi^n$. We note that $s_{j+1}^{n+1}(\sigma_2) \circ s_i^n(\sigma_1) = s_{i,j}^n(\sigma_1,\sigma_2)$.
We define the homomorphisms 
$$ \begin{array}{l} h^n_{\sigma} : C^{n+1}(\Pi, M) \longrightarrow C^{n}(\Pi, M),\\
H^n_{\sigma_1, \sigma_2} : C^{n+2}(\Pi, M) \longrightarrow C^{n}(\Pi, M)
\end{array}
$$
by
$$\begin{array}{l} {\displaystyle h^n_{\sigma}(\alpha^{n+1})  := \underset{0 \leq i \leq n}{\sum} (-1)^i (\alpha^{n+1} \circ s^n_i(\sigma)), } \\
{\displaystyle H^n_{\sigma_1,\sigma_2} (\alpha^{n+2}) := \underset{0 \leq i \leq j \leq n}{\sum} (-1)^{i+j} (\alpha^{n+2} \circ s^n_{i,j}(\sigma_1,\sigma_2))} 
\end{array} \leqno{(1.2.5)}  $$
for  $\alpha^{n+1} \in C^{n+1}(\Pi,M)$ and $\alpha^{n+2} \in C^{n+2}(\Pi,M)$. For example,  explicit forms of $h^n_{\sigma}(\alpha^{n+1}), H^n_{\sigma_1,\sigma_2} (\alpha^{n+2})$ for $n = 1, 2$ are given as follows:
$$ \begin{array}{l}
h^1_{\sigma}(\alpha^2)(g) = \alpha^2(\sigma, \sigma^{-1}g\sigma) - \alpha^2(g,\sigma).\\
h^2_{\sigma}(\alpha^3)(g_1,g_2) = \alpha^3(\sigma, \sigma^{-1}g_1\sigma, \sigma^{-1}g_2\sigma) - \alpha^3(g_1, \sigma, \sigma^{-1}g_2\sigma) +  \alpha^3(g_1,g_2,\sigma).\\
H^1_{\sigma_1, \sigma_2}(\alpha^3)(g) = \alpha^3(\sigma_1, \sigma_2, (\sigma_1\sigma_2)^{-1}g\sigma_1\sigma_2) -  \alpha^3(\sigma_1, \sigma_1^{-1}g\sigma_1, \sigma_2) + \alpha^3(g,\sigma_1,\sigma_2) \\
H^2_{\sigma_1, \sigma_2}(\alpha^4)(g_1,g_2) = \alpha^4(\sigma_1, \sigma_2, (\sigma_1\sigma_2)^{-1}g_1\sigma_1\sigma_2, (\sigma_1\sigma_2)^{-1}g_2\sigma_1\sigma_2) \\
\;\;\;\;\;\;\; \;\;\; - \alpha^4(\sigma_1,\sigma_1^{-1}g_1\sigma_1,\sigma_2,(\sigma_1\sigma_2)^{-1}g_2\sigma_1\sigma_2) 
 +  \alpha^4(\sigma_1,\sigma_1^{-1}g_1\sigma_1, \sigma_1^{-1}g_2\sigma_1,\sigma_2)  \\
\;\;\;\;\;+ \alpha^4(g_1,\sigma_1,\sigma_2,(\sigma_1\sigma_2)^{-1}g_2\sigma_1\sigma_2) 
  - \alpha^4(g_1,\sigma_1,\sigma_1^{-1}g_2\sigma_1,\sigma_2) + \alpha^4(g_1,g_2,\sigma_1,\sigma_2)
\end{array}$$
We call $ h^n_{\sigma}, H^n_{\sigma_1, \sigma_2}$ the {\em transgression} homomorphisms, which play roles similar to the transgression homomorphisms in [G].

The following Theorem 1.2.6 and Corollary 1.2.7 were shown in Appendices A and B of [Ki]. Here we give an elementary direct proof. See also Remark 1.2.8 below for the background of the proof. \\
\\
{\bf Theorem 1.2.6.}  {\em Notations being as above,  we have the following equalities.}
$$\sigma.\alpha^n - \alpha^n = h^n_{\sigma} (d^{n+1} (\alpha^n)) + d^n (h^{n-1}_{\sigma} (\alpha^n) ), $$
$$\sigma_1 . h^n_{\sigma_2} (\alpha^{n+1})-h^n_{\sigma_1 \sigma_2} (\alpha^{n+1}) + h^n_{\sigma_1}(\alpha^{n+1}) = H^n_{\sigma_1,\sigma_2}(d^{n+2} (\alpha^{n+1})) - d^n(H^{n-1}_{\sigma_1,\sigma_2} (\alpha^{n+1})).$$
{\em for $\alpha^n \in C^n(\Pi,M)$ and $\alpha^{n+1} \in C^{n+1}(\Pi,M)$ $(n\geq 1)$.}\\
\\
{\em Proof.}  
By (1.2.4), we can see
$$ \begin{array}{l}
s_i \circ t_k = 
\left\{ \begin{array}{l}
t_k \circ s_{i+1} \;\; (k \leq i)\\
t_{k+1} \circ s_i \;\; (i<k), 
\end{array} 
\right. \\
s_{i,j} \circ t_k = 
\left\{ \begin{array}{l}
t_k \circ s_{i+1,j+1}  \;\;(k \leq i)\\
t_{k+1} \circ s_{i,j+1} \;\;  (i < k \leq j)\\
t_{k+2} \circ s_{i,j} \;\; (j<k). 
\end{array} \right.
\end{array}
\leqno{(1.2.6.1)} 
$$
We note that $t_{i+1} \circ s_{i+1} = t_{i+1}\circ s_i$. 
By (1.2.1) and (1.2.5), we have, for any $(g_1,g_2, \dots, g_n) \in \Pi^n$,
$$ \begin{array}{ll} 
h^n_{\sigma}(d^{n+1} (\alpha^n))(g_1,\dots , g_n) & =(\sigma.\alpha^n)(g_1,\dots,g_n) \\
& + \underset{1 \leq i \leq n}{\sum} (-1)^i g_1(\alpha^n \circ s_{i-1}) (g_2, \dots,g_n)  \\
& + \underset{0 \leq i \leq n, 1 \leq k \leq n}{\sum} (-1)^{i+k} (\alpha^n \circ t_k \circ s_i) (g_1,\dots,g_n)  \\
& +(-1)^{n+n+1} \alpha^n(g_1,\dots,g_n) \\
& + \underset{0 \leq i \leq n-1} {\sum} (-1)^{i+n+1} (\alpha^n \circ s_i)(g_1,\dots,g_{n-1}),\\
d^n(h^{n-1}_{\sigma} (\alpha^n))(g_1,\dots , g_n)
& = \underset{0 \leq i \leq n-1}{\sum} (-1)^i g_1(\alpha^n \circ s_{i}) (g_2, \dots,g_n)  \\
& + \underset{0 \leq i \leq n-1, 1 \leq k \leq n-1}{\sum} (-1)^{i+k} (\alpha^n \circ s_i \circ t_k) (g_1,\dots,g_n)\\
& + \underset{0 \leq i \leq n-1} {\sum} (-1)^{i+n} (\alpha^n \circ s_i)(g_1,\dots,g_{n-1}),
\end{array}$$
and
$$ \begin{array}{l}  H^n_{\sigma_1,\sigma_2}(d^{n+2} (\alpha^{n+1}))(g_1,\dots,g_n)\\
=   (\sigma_1.h^n_{\sigma_2} (\alpha^{n+1}))(g_1,\dots,g_n)+ \underset{0 < i \leq j \leq n}{\sum} (-1)^{i+j} g_1(\alpha^{n+1} \circ s_{i-1,j-1}) (g_2,\dots,g_n)   \\
  -h^n_{\sigma_1 \sigma_2}(\alpha^{n+1})(g_1,\dots,g_n) + \underset{\underset{ {i \neq  j\ {\rm or}\ k \neq i+1} }{0 \leq i \leq j \leq n, 1 \leq k \leq n+1}}{\sum}(-1)^{i+j+k}(\alpha^{n+1} \circ t_k \circ s_{i,j}) (g_1,\dots,g_n)   \\
+   h^n_{\sigma_1}(\alpha^{n+1})(g_1,\dots,g_n) + \underset{0 \leq i \leq j \leq n-1}{\sum} (-1)^{i+j+n+2} (\alpha^{n+1} \circ s_{i,j})(g_1,\dots,g_{n-1}), \\
d^n(H^{n-1}_{\sigma_1,\sigma_2} (\alpha^{n+1}))(g_1,\dots,g_n) \\
=   \underset{0 \leq i \leq j \leq n-1}{\sum} (-1)^{i+j} g_1. (\alpha^{n+1} \circ s_{i,j}) (g_2,\dots,g_n)   \\
+  \underset{0 \leq i \leq j \leq n-1,1 \leq k \leq n-1}{\sum}(-1)^{i+j+k}(\alpha^{n+1} \circ s_{i,j} \circ t_k) (g_1,\dots,g_n)  \\
+  \underset{0 \leq i \leq j \leq n-1}{\sum} (-1)^{i+j+n} (\alpha^{n+1} \circ s_{i,j})(g_1,\dots,g_{n-1}).
\end{array}$$ 
Hence we have 
$$ \begin{array}{l} h^n_{\sigma}(d^{n+1} (\alpha^n))(g_1,\dots , g_n) + d^n(h^{n-1}_{\sigma} (\alpha^n))(g_1,\dots , g_n) \\
= (\sigma.\alpha^n)(g_1,\dots,g_n)  - \alpha^n(g_1,\dots,g_n) \\ 
 +\underset{0 \leq i \leq n, 1 \leq k \leq n}{\sum} (-1)^{i+k} (\alpha^n \circ t_k \circ s_i) (g_1,\dots,g_n) \\
 + \underset{0 \leq i \leq n-1, 1 \leq k \leq n-1}{\sum} (-1)^{i+k} (\alpha^n \circ s_i \circ t_k) (g_1,\dots,g_n),
\end{array}$$\\
and
$$ \begin{array}{l} H^n_{\sigma_1,\sigma_2}(d^{n+2} (\alpha^{n+1}))(g_1,\dots,g_n) -  d^n(H^{n-1}_{\sigma_1,\sigma_2}(\alpha^{n+1}))(g_1,\dots,g_n)\\
= \sigma_1 . h^n_{\sigma_2} (\alpha^{n+1})(g_1,\dots,g_n) -h^n_{\sigma_1 \sigma_2}(\alpha^{n+1})(g_1,\dots,g_n)  + h^n_{\sigma_1}(\alpha^{n+1})(g_1,\dots,g_n) \\
+\underset{\underset{ {i \neq  j\ {\rm or}\ k \neq i+1} }{0 \leq i \leq j \leq n, 1 \leq k \leq n+1}}{\sum}(-1)^{i+j+k}(\alpha^{n+1} \circ t_k \circ s_{i,j}) (g_1,\dots,g_n) \\
- \underset{0 \leq i \leq j \leq n-1,1 \leq k \leq n-1}{\sum}(-1)^{i+j+k}(\alpha^{n+1} \circ s_{i,j} \circ t_k) (g_1,\dots,g_n).
\end{array}$$
By (1.2.6.1), we obtain the required equalities. $\;\; \Box$
\\
\\
By (1.2.3), $\Pi$ acts on $Z^n(\Pi, M)$ from the left. This action is described by Theorem 1.2.6 as follows.\\
\\
{\bf Corollary 1.2.7.} {\em  Suppose $\alpha \in Z^{n}(\Pi, M)$ $(n \geq 1)$. For $\sigma \in \Pi$, we let}
$$\beta_{\sigma}  := h^{n-1}_{\sigma}(\alpha).$$
{\em Then we have}
$$  \sigma.\alpha = \alpha + d^{n}\beta_{\sigma}.$$
For $\sigma, \sigma' \in \Pi$, {\em we have}
$$ \beta_{\sigma \sigma'} = \beta_{\sigma} + \sigma.\beta_{\sigma'} \;\; \mbox{mod}\; B^{n-1}(\Pi,M),$$
{\em namely, the map $\Pi \ni \sigma \mapsto \beta_{\sigma} \; \mbox{mod} \; B^{n-1}(\Pi,M) \in C^{n-1}(\Pi,M)/B^{n-1}(\Pi,M)$ is a $1$-cocycle.}\\
\\
{\em Proof.} The both equalities are obtained immediately from Theorem 1.2.6, since $d^{n+1}(\alpha) = 0$ by  $\alpha \in Z^{n}(\Pi, M)$ $(n \geq 1)$. $\;\; \Box$
\\
\\
{\bf Remark 1.2.8} (Algebro-topological proof of Theorem 1.2.6). For $\sigma \in \Pi$, let $\sigma^{\bullet}$ denote the automorphism of the cochain complex $(C^{\bullet}(\Pi,M), d^{\bullet})$ defined by $\sigma^n(\alpha) := \sigma.\alpha$ for $\alpha \in C^n(\Pi,M)$. Then Theorem 1.2.6 asserts that the family of homomorphisms $\{ h^{n}_{\sigma} : C^{n+1}(\Pi,M) \rightarrow C^n(\Pi,M) \}$ gives a homotopy connecting  $\sigma^{\bullet}$ and ${\rm id}_{C^{\bullet}(\Pi,M)}$. Actually our explicit definition  (1.2.5) is obtained by making the following algebro-topological proof concrete:
 We may assume $\Pi$ is finite by the limit argument. Let ${\cal E}$ be the one-object category whose morphisms are the elements of $\Pi$.  We consider two functors ${\rm id}_{\cal E}, \widehat{\sigma} : {\cal E} \rightarrow {\cal E}$ defined by ${\rm id}_{\cal E}(g) := g, \widehat{\sigma}(g) := \sigma^{-1}g\sigma$ for each morphism $g \in \Pi$. Let ${\cal N} : {\rm Cat} \rightarrow {\rm Fct}(\Delta^{\rm op}, {\rm Set})$ denote the nerve functor, where ${\rm Cat}$ is the category of small categories and ${\rm Fct}(\Delta^{\rm op}, {\rm Set})$ is the category of simplicial sets. Define the natural transformation $\eta : \widehat{\sigma} \rightarrow {\rm id}_{\cal E}$ by $\eta(*) := \sigma$ ($*$ is the unique object of ${\cal E}$).  
Then $\eta$ induces a corresponding funcor $h_{\eta} : {\cal E} \times \underline{1} \rightarrow {\cal E}$, where $\underline{n}$ denotes the category 
defined by the set $\{ 0,1, \dots , n\}$ and its order. Then ${\cal N} h_{\eta} : {\cal N}{\cal E} \times {\cal N}\underline{1} \rightarrow {\cal N}{\cal E}$ is a homotopy connecting the two simplicial maps ${\cal N}\widehat{\sigma},  {\cal N}{\rm id}_{\cal E} : {\cal N}{\cal E} \rightarrow {\cal N}{\cal E}$. Let $C_n({\cal N}{\cal E}) = \mathbb{Z}[{\cal N}{\cal E}(\underline{n})]$ be the group of $n$-chains of the simplicial set ${\cal N}{\cal E}$.    By [My; Proposition 5.3] and [My; Proposition 6.2], ${\cal N} h_{\eta}$ induces a homotopy $\{ h_n^{\sigma} : C_n({\cal N}{\cal E})  \rightarrow C_{n+1}({\cal N}{\cal E})\}$ connecting  two chain maps $({\cal N} \hat{\sigma})_{\bullet}, ({\cal N}{\rm id}_{\cal E})_{\bullet} : C_{\bullet}({\cal N}{\cal E}) \rightarrow C_{\bullet}({\cal N}{\cal E})$. For the groups of $n$-cochains $C^n({\cal N}{\cal E},M)={\rm Hom}(C_n({\cal N}{\cal E}),M)$, the homotopy $\{ h_ n^{\sigma} \}$ induces the homotopy $\{ h^n_{\sigma} : C^{n+1}({\cal N}{\cal E},M) \rightarrow C^n({\cal N}{\cal E},M) \}$ connecting the two cochain maps $({\cal N}\widehat{\sigma})^{\bullet}, ({\cal N}{\rm id}_{\cal E})^{\bullet} :  C^{\bullet}({\cal N}{\cal E},M) \rightarrow C^{\bullet}({\cal N}{\cal E},M)$. Since  ${\cal N}{\cal E}(n)$ is  $\Pi^n$, we have the isomorphisms for $i\geq 0$
$$ C^n({\cal N}{\cal E}, M) \simeq {\rm Map}(\Pi^n, M) = C^n(\Pi,M).$$
Under the above isomorphisms, $({\cal N}\widehat{\sigma})^{\bullet}$ and $({\cal N}{\rm id}_{\cal E})^{\bullet}$ are identified with $\sigma^{\bullet}$ and ${\rm id}_{C^{\bullet}(\Pi,M)}$, respectively,  and hence $\{ h^n_{\sigma} \}$ gives a homotopy connecting  $\sigma^{\bullet}$ and ${\rm id}_{C^{\bullet}(\Pi,M)}$. $\;\; \Box$
\\
 \\
\begin{center}
{\bf 2. Classical theory}
\end{center}

In this section, we construct the arithmetic prequantization bundle and the arithmetic Chern-Simons $1$-cocycle for $\partial V_S := \sqcup_{i=1}^r {\rm Spec}(k_{\frak{p}_i})$, where $S = \{ \frak{p}_1, \dots , \frak{p}_r \}$ is a finite set of finite primes of an algebraic number field $k$ of finite degree over $\mathbb{Q}$, and the arithmetic Chern-Simons functional over a space of  Galois representations unramified outside $S$. These constructions correspond to the classical theory of topological Dijkgraaf-Witten TQFT. 

Throughout the rest of this paper, we fix a natural number $N >1$and let $\mu_N$ be the group of $N$-th roots of unity in the field $\mathbb{C}$ of complex numbers. We fix a primitive $N$-th root of unity $\zeta_N$ and the isomorphism $\mathbb{Z}/N\mathbb{Z} \simeq \mu_N;\; m \mapsto \zeta_N^m$.  The base  number field $k$ (in $\mathbb{C}$) is supposed to contain $\mu_N$.  Let $G$ be a finite group and let $c$ be a fixed 3-cocycle of $G$ with coefficients in $\mathbb{Z}/N\mathbb{Z}$, $c \in Z^3(G,\mathbb{Z}/N\mathbb{Z})$, where $G$ acts on $\mathbb{Z}/N\mathbb{Z}$ trivially. \\
\\
{\bf 2.1. Arithmetic  prequantization  bundles and arithmetic Chern-Simons 1-cocycles.}  We firstly develop a local theory at a finite prime. Let $\frak{p}$ be a finite prime of $k$ and let $k_{\frak{p}}$ be the $\frak{p}$-adic field. We let $\partial V_{\frak{p}} := {\rm Spec}(k_{\frak{p}})$, which play a role analogous to the boundary of a tubular neighborhood of  a knot
 (see the dictionary of the analogies in Introduction). Let $\Pi_{\frak{p}}$ denote the \'{e}tale fundamental group of $\partial V_{\frak{p}}$ with base point ${\rm Spec}(\overline{k}_{\frak{p}})$ ($\overline{k}_{\frak{p}}$ being an algebraic closure of $k_{\frak{p}}$), which is the absolute Galois group ${\rm Gal}(\overline{k}_{\frak{p}}/k_{\frak{p}})$.
 
Let ${\cal F}_{\frak{p}}$ be  the set of continuous homomorphisms of $\Pi_{\frak{p}}$ to $G$:
$$ {\cal F}_{\frak{p}} := {\rm Hom}_{\rm cont}(\Pi_{\frak{p}}, G).$$
 It is a finite set on which $G$ acts  from the right by
$$ {\cal F}_{\frak{p}} \times G \rightarrow {\cal F}_{\frak{p}}; \;\; (\rho_{\frak{p}}, g) \mapsto \rho_{\frak{p}}.g := g^{-1} \rho_{\frak{p}} g.  \leqno{(2.1.1)}$$
 Let ${\cal M}_{\frak{p}}$ denote the quotient space by this action:
$$ {\cal M}_{\frak{p}} := {\cal F}_{\frak{p}}/G.$$
Let ${\rm Map}({\cal F}_{\frak{p}}, \mathbb{Z}/N\mathbb{Z})$ denote the additive group consisting of maps from ${\cal F}_{\frak{p}}$ to $\mathbb{Z}/N\mathbb{Z}$, on which $G$ acts from the left by
$$ (g.\psi_{\frak{p}})(\rho_{\frak{p}}) := \psi_{\frak{p}}(\rho_{\frak{p}}.g)      \leqno{(2.1.2)}$$
for $g \in G, \psi_{\frak{p}} \in {\rm Map}({\cal F}_{\frak{p}}, \mathbb{Z}/N\mathbb{Z})$ and $\rho_{\frak{p}} \in {\cal F}_{\frak{p}}$. For $\rho_{\frak{p}} \in {\cal F}_{\frak{p}}$ and $\alpha \in C^n(G,\mathbb{Z}/N\mathbb{Z})$, we denote by $\alpha \circ \rho_{\frak{p}}$ the $n$-cochain of $\Pi_{\frak{p}}$ with coefficients in $\mathbb{Z}/N\mathbb{Z}$ defined by
$$ (\alpha \circ \rho_{\frak{p}})(\gamma_1, \dots , \gamma_n) := \alpha( \rho_{\frak{p}}(\gamma_1),\dots , \rho_{\frak{p}}(\gamma_n)). 
$$
By (1.2.2) and (2.1.1), we have
$$ (g.\alpha)\circ \rho_{\frak{p}} = \alpha \circ (\rho_{\frak{p}}.g)    \leqno{(2.1.3)}$$
for $g \in G, \alpha \in C^n(G, \mathbb{Z}/N\mathbb{Z})$ and $\rho_{\frak{p}} \in {\cal F}_{\frak{p}}$.

Firstly, we shall construct an arithmetic analog for $\partial V_{\frak{p}} := {\rm Spec}(k_{\frak{p}})$ of the prequantization bundle, using the given $3$-cocycle $c \in Z^3(G,\mathbb{Z}/N\mathbb{Z})$. The key idea for this is due to  Kim ([Ki]), who uses the conjugate $G$-action on $c$ and the $2$nd Galois cohomology group (Brauer group) of the local field $k_{\frak{p}}$.

Let $\rho_{\frak{p}} \in {\cal F}_{\frak{p}}$ and so $c \circ \rho_{\frak{p}} \in Z^3(\Pi_{\frak{p}},\mathbb{Z}/N\mathbb{Z})$. Let $d$ denote the coboundary homomorphism $C^2(\Pi_{\frak{p}},\mathbb{Z}/N\mathbb{Z}) \rightarrow C^3(\Pi_{\frak{p}},\mathbb{Z}/N\mathbb{Z})$. We define ${\cal L}_{\frak{p}}(\rho_{\frak{p}})$ 
by the quotient set
$$ {\cal L}_{\frak{p}}(\rho_{\frak{p}}) := d^{-1}(c\circ \rho_{\frak{p}})/B^2(\Pi_{\frak{p}},\mathbb{Z}/N\mathbb{Z}). \leqno{(2.1.4)}$$
Here we note that $d^{-1}(c \circ \rho_{\frak{p}})$ is non-empty, because the cohomological dimension of $\Pi_{\frak{p}}$ is 2 ([NSW; Theorem 7.1.8], [S1; Chapitre II, 5.3, Proposition 15]) and so $H^3(\Pi_{\frak{p}}, \mathbb{Z}/N\mathbb{Z}) = 0$. Thus $d^{-1}(c \circ \rho_{\frak{p}})$ is a $Z^2(\Pi_{\frak{p}},\mathbb{Z}/N\mathbb{Z})$-torsor in the obvious manner and so ${\cal L}_{\frak{p}}(\rho_{\frak{p}})$ is an $H^2(\Pi_{\frak{p}},\mathbb{Z}/N\mathbb{Z})$-torsor by (2.1.4) and Lemma 1.1.4 (4). Since $k_{\frak{p}}$ contains $\mu_N$ and so $H^2(\Pi_{\frak{p}},\mathbb{Z}/N\mathbb{Z}) = H^2(k_{\frak{p}},\mu_N)$, the theory of Brauer groups (cf. [S2; Chapitre XII]) tells us that there is the canonical isomorphism
$$ {\rm inv}_{\frak{p}} : H^2(\Pi_{\frak{p}},\mathbb{Z}/N\mathbb{Z}) \stackrel{\sim}{\longrightarrow} \mathbb{Z}/N\mathbb{Z}$$
and hence ${\cal L}_{\frak{p}}(\rho_{\frak{p}})$ is a $\mathbb{Z}/N\mathbb{Z}$-torsor via ${\rm inv}_{\frak{p}}$. 

Let ${\cal L}_{\frak{p}}$ be the disjoint union of  ${\cal L}_{\frak{p}}(\rho_{\frak{p}})$ over all $\rho_{\frak{p}} \in {\cal F}_{\frak{p}}$:
$$ {\cal L}_{\frak{p}} := \bigsqcup_{\rho_{\frak{p}} \in {\cal F}_{\frak{p}}} {\cal L}_{\frak{p}}(\rho_{\frak{p}})  $$
and consider the projection
$$ \varpi_{\frak{p}} : {\cal L}_{\frak{p}} \longrightarrow {\cal F}_{\frak{p}}; \;\; \alpha_{\frak{p}} \mapsto \rho_{\frak{p}} \; \mbox{if}\; \alpha_{\frak{p}} \in {\cal L}_{\frak{p}}(\rho_{\frak{p}}).   $$
Since each fiber $\varpi_{\frak{p}}^{-1}(\rho_{\frak{p}}) = {\cal L}_{\frak{p}}(\rho_{\frak{p}})$ is a $\mathbb{Z}/N\mathbb{Z}$-torsor, we may regard ${\cal L}_{\frak{p}}$  as a principal $\mathbb{Z}/N\mathbb{Z}$-bundle over ${\cal F}_{\frak{p}}$.

Let $g \in G$. Using the transgression map $h^2_g$, we define $h_g \in C^2(G, \mathbb{Z}/N\mathbb{Z})/B^2(G, \mathbb{Z}/N\mathbb{Z})$ by
$$  h_g := h^2_g(c) \; \mbox{mod} \; B^2(G, \mathbb{Z}/N\mathbb{Z}),$$
where $h^2_g(c)$ is the 2-cochain defined as in (1.2.5), namely,
$$ h^2_g(c)(g_1,g_2) := c(g, g^{-1}g_1 g, g^{-1}g_2 g) - c(g_1, g, g^{-1}g_2 g) +  c(g_1,g_2, g),$$
where $g_1, g_2 \in G$. By Corollary 1.2.7, we have
$$ g.c = c + d h_g \leqno{(2.1.5)}$$
and 
$$ h_{gg'} = h_g + g.h_{g'} \leqno{(2.1.6)}$$
for $g, g' \in G$.   
By (2.1.3), (2.1.4) and (2.1.5), we have
$$  
d(\alpha + h_g\circ \rho_{\frak{p}})  = c\circ \rho_{\frak{p}} + (g.c - c) \circ \rho_{\frak{p}} = (g.c) \circ \rho_{\frak{p}} = c \circ (\rho_{\frak{p}}.g) 
$$
for $\alpha_{\frak{p}} \in {\cal L}_{\frak{p}}(\rho_{\frak{p}})$ and so we have the isomorphism of $\mathbb{Z}/N\mathbb{Z}$-torsors
$$ f_{\frak{p}}(g,\rho_{\frak{p}}) : {\cal L}_{\frak{p}}(\rho_{\frak{p}}) \stackrel{\sim}{\longrightarrow} {\cal L}_{\frak{p}}(\rho_{\frak{p}}.g); \;\; \alpha_{\frak{p}} \mapsto \alpha_{\frak{p}} + h_g \circ \rho_{\frak{p}}. \leqno{(2.1.7)}$$
By (2.1.3) and (2.1.6), we have
$$ \begin{array}{ll} \alpha_{\frak{p}} + h_{gg'}\circ \rho_{\frak{p}} & = \alpha_{\frak{p}} + (h_g + g.h_{g'})\circ \rho_{\frak{p}} \\
                                                                                                     & = \alpha_{\frak{p}} + h_g \circ \rho_{\frak{p}} + h_{g'} \circ (\rho_{\frak{p}}.g) 
                                                                                                     \end{array}$$
for $g, g' \in G$. It means that  $G$ acts on ${\cal L}_{\frak{p}}$ from the right by
$$ {\cal L}_{\frak{p}} \times G \rightarrow {\cal L}_{\frak{p}}; \;\; \alpha_{\frak{p}} \mapsto \alpha_{\frak{p}}.g := f(g,\rho_{\frak{p}})(\alpha_{\frak{p}}). \leqno{(2.1.8)}$$
By (2.1.7), (2.1.8) and the way of the $\mathbb{Z}/N\mathbb{Z}$-action on ${\cal L}_{\frak{p}}$,  we have the following commutative diagram 
$$  \begin{array}{rcl}
 {\cal L}_{\frak{p}} & \stackrel{.g}{\longrightarrow} & {\cal L}_{\frak{p}} \curvearrowleft \mathbb{Z}/N\mathbb{Z} \\
\varpi_{\frak{p}} \downarrow & & \downarrow \varpi_{\frak{p}} \\
{\cal F}_{\frak{p}} & \stackrel{.g}{\longrightarrow} & {\cal F}_{\frak{p}},\\
\end{array} 
$$
namely, 
$$(\alpha_{\frak{p}}.m).g = (\alpha_{\frak{p}}.g).m, \;\; \varpi_{\frak{p}}(\alpha_{\frak{p}}.g) = \varpi_{\frak{p}}(\alpha_{\frak{p}}).g 
 \leqno{(2.1.9)} 
$$
for $\alpha_{\frak{p}} \in {\cal F}_{\frak{p}}$, $m \in \mathbb{Z}/N\mathbb{Z}, g \in G$. So ${\cal L}_{\frak{p}}$ is a $G$-equivariant principal $\mathbb{Z}/N\mathbb{Z}$-bundle over ${\cal F}_{\frak{p}}$.  
Taking the quotient by the action of $G$, we have the principal $\mathbb{Z}/N\mathbb{Z}$-bundle $\overline{\varpi}_{\frak{p}} : \overline{\cal L}_{\frak{p}} \rightarrow {\cal M}_{\frak{p}}$. We call  $\varpi_{\frak{p}} : {\cal L}_{\frak{p}} \rightarrow {\cal F}_{\frak{p}}$ or $\overline{\varpi}_{\frak{p}} : \overline{\cal L}_{\frak{p}} \rightarrow {\cal M}_{\frak{p}}$ the {\it arithmetic prequantization $\mathbb{Z}/N\mathbb{Z}$-bundle} for $\partial V_\frak{p} := {\rm Spec}(k_{\frak{p}})$. 

Let us choose a section $x_{\frak{p}} \in \Gamma({\cal F}_{\frak{p}}, {\cal L}_{\frak{p}})$, namely,  the map 
$$x_{\frak{p}} : {\cal F}_{\frak{p}} \longrightarrow {\cal L}_{\frak{p}}  \; \; \mbox{such that} \;\varpi_{\frak{p}} \circ x_{\frak{p}} = {\rm id}_{{\cal F}_{\frak{p}}}.   $$ 
This means that we fix a ``coordinate" on ${\cal L}_{\frak{p}}$. In fact, by the trivialization at $x_{\frak{p}}(\rho_{\frak{p}})$ in (1.1.3), we may identify each fiber ${\cal L}_{\frak{p}}(\rho_{\frak{p}})$ over $\rho_{\frak{p}}$ with $\mathbb{Z}/N\mathbb{Z}$:
$$ \varphi_{x_{\frak{p}}(\rho_{\frak{p}})} :  {\cal L}_{\frak{p}}(\rho_{\frak{p}}) \stackrel{\sim}{\longrightarrow} \mathbb{Z}/N\mathbb{Z}; \; \alpha_{\frak{p}} \mapsto \alpha_{\frak{p}} - x_{\frak{p}}(\rho_{\frak{p}}). $$
For $g \in G$ and $\rho_{\frak{p}} \in {\cal F}_{\frak{p}}$, we let 
$$ \lambda_{\frak{p}}^{x_{\frak{p}}}(g, \rho_{\frak{p}}) := f_{\frak{p}}(g, \rho_{\frak{p}})(x_{\frak{p}}(\rho_{\frak{p}})) - x_{\frak{p}}(\rho_{\frak{p}}.g) = x_{\frak{p}}(\rho_{\frak{p}}).g - x_{\frak{p}}(\rho_{\frak{p}}.g)  \leqno{(2.1.10)}$$
 so that we have the following commutative 
diagram by Lemma 1.1.4 (3):
$$ \begin{array}{rcl}
{\cal L}_{\frak{p}}(\rho_{\frak{p}})   &  \stackrel{f_{\frak{p}}(g,\rho_{\frak{p}})}{\longrightarrow} &  {\cal L}_{\frak{p}}(\rho_{\frak{p}}.g)  \\
\varphi_{x_{\frak{p}}(\rho_{\frak{p}})} \downarrow & & \downarrow \varphi_{x_{\frak{p}}(\rho_{\frak{p}}.g)} \\
\mathbb{Z}/N\mathbb{Z} & \stackrel{+ \lambda_{\frak{p}}^{x_{\frak{p}}}(g,\rho_{\frak{p}})}{\longrightarrow}  &  \mathbb{Z}/N\mathbb{Z},
\end{array}
$$
namely, for $\alpha_{\frak{p}} \in {\cal L}_{\frak{p}}(\rho_{\frak{p}})$, we have
$$ \alpha_{\frak{p}}.g - x_{\frak{p}}(\rho_{\frak{p}}.g) = (\alpha_{\frak{p}} - x_{\frak{p}}(\rho_{\frak{p}})) + \lambda_{\frak{p}}^{x_{\frak{p}}}(g,\rho_{\frak{p}}). \leqno{(2.1.11)}$$
We define the map $\lambda_{\frak{p}}^{x_{\frak{p}}} : G \rightarrow {\rm Map}({\cal F}_{\frak{p}},\mathbb{Z}/N\mathbb{Z})$ by
$$ \lambda_{\frak{p}}^{x_{\frak{p}}}(g)(\rho_{\frak{p}}) := \lambda_{\frak{p}}^{x_{\frak{p}}}(g,\rho_{\frak{p}}) \leqno{(2.1.12)}$$
for $g \in G$ and $\rho_{\frak{p}} \in {\cal F}_{\frak{p}}$. 
\\
\\
{\bf Theorem 2.1.13.}  {\em For $g, g' \in G$, we have }
$$ \lambda_{\frak{p}}^{x_{\frak{p}}}(gg') = \lambda_{\frak{p}}^{x_{\frak{p}}}(g) + (g.\lambda_{\frak{p}}^{x_{\frak{p}}})(g').$$
{\em Namely, the map $\lambda_{\frak{p}}^{x_{\frak{p}}}$ is a $1$-cocycle}:
$$ \lambda_{\frak{p}}^{x_{\frak{p}}} \in Z^1(G, {\rm Map}({\cal F}_{\frak{p}},\mathbb{Z}/N\mathbb{Z})).$$
{\em Proof.}  For $g, g' \in G$ and $\rho_{\frak{p}} \in {\cal F}_{\frak{p}}$, we have 
$$ \begin{array}{ll}
\lambda_{\frak{p}}^{x_{\frak{p}}}(gg', \rho_{\frak{p}}) & = f_{\frak{p}}(gg',\rho_{\frak{p}})(x_{\frak{p}}(\rho_{\frak{p}})) - x_{\frak{p}}(\rho_{\frak{p}}(gg')) \; \; \mbox{by}\; (2.1.10) \\
                                                            & = (x_{\frak{p}}(\rho_{\frak{p}}) + h_{gg'}\circ \rho_{\frak{p}}) - x_{\frak{p}}(\rho_{\frak{p}}.(gg')) \;\; \mbox{by} \; (2.1.7)\\
                                                             & = (x_{\frak{p}}(\rho_{\frak{p}}) + h_g\circ \rho_{\frak{p}} + h_{g'} \circ (\rho_{\frak{p}}.g)) -  x_{\frak{p}}(\rho_{\frak{p}}.(gg')) \;\;  \mbox{by}\; (2.1.3), (2.1.6).\\
                                                            
\end{array} $$
By Lemma 1.1.4 (1), we have
$$ \begin{array}{l} (x_{\frak{p}}(\rho_{\frak{p}}) + h_g\circ \rho_{\frak{p}} + h_{g'} \circ (\rho_{\frak{p}}.g)) -  x_{\frak{p}}(\rho_{\frak{p}}.(gg'))\\
 \;\;\;\;\;\;\;\;\;\; = \{ (x_{\frak{p}}(\rho_{\frak{p}}) + h_g \circ \rho_{\frak{p}})  - x_{\frak{p}}(\rho_{\frak{p}}.g) \}  + \{ (x_{\frak{p}}(\rho_{\frak{p}}.g) + h_{g'}\circ (\rho_{\frak{p}}.g)) - x_{\frak{p}}(\rho_{\frak{p}}.(gg'))    \}. 
\end{array}
$$
Here we see by (2.1.7), (2.1.10) that
$$ \begin{array}{l}
(x_{\frak{p}}(\rho_{\frak{p}}) + h_g \circ \rho_{\frak{p}})  - x_{\frak{p}}(\rho_{\frak{p}}.g) = \lambda_{\frak{p}}^{x_{\frak{p}}}(g, \rho_{\frak{p}}),\\
( x_{\frak{p}}(\rho_{\frak{p}}.g) + h_{g'}\circ (\rho_{\frak{p}}.g) ) - x_{\frak{p}}(\rho_{\frak{p}}.(gg')) = \lambda_{\frak{p}}^{x_{\frak{p}}}(g', \rho_{\frak{p}}.g).
\end{array}
$$
Combining these, we have 
$$ \lambda_{\frak{p}}^{x_{\frak{p}}}(gg',\rho_{\frak{p}}) =  \lambda_{\frak{p}}^{x_{\frak{p}}}(g, \rho_{\frak{p}}) + \lambda_{\frak{p}}^{x_{\frak{p}}}(g', \rho_{\frak{p}}.g)$$
for any $\rho_{\frak{p}} \in {\cal F}_{\frak{p}}$. By (2.1.2) and (2.1.12), we obtain the assertion. $\;\; \Box$ \\
\\
We call $\lambda_{\frak{p}}^{x_{\frak{p}}}$ the {\it Chern-Simons 1-cocycle} for $\partial V_\frak{p}$ with respect to the section $x_{\frak{p}}$.\\
\\
For a section $x_{\frak{p}} \in \Gamma({\cal F}_{\frak{p}}, {\cal L}_{\frak{p}})$, we define ${\cal L}_{\frak{p}}^{x_{\frak{p}}}$ by the product (trivial)  principal $\mathbb{Z}/N\mathbb{Z}$-bundle over ${\cal F}_{\frak{p}}$:
$$ {\cal L}_{\frak{p}}^{x_{\frak{p}}} := {\cal F}_{\frak{p}} \times \mathbb{Z}/N\mathbb{Z},  $$
on which $G$ acts from the right by 
$$   {\cal L}_{\frak{p}}^{x_{\frak{p}}} \times G \rightarrow {\cal L}_{\frak{p}}^{x_{\frak{p}}}; \;\; ((\rho_{\frak{p}},m),g) \mapsto (\rho_{\frak{p}}.g, m+\lambda_{\frak{p}}^{x_{\frak{p}}}(g,\rho_{\frak{p}})), \leqno{(2.1.14)}$$
and so the projection
$$ \varpi_{\frak{p}}^{x_{\frak{p}}} : {\cal L}_{\frak{p}}^{x_{\frak{p}}}  \longrightarrow {\cal F}_{\frak{p}}$$
is $G$-equivariant.\\
\\
{\bf Proposition 2.1.15.} {\em We have the following isomorphism of $G$-equivariant principal  $\mathbb{Z}/N\mathbb{Z}$-bundles}
$$ \Phi_{\frak{p}}^{x_{\frak{p}}} : {\cal L}_{\frak{p}} \stackrel{\sim}{\longrightarrow} {\cal L}_{\frak{p}}^{x_{\frak{p}}}; \;\; \alpha_{\frak{p}} \mapsto (\varpi_{\frak{p}}(\alpha_{\frak{p}}), \alpha_{\frak{p}} - x_{\frak{p}}(\varpi_{\frak{p}}(\alpha_{\frak{p}}))).$$
{\it In particular, the isomorphism class of ${\cal L}_{\frak{p}}^{x_{\frak{p}}}$ is independent of the choice of a section $x_{\frak{p}}$. In other words, for another section $x'_{\frak{p}}  \in \Gamma({\cal F}_{\frak{p}}, {\cal L}_{\frak{p}})$, we have ${\cal L}_{\frak{p}}^{x'_{\frak{p}}} \simeq {\cal L}_{\frak{p}}^{x_{\frak{p}}}$ as $G$-equivariant principal $\mathbb{Z}/N\mathbb{Z}$-bundles.}\\
\\
{\em Proof.} (i) It is easy to see that $\varpi_{\frak{p}}^{x_{\frak{p}}}\circ \Phi_{\frak{p}}^{x_{\frak{p}}} = \varpi_{\frak{p}}$.\\ 
(ii) For $\alpha_{\frak{p}} \in {\cal L}_{\frak{p}}$ and $m \in \mathbb{Z}/N\mathbb{Z}$,  we have
$$ \begin{array}{ll}
\Phi_{\frak{p}}^{x_{\frak{p}}}(\alpha_{\frak{p}}.m)& = (\varpi_{\frak{p}}(\alpha_{\frak{p}}.m), \alpha_{\frak{p}}.m - x_{\frak{p}}(\varpi_{\frak{p}}(\alpha_{\frak{p}}.m)))\\
& = (\varpi_{\frak{p}}(\alpha_{\frak{p}}), \alpha_{\frak{p}}.m - x_{\frak{p}}(\varpi_{\frak{p}}(\alpha_{\frak{p}})))\\
& = (\varpi_{\frak{p}}(\alpha_{\frak{p}}), (\alpha_{\frak{p}} - x_{\frak{p}}(\varpi_{\frak{p}}(\alpha_{\frak{p}}))) + m) \; \mbox{by Lemma 1.1.4 (1)}\\
& = \Phi_{\frak{p}}^{x_{\frak{p}}}(\alpha_{\frak{p}}).m.
\end{array}
$$
(iii) $\Phi_{\frak{p}}^{x_{\frak{p}}}$ has the inverse defined by $(\Phi_{\frak{p}}^{x_{\frak{p}}})^{-1}((\rho_{\frak{p}},m)) := x_{\frak{p}}(\rho_{\frak{p}}).m$ for $(\rho_{\frak{p}},m) \in {\cal F}_{\frak{p}} \times \mathbb{Z}/N\mathbb{Z}$.\\
By (i), (ii), (iii), $\Phi_{\frak{p}}^{x_{\frak{p}}}$  is an isomorphism of principal $\mathbb{Z}/N\mathbb{Z}$-bundles. So it suffices to show that $\Phi^{x_{\frak{p}}}$ is $G$-equivariant. It follows from that 
$$ \begin{array}{ll} \Phi_{\frak{p}}^{x_{\frak{p}}}(\alpha_{\frak{p}}.g) & = (\varpi_{\frak{p}}(\alpha_{\frak{p}}.g), \alpha_{\frak{p}}.g - x_{\frak{p}}(\varpi_{\frak{p}}(\alpha_{\frak{p}}.g)))\\
                                                                                         & = (\varpi_{\frak{p}}(\alpha_{\frak{p}}).g, (\alpha_{\frak{p}} - x_{\frak{p}}(\varpi_{\frak{p}}(\alpha_{\frak{p}}))) + \lambda_{\frak{p}}^{x_{\frak{p}}}(g, \varpi_{\frak{p}}(\alpha_{\frak{p}}))) \;\; \mbox{by (2.1.9), (2.1.11)}\\
                                                                                         & = \Phi_{\frak{p}}^{x_{\frak{p}}}(\alpha).g \;\; \;\;  \mbox{by (2.1.14)}   \;\; \Box                    \end{array}
$$
\\
Taking the quotient of $\varpi_{\frak{p}}^{x_{\frak{p}}} : {\cal L}_{\frak{p}}^{x_{\frak{p}}} \rightarrow {\cal F}_{\frak{p}}$ by the action of $G$, we have the principal $\mathbb{Z}/N\mathbb{Z}$-bundle $\overline{\varpi}_{\frak{p}}^{x_{\frak{p}}} : \overline{\cal L}_{\frak{p}}^{x_{\frak{p}}} \rightarrow {\cal M}_{\frak{p}}$. We call  $\varpi_{\frak{p}}^{x_{\frak{p}}} : {\cal L}_{\frak{p}}^{x_{\frak{p}}} \rightarrow {\cal F}_{\frak{p}}$ or $\overline{\varpi}_{\frak{p}}^{x_{\frak{p}}} : \overline{\cal L}_{\frak{p}}^{x_{\frak{p}}} \rightarrow {\cal M}_{\frak{p}}$ the {\it arithmetic prequantization principal $\mathbb{Z}/N\mathbb{Z}$-bundle} for $\partial V_\frak{p}$ with respect to the section $x_{\frak{p}}$.\\
\\
For $x_{\frak{p}}, x'_{\frak{p}} \in \Gamma({\cal F}_{\frak{p}}, {\cal L}_{\frak{p}})$, we define the map $\delta_{\frak{p}}^{x_{\frak{p}},x'_{\frak{p}}} : {\cal F}_{\frak{p}} \rightarrow \mathbb{Z}/N\mathbb{Z}$ by
$$ \delta_{\frak{p}}^{x_{\frak{p}},x'_{\frak{p}}}(\rho_{\frak{p}}) := x_{\frak{p}}(\rho_{\frak{p}}) - x'_{\frak{p}}(\rho_{\frak{p}}) \leqno{(2.1.16)}$$
for $\rho_{\frak{p}} \in {\cal F}_{\frak{p}}$. \\
\\
{\bf Lemma 2.1.17} {\it For $x_{\frak{p}}, x'_{\frak{p}}, x''_{\frak{p}} \in \Gamma({\cal F}_{\frak{p}}, {\cal L}_{\frak{p}})$, we have}
$$ \delta_{\frak{p}}^{x_{\frak{p}},x_{\frak{p}}} = 0, \;  \delta_{\frak{p}}^{x'_{\frak{p}},x_{\frak{p}}} = - \delta_{\frak{p}}^{x_{\frak{p}},x'_{\frak{p}}}, \; \delta_{\frak{p}}^{x_{\frak{p}},x'_{\frak{p}}} + \delta_{\frak{p}}^{x'_{\frak{p}},x''_{\frak{p}}} = \delta_{\frak{p}}^{x_{\frak{p}},x''_{\frak{p}}}.$$
\\
{\it Proof.} These equalities follow from Lemma 1.1.4 (1). $\;\; \Box$\\
\\
 The following proposition tells us how $\lambda_{\frak{p}}^{x_{\frak{p}}}$ is changed when we change the section $x_{\frak{p}}$.\\
 \\
{\bf Proposition 2.1.18.} {\it   For $x_{\frak{p}}, x'_{\frak{p}} \in \Gamma({\cal F}_{\frak{p}}, {\cal L}_{\frak{p}})$, we have }
$$ \lambda_{\frak{p}}^{x'_{\frak{p}}}(g) - \lambda_{\frak{p}}^{x_{\frak{p}}}(g) = g.\delta_{\frak{p}}^{x_{\frak{p}},x'_{\frak{p}}} - \delta_{\frak{p}}^{x_{\frak{p}},x'_{\frak{p}}}$$
{\em for any $g \in G$. So the cohomology class $[\lambda_{\frak{p}}^{x_{\frak{p}}}] \in H^1(G, {\rm Map}({\cal F}_{\frak{p}},\mathbb{Z}/N\mathbb{Z}))$ is independent of the choice of a section $x_{\frak{p}}$.}\\
\\
{\it Proof.}  By (2.1.10) and  Lemma 1.1.4 (1), (2), we have
$$ \begin{array}{ll}
\lambda_{\frak{p}}^{x'_{\frak{p}}}(g,\rho_\frak{p}) - \lambda_{\frak{p}}^{x_{\frak{p}}}(g,\rho_{\frak{p}}) & = 
(f_{\frak{p}}(g, \rho_{\frak{p}})(x'_{\frak{p}}(\rho_{\frak{p}})) - x'_{\frak{p}}(\rho_{\frak{p}}.g) ) - ( f_{\frak{p}}(g, \rho_{\frak{p}})(x_{\frak{p}}(\rho_{\frak{p}})) - x_{\frak{p}}(\rho_{\frak{p}}.g)   )\\
& =  ( x_{\frak{p}}(\rho_{\frak{p}}.g) -  x'_{\frak{p}}(\rho_{\frak{p}}.g) )  + (  f_{\frak{p}}(g, \rho_{\frak{p}})(x'_{\frak{p}}(\rho_{\frak{p}}))  -  f_{\frak{p}}(g, \rho_{\frak{p}})(x_{\frak{p}}(\rho_{\frak{p}}))  ) \\
& = ( x_{\frak{p}}(\rho_{\frak{p}}.g) -  x'_{\frak{p}}(\rho_{\frak{p}}.g) ) + (x'_{\frak{p}}(\rho_{\frak{p}}) - x_{\frak{p}}(\rho_{\frak{p}}))  \\
& =  (g.\delta_{\frak{p}}^{x_{\frak{p}},x'_{\frak{p}}})(\rho_{\frak{p}}) - \delta_{\frak{p}}^{x_{\frak{p}},x'_{\frak{p}}}(\rho_{\frak{p}}) \;\; \mbox{by (2.1.2)}
\end{array}
$$
for any $g \in G$ and $\rho_{\frak{p}} \in {\cal F}_{\frak{p}}$, hence the assertion. $\;\;\;\;\;\;\;\; \Box$\\
\\
By Proposition 2.1.18, we denote the cohomology class $[\lambda_{\frak{p}}^{x_{\frak{p}}}]$ by $[\lambda_\frak{p}]$, which we call  the {\it arithmetic Chern-Simons 1st cohomology class}  for $\partial V_\frak{p}$. As a corollary of Proposition 2.1.18, we can make the latter statement of Proposition 2.1.15 more precise 
as follows.\\
\\
{\bf Corollary 2.1.19.} (1) {\em For $x_{\frak{p}}, x'_{\frak{p}} \in \Gamma({\cal F}_{\frak{p}}, {\cal L}_{\frak{p}})$,  we have the following isomorphism of $G$-equivariant principal $\mathbb{Z}/N\mathbb{Z}$-bundles over ${\cal F}_{\frak{p}}$:}
$$ \Phi_{\frak{p}}^{x_{\frak{p}}, x'_{\frak{p}}} : {\cal L}_{\frak{p}}^{x_{\frak{p}}} \stackrel{\sim}{\longrightarrow} {\cal L}_{\frak{p}}^{x'_{\frak{p}}}; \;\; (\rho_{\frak{p}}, m)  \mapsto (\rho_{\frak{p}}, m + \delta_{\frak{p}}^{x_{\frak{p}}, x'_{\frak{p}}}(\rho_{\frak{p}})),$$
{\it where $\delta_{\frak{p}}^{x_{\frak{p}}, x'_{\frak{p}}} : {\cal F}_{\frak{p}} \rightarrow \mathbb{Z}/N\mathbb{Z}$ is the map defined in  (2.1.16). }\\
(2) {\em For $x_{\frak{p}}, x'_{\frak{p}}, x''_{\frak{p}} \in \Gamma({\cal F}_{\frak{p}}, {\cal L}_{\frak{p}})$, we have}
$$  \left\{ \begin{array}{l} \Phi_{\frak{p}}^{x_{\frak{p}}, x'_{\frak{p}} } \circ \Phi_{\frak{p}}^{x_{\frak{p}}} = \Phi_{\frak{p}}^{x'_{\frak{p}}}, \\
\Phi_{\frak{p}}^{x_{\frak{p}}, x_{\frak{p}} } = {\rm id}_{{\cal L}_{\frak{p}}^{x_{\frak{p}}}}, \;  \Phi_{\frak{p}}^{x'_{\frak{p}}, x_{\frak{p}} } =  (\Phi_{\frak{p}}^{x_{\frak{p}}, x'_{\frak{p}} })^{-1}, \;  \Phi_{\frak{p}}^{x'_{\frak{p}}, x''_{\frak{p}} } \circ  \Phi_{\frak{p}}^{x_{\frak{p}}, x'_{\frak{p}} } =  \Phi_{\frak{p}}^{x_{\frak{p}}, x''_{\frak{p}} }
\end{array} \right.$$
\\
{\it Proof.} (1) We easily see that $\Phi_{\frak{p}}^{x_{\frak{p}}, x'_{\frak{p}} }$ is isomorphism of principal $\mathbb{Z}/N\mathbb{Z}$-bundles and so it suffices to show that $\Phi_{\frak{p}}^{x_{\frak{p}}, x'_{\frak{p}} }$ is $G$-equivariant. This follows from
$$\begin{array}{ll}
\Phi_{\frak{p}}^{x_{\frak{p}}, x'_{\frak{p}} }((\rho_{\frak{p}},m).g) & = \Phi_{\frak{p}}^{x_{\frak{p}}, x'_{\frak{p}} }((\rho_{\frak{p}}.g, m + \lambda_{\frak{p}}^{x_{\frak{p}}}(g, \rho_{\frak{p}})))  \;\; \mbox{by (2.1.14)}\\
                                                       & = (\rho_{\frak{p}}.g, m + \lambda_{\frak{p}}^{x_{\frak{p}}}(g,\rho_{\frak{p}}) + \delta_{\frak{p}}^{x_{\frak{p}}, x'_{\frak{p}}}(\rho_{\frak{p}}.g))\\
                                                       & = (\rho_{\frak{p}}.g, m + \delta_{\frak{p}}^{x_{\frak{p}}, x'_{\frak{p}}}(\rho_{\frak{p}}) + \lambda_{\frak{p}}^{x'_{\frak{p}}}(g, \rho_{\frak{p}})) \;\; \mbox{by Proposition 2.1.18} \\
                                                       & = \Phi_{\frak{p}}^{x_{\frak{p}},x'_{\frak{p}}}(\rho_{\frak{p}},m).g.  
\end{array}
$$ 
(2) The first equality follows from the definitions of $\Phi_{\frak{p}}^{x_{\frak{p}}}, \Phi_{\frak{p}}^{x_{\frak{p}}, x'_{\frak{p}}}$. The latter equalities follow from Lemma 2.1.17. $\;\; \Box$\\

Let $F$ be a field containing $\mu_N$. Let $L_{\frak{p}}$ be the $F$-line bundle over ${\cal F}_{\frak{p}}$ associated to the principal $\mathbb{Z}/N\mathbb{Z}$-bundle ${\cal L}_{\frak{p}}$ and the homomorphism $\mathbb{Z}/N\mathbb{Z} \hookrightarrow F^{\times}$; $m \mapsto \zeta_N^m$, namely,
$$ \begin{array}{ll} L_{\frak{p}} & := \displaystyle{ {\cal L}_{\frak{p}} \times_{\mathbb{Z}/N\mathbb{Z}} F }\\
 & := ({\cal L}_{\frak{p}} \times F)/(\alpha_{\frak{p}}, z)\sim (\alpha_{\frak{p}}.m, \zeta_N^{-m}z) \;\; (\alpha_{\frak{p}} \in {\cal L}_{\frak{p}}, m \in \mathbb{Z}/N\mathbb{Z}, z \in F),
\end{array}   \leqno{(2.1.20)}
$$ 
on which $G$ acts from the right by
$$ L_{\frak{p}} \times G \rightarrow L_{\frak{p}}; \;\; ([(\alpha_{\frak{p}}, z)],g) \mapsto  [(\alpha_{\frak{p}}.g,z)]. \leqno{(2.1.21)}$$
The projection 
$$ \varpi_{\frak{p},F} : L_{\frak{p}} \longrightarrow {\cal F}_{\frak{p}}; \; [(\alpha_{\frak{p}}, z)] \mapsto \varpi_{\frak{p}}(\alpha_{\frak{p}}) $$
is a $G$-equivariant $F$-line bundle. We denote the fiber $\varpi_{\frak{p},F}^{-1}(\rho_{\frak{p}})$ over $\rho_{\frak{p}}$ by $L_{\frak{p}}(\rho_{\frak{p}})$:
$$ L_{\frak{p}}(\rho_{\frak{p}}) := \{ [(\alpha_{\frak{p}},z)] \in L_{\frak{p}} \, | \, \varpi_{\frak{p}}(\alpha_{\frak{p}}) = \rho_{\frak{p}}, z \in F \}  \leqno{(2.1.22)}$$
We have a non-canonical bijection by fixing an $\alpha_{\frak{p}} \in {\cal L}_{\frak{p}}(\rho_{\frak{p}})$:
$$ L_{\frak{p}}(\rho_{\frak{p}}) \stackrel{\sim}{\longrightarrow} F; \; [(\alpha_{\frak{p}}, z)] \mapsto z. $$
Taking the quotient by the action of $G$, we obtain the $F$-line bundles $\overline{\varpi}_{\frak{p}, F} : \overline{L}_{\frak{p}} \rightarrow {\cal M}_{\frak{p}}$. We call  $\varpi_{\frak{p}, F} : L_{\frak{p}} \rightarrow {\cal F}_{\frak{p}}$ or $\overline{\varpi}_{\frak{p},F} : \overline{L}_{\frak{p}} \rightarrow {\cal M}_{\frak{p}}$ the {\it arithmetic prequantization $F$-line bundle} for $\partial V_\frak{p}$.

Let $L_{\frak{p}}^{x_{\frak{p}}}$ be the product $F$-line bundle over ${\cal F}_{\frak{p}}$:
$$ L_{\frak{p}}^{x_{\frak{p}}} := {\cal F}_{\frak{p}} \times F, $$ 
on which $G$ acts from the right by
$$ L_{\frak{p}}^{x_{\frak{p}}}  \times G \rightarrow L_{\frak{p}}^{x_{\frak{p}}}; \;\; ((\rho_{\frak{p}},z),g) \mapsto (\rho_{\frak{p}}.g, z \zeta_N^{\lambda_{\frak{p}}^{x_{\frak{p}}}(g,\rho_{\frak{p}})}), \leqno{(2.1.23)} $$
and the projection 
$$ \varpi_{\frak{p},F}^{x_{\frak{p}}} : L_{\frak{p}}^{x_{\frak{p}}} \longrightarrow {\cal F}_{\frak{p}}$$
is $G$-equivariant. 
Then we have the following Proposition  similar to Proposition 2.1.15 and Corollary 2.1.19.\\
\\
{\bf Proposition 2.1.24.} {\em We have the following isomorphism of $G$-equivariant $F$-line bundles over ${\cal F}_{\frak{p}}$}
$$\Phi_{\frak{p}, F}^{x_{\frak{p}}} : L_{\frak{p}}  \stackrel{\sim}{\longrightarrow} L_{\frak{p}}^{x_{\frak{p}}}; [(\alpha_{\frak{p}},z)] \mapsto (\varpi_{\frak{p}}(\alpha_{\frak{p}}), z\zeta_N^{\alpha_{\frak{p}} - x_{\frak{p}}(\varpi_{\frak{p}}(\alpha_{\frak{p}}))}).   $$
{\em For another section $x'_{\frak{p}}$, we have the following isomorphism of $G$-equivariant $F$-line bundles over ${\cal F}_{\frak{p}}$}
$$ \Phi_{\frak{p}, F}^{x_{\frak{p}}, x'_{\frak{p}}} : L_{\frak{p}}^{x_{\frak{p}}} \stackrel{\sim}{\longrightarrow} L_{\frak{p}}^{x'_{\frak{p}}} : (\rho_{\frak{p}},z) \mapsto (\rho_{\frak{p}}, z\zeta_N^{\delta_{\frak{p}}^{x_{\frak{p}}, x'_{\frak{p}}}(\rho_{\frak{p}})}),$$
{\em where $\delta_{\frak{p}}^{x_{\frak{p}}, x'_{\frak{p}}} : {\cal F}_{\frak{p}} \rightarrow \mathbb{Z}/N\mathbb{Z}$ is the map in (2.1.16), and we have the equalities}
$$ \left\{ \begin{array}{l} \Phi_{\frak{p},F}^{x_{\frak{p}}, x'_{\frak{p}}} \circ \Phi_{\frak{p},F}^{x_{\frak{p}}} = \Phi_{\frak{p},F}^{x'_{\frak{p}}} \\ 
\Phi_{\frak{p},F}^{x_{\frak{p}}, x_{\frak{p}} } = {\rm id}_{L_{\frak{p},F}^{x_{\frak{p}}}}, \;  \Phi_{\frak{p},F}^{x'_{\frak{p}}, x_{\frak{p}} } =  (\Phi_{\frak{p},F}^{x_{\frak{p}}, x'_{\frak{p}} })^{-1}, \;  \Phi_{\frak{p},F}^{x'_{\frak{p}}, x''_{\frak{p}} } \circ  \Phi_{\frak{p},F}^{x_{\frak{p}}, x'_{\frak{p}} } =  \Phi_{\frak{p},F}^{x_{\frak{p}}, x''_{\frak{p}} }
\end{array} \right. $$
{\em for} $x_{\frak{p}}, x'_{\frak{p}}, x''_{\frak{p}} \in \Gamma({\cal F}_{\frak{p}}, {\cal L}_{\frak{p}})$.\\
\\
{\it Proof.} (i) It is easy to see that $\varpi_{\frak{p},F}^{x_{\frak{p}}} \circ \Phi_{\frak{p},F}^{x_{\frak{p}}} = \varpi_{\frak{p},F}$.\\
(ii) For $\rho_{\frak{p}} \in {\cal F}_{\frak{p}}$, we let
$$  L_{\frak{p}}^{x_{\frak{p}}}(\rho_{\frak{p}})  := (\varpi_{\frak{p},F}^{x_{\frak{p}}})^{-1}(\rho_{\frak{p}})\\
                                                                  = \{ (\rho_{\frak{p}}, z) \; | \; z \in F\} \simeq F.
                                                                 $$
So $\Phi_{\frak{p},F}^{x_{\frak{p}}}$ restricted to a fiber over $\rho_{\frak{p}}$ 
$$ \Phi_{\frak{p},F}^{x_{\frak{p}}}|_{L_{\frak{p}}(\rho_{\frak{p}})} : L_{\frak{p}}(\rho_{\frak{p}}) \longrightarrow L_{\frak{p}}^{x_{\frak{p}}}(\rho_{\frak{p}}); \;\; [(\alpha_{\frak{p}},z)] \mapsto (\rho_{\frak{p}},z\zeta_N^{\alpha_{\frak{p}} - x_{\frak{p}}(\rho_{\frak{p}})})$$
is $F$-linear.\\
(iv) For $g \in G$, we have
$$ \begin{array}{ll} \Phi_{\frak{p}, F}^{x_{\frak{p}}}([(\alpha_{\frak{p}},z)].g) & = \Phi_{\frak{p},F}^{x_{\frak{p}}}([(\alpha_{\frak{p}}.g,z)]) \;\; \mbox{by (2.1.21)}\\
 & = (\varpi_{\frak{p}}(\alpha_{\frak{p}}.g), z\zeta_N^{\alpha_{\frak{p}}.g - x_{\frak{p}}(\varpi_{\frak{p}}(\alpha_{\frak{p}}.g))} )\\
          & = (\varpi_{\frak{p}}(\alpha_{\frak{p}}).g, z\zeta^{(\alpha_{\frak{p}} - x_{\frak{p}}(\rho_{\frak{p}})) + \lambda_{\frak{p}}^{x_{\frak{p}}}(g, \rho_{\frak{p}})}) \;\; \mbox{by (2.1.11)}\\
          & = \Phi_{\frak{p},F}^{x_{\frak{p}}}([\alpha_{\frak{p}},z)]).g \;\; \mbox{by (2.1.23).}

\end{array}$$
Hence $\Phi_{\frak{p},F}^{x_{\frak{p}}}$ is the isomorphism of $G$-equivariant $F$-line bundles over ${\cal F}_{\frak{p}}$. 

The proofs of the latter parts are similar to those of Corollary 2.1.19 (1), (2). $\;\; \Box$\\
\\
Taking the quotient of $\varpi_{\frak{p}, F}^{x_{\frak{p}}} : L_{\frak{p}}^{x_{\frak{p}}} \rightarrow {\cal F}_{\frak{p}}$ by the action of $G$, we have the $F$-line bundle $\overline{\varpi}_{\frak{p}, F}^{x_{\frak{p}}} : \overline{L}_{\frak{p}}^{x_{\frak{p}}} \rightarrow {\cal M}_{\frak{p}}$. We call $\varpi_{\frak{p}, F}^{x_{\frak{p}}} : L_{\frak{p}}^{x_{\frak{p}}} \rightarrow {\cal F}_{\frak{p}}$ or $\overline{\varpi}_{\frak{p},F}^{x_{\frak{p}}} : \overline{L}_{\frak{p}}^{x_{\frak{p}}} \rightarrow {\cal M}_{\frak{p}}$ the {\it arithmetic prequantization $F$-line bundle} for $\partial V_\frak{p}$ with respect to the section $x_{\frak{p}}$. \\

Let $S = \{ \frak{p}_1, \dots , \frak{p}_r \}$ be a finite set of finite primes of $k$ and let $\partial V_S := \partial V_{\frak{p}_1} \sqcup \cdots \sqcup \partial V_{\frak{p}_r}$. Let ${\cal F}_S$ be the direct product of ${\cal F}_{\frak{p}_i}$'s:
$$ {\cal F}_S := {\cal F}_{\frak{p}_1} \times \cdots \times {\cal F}_{\frak{p}_r}.$$
It is a finite set on which $G$ acts diagonally from the right, namely, 
$$ {\cal F}_S \times G \rightarrow {\cal F}_S;  \;\;  (\rho_S, g) \mapsto \rho_S. g := (\rho_{\frak{p}_1}. g, \dots , \rho_{\frak{p}_r}.g) \leqno{(2.1.25)} $$
for  $\rho_S = (\rho_{\frak{p}_1}, \dots , \rho_{\frak{p}_r}) \in {\cal F}_S$
and let ${\cal M}_{S}$ denote the quotient space by this action
$$ {\cal M}_{S} := {\cal F}_{S}/G.$$
Let ${\rm Map}({\cal F}_S, \mathbb{Z}/N\mathbb{Z})$ be the additive group of maps from ${\cal F}_S$ to $\mathbb{Z}/N\mathbb{Z}$, on which $G$ acts from the left by
$$(g.\psi_S)(\rho_S) := \psi_S(\rho_S.g)  \leqno{(2.1.26)} $$
for $\psi_S \in {\rm Map}({\cal F}_S,\mathbb{Z}/N\mathbb{Z}), g \in G$ and $\rho_S \in {\cal F}_S$.

For $\rho_S = (\rho_{\frak{p}_1},\dots , \rho_{\frak{p}_r}) \in {\cal F}_S$, let ${\cal L}_S(\rho_S)$ be the quotient space of the product ${\cal L}_{\frak{p}_1}(\rho_{\frak{p}_1}) \times \cdots \times {\cal L}_{\frak{p}_r}(\rho_{\frak{p}_r})$: 
$$ {\cal L}_S(\rho_S) := ({\cal L}_{\frak{p}_1}(\rho_{\frak{p}_1})  \times \cdots \times {\cal L}_{\frak{p}_r}(\rho_{\frak{p}_r}))/\sim, \leqno{(2.1.27)}$$
where the equivalence relation $\sim$ is  defined by
$$ (\alpha_{\frak{p}_1}, \dots , \alpha_{\frak{p}_r}) \sim (\alpha'_{\frak{p}_1}, \dots , \alpha'_{\frak{p}_r})  \Longleftrightarrow \sum_{i=1}^r (\alpha_{\frak{p}_i} - \alpha'_{\frak{p}_i}) = 0. \leqno{(2.1.28)}$$
We see easily that ${\cal L}_S(\rho_S)$ is equipped with the simply transitive action of  $\mathbb{Z}/N\mathbb{Z}$ defined  by
$$ \begin{array}{l} {\cal L}_S(\rho_S) \times \mathbb{Z}/N\mathbb{Z} \longrightarrow {\cal L}_S(\rho_S); \\
([\alpha_S], m) \mapsto [\alpha_S].m := [(\alpha_{\frak{p}_1}.m, \dots , \alpha_{\frak{p}_r})]=\cdots = [(\alpha_{\frak{p}_1}, \dots , \alpha_{\frak{p}_r}.m)] 
\end{array}$$
for $\alpha_S = (\alpha_{\frak{p}_1}, \dots , \alpha_{\frak{p}_r})$ and hence ${\cal L}_S(\rho_S)$ is a $\mathbb{Z}/N\mathbb{Z}$-torsor.

Let ${\cal L}_S$ be the disjoint union of  ${\cal L}_{\frak{p}}(\rho_{S})$ for $\rho_{S} \in {\cal F}_{S}$:
$$ {\cal L}_S := \bigsqcup_{\rho_S \in {\cal F}_S} {\cal L}_{S}(\rho_S),   \leqno{(2.1.29)}$$
on which $G$ acts diagonally from the right by
$$ {\cal L}_S \times G \longrightarrow {\cal L}_S; \; ([(\alpha_{\frak{p}_1}, \dots , \alpha_{\frak{p}_r})],g) \mapsto [(\alpha_{\frak{p}_1}.g, \dots, \alpha_{\frak{p}_r}.g)].  \leqno{(2.1.30)}$$
Consider the projection
$$ \varpi_S : {\cal L}_S \longrightarrow {\cal F}_S; \; [\alpha_S] = [(\alpha_{\frak{p}_i})] \mapsto (\varpi_{\frak{p}_i}(\alpha_{\frak{p}_i})),   $$
which is $G$-equivariant. Since each fiber $\varpi_{\frak{p}}^{-1}(\rho_S) = {\cal L}_{S}(\rho_S)$ is a $\mathbb{Z}/N\mathbb{Z}$-torsor, we may regard $ \varpi_S : {\cal L}_S \longrightarrow {\cal F}_S$ as a $G$-equivariant principal $\mathbb{Z}/N\mathbb{Z}$-bundle. Taking the quotient by the action of $G$, we have the principal $\mathbb{Z}/N\mathbb{Z}$-bundle $\overline{\varpi}_{S} : \overline{\cal L}_{S} \rightarrow {\cal M}_{S}$. We call  $\varpi_{S} : {\cal L}_{S} \rightarrow {\cal F}_{S}$ or $\overline{\varpi}_{S} : \overline{\cal L}_{S} \rightarrow {\cal M}_{S}$ the {\it arithmetic prequantization $\mathbb{Z}/N\mathbb{Z}$-bundle} for $\partial V_S = {\rm Spec}(k_{\frak{p}_1}) \sqcup \cdots \sqcup {\rm Spec}(k_{\frak{p}_r})$. 

Let $x_S$ be a section of $\varpi_S$, $ x_S \in \Gamma({\cal F}_S, {\cal L}_S)$. By (2.1.27) and (2.1.29), it is written as $x_S = [(x_{\frak{p}_1}, \dots , x_{\frak{p}_r})]$, where $x_{\frak{p}_i} \in \Gamma({\cal F}_{\frak{p}_i}, {\cal L}_{\frak{p}_i})$ for $1\leq i \leq r$. For $g \in G$ and $\rho_S = (\rho_{\frak{p}_i}) \in {\cal F}_S$, we set  
$$  \lambda_S^{x_S}(g, \rho_S) := \lambda_{\frak{p}_1}^{x_{\frak{p}_1}}(g,\rho_{\frak{p}_1}) + \cdots + \lambda_{\frak{p}_r}^{x_{\frak{p}_r}}(g, \rho_{\frak{p}_r}) \leqno{(2.1.31)}$$
and define the map $\lambda_S^{x_S} : G \rightarrow {\rm Map}({\cal F}_S, \mathbb{Z}/N\mathbb{Z})$ by
$$ \lambda_S^{x_S}(g)(\rho_S) := \lambda_S^{x_S}(g,\rho_S)   \leqno{(2.1.32)}$$
for $g \in G$ and $\rho_S \in {\cal F}_S$. 
\\  
\\
{\bf Lemma 2.1.33.} (1) {\it Let $x'_{\frak{p}_i} \in \Gamma({\cal F}_{\frak{p}_i},{\cal L}_{\frak{p}_i})$ be another section for $1\leq i \leq r$ such that $[(x'_{\frak{p}_1}, \dots , x'_{\frak{p}_r})] = x_S$. Then we have}
$$ \sum_{i=1}^r  \lambda_{\frak{p}_i}^{x_{\frak{p}_i}}(g,\rho_{\frak{p}_i}) =  \sum_{i=1}^r  \lambda_{\frak{p}_i}^{x'_{\frak{p}_i}}(g,\rho_{\frak{p}_i})$$
{\it for $g \in G$ and $\rho_{\frak{p}_i} \in {\cal F}_{\frak{p}_i}$. So $\lambda_S^{x_S}(g, \rho_S)$ is independent of the choic of $x_{\frak{p}_i}$'s such that $x_S = [(x_{\frak{p}_1}, \dots , x_{\frak{p}_r})]$.}\\
(2) {\em The map $\lambda_S^{x_S}$ is a 1-cocycle}: 
$$ \lambda_S^{x_S} \in Z^1(G, {\rm Map}({\cal F}_S, \mathbb{Z}/N\mathbb{Z})).$$
\\
{\it Proof.} (1) Since $(x_{\frak{p}_1}(\rho_{\frak{p}_1}), \dots , x_{\frak{p}_r}(\rho_{\frak{p}_r})) \sim (x'_{\frak{p}_1}(\rho_{\frak{p}_1}), \dots , x'_{\frak{p}_r}(\rho_{\frak{p}_r}))$, by (2.1.28), we have
$$ \sum_{i=1}^r (x_{\frak{p}_i}(\rho_{\frak{p}_i}) - x'_{\frak{p}_i}(\rho_{\frak{p}_i})) = 0$$
for any $\rho_{\frak{p}_i} \in {\cal F}_{\frak{p}_i}$. Therefore we have
$$ \begin{array}{ll}   \displaystyle{ \sum_{i=1}^r  \lambda_{\frak{p}_i}^{x_{\frak{p}_i}}(g,\rho_{\frak{p}_i}) } & =  \displaystyle{  \sum_{i=1}^r  ( f_{\frak{p}_i}(g, \rho_{\frak{p}_i})(x_{\frak{p}_i}(\rho_{\frak{p}_i})) - x_{\frak{p}_i}(\rho_{\frak{p}_i}.g) ) } \;\;\mbox{by (2.1.10)} \\
&= \displaystyle{  \sum_{i=1}^r  (  (f_{\frak{p}_i}(g, \rho_{\frak{p}_i})(x_{\frak{p}_i}(\rho_{\frak{p}_i})) - f_{\frak{p}_i}(g, \rho_{\frak{p}_i})(x'_{\frak{p}_i}(\rho_{\frak{p}_i}))) }\\
&  \;\;  +  \displaystyle{ \sum_{i=1}^r   (f_{\frak{p}_i}(g, \rho_{\frak{p}_i})(x'_{\frak{p}_i}(\rho_{\frak{p}_i}))    -   x'_{\frak{p}_i}(\rho_{\frak{p}_i}.g))  }\\
&   \;\;  + \displaystyle{  \sum_{i=1}^r  (x'_{\frak{p}_i}(\rho_{\frak{p}_i}.g) - x_{\frak{p}_i}(\rho_{\frak{p}_i}.g))  } \;\;\;\;\;\;\;\; \mbox{by Lemma 1.1.4 (1)}\\
& = \displaystyle{ \sum_{i=1}^r   (f_{\frak{p}_i}(g, \rho_{\frak{p}_i})(x'_{\frak{p}_i}(\rho_{\frak{p}_i}))    -   x'_{\frak{p}_i}(\rho_{\frak{p}_i}.g))  }\;\;\;\; \mbox{by Lemma 1.1.4 (2)} \\
& = \displaystyle{ \sum_{i=1}^r \lambda_{\frak{p}_i}^{x'_{\frak{p}_i}}(g,\rho_{\frak{p}_i})}
\end{array}
$$
for $g \in G$ and $\rho_{\frak{p}_i} \in {\cal F}_{\frak{p}_i}$.\\
(2) By Theorem 2.1.13, (2.1.26), (2.1.31) and (2.1.32), we have
$$ \begin{array}{ll} \lambda_S^{x_S}(gg', \rho_S)  & = \displaystyle{ \sum_{i=1}^r \lambda_{\frak{p}_i}^{x_{\frak{p}_i}}(gg', \rho_{\frak{p}_i}) }\\
                                                                              & = \displaystyle{ \sum_{i=1}^r  \lambda_{\frak{p}_i}^{x_{\frak{p}_i}}(g, \rho_{\frak{p}_i}) + \sum_{i=1}^r  \lambda_{\frak{p}_i}^{x_{\frak{p}_i}}(g',  \rho_{\frak{p}_i}.g) }\\
                                                                              & = \lambda_S^{x_S}(g,\rho_S) + \lambda_S^{x_S}(g',\rho_S.g)\\
                                                                             & = (\lambda_S^{x_S}(g) + (g.\lambda_S^{x_S})(g'))(\rho_S)
                                                                             \end{array}$$
for $g \in G$ and $\rho_S = (\rho_{\frak{p}_i}) \in {\cal F}_S$.  Thus we obtain
$$ \lambda_S^{x_S}(gg') = \lambda_S^{x_S}(g) + (g.\lambda_S^{x_S})(g') \;\; \Box$$
\\
We call $\lambda_S^{x_S}$ the {\em arithmetic Chern-Simons 1-cocycle} for $\partial V_S$ with respect to $x_S$. \\
\\
{\bf Proposition 2.1.34.}  {\it   Let $x'_S = [(x'_{\frak{p}_1}, \dots , x'_{\frak{p}_r})] \in \Gamma({\cal F}_S,{\cal L}_S)$ be another section of $\varpi_S$. We define the map
$\delta_S^{x_S, x'_S} : {\cal F}_S \rightarrow \mathbb{Z}/N\mathbb{Z}$ by}
$$ \delta_S^{x_S, x'_S}(\rho_S) := \sum_{i=1}^r \delta_{\frak{p}_i}^{x_{\frak{p}_i},x'_{\frak{p}_i}}(\rho_{\frak{p}_i})$$
{\em for $\rho_S = (\rho_{\frak{p}_i}) \in {\cal F}_S$, where $\delta_{\frak{p}_i}^{x_{\frak{p}_i},x'_{\frak{p}_i}}$ is the map defined in (2.1.16). Then we have}
$$ \lambda_S^{x'_S}(g) - \lambda_S^{x_S}(g) = g.\delta_S^{x_S, x'_S} - \delta_S^{x_S, x'_S}$$
{\it for $g \in G$. So the cohomology class $[\lambda_S^{x_S}] \in H^1(G, {\rm Map}({\cal F}_S, \mathbb{Z}/N\mathbb{Z}))$ is independent of the choice of $x_S$.}\\
\\
{\it Proof.} First, note that $\delta_S^{x_S,x'_S}$ is proved to be  independent of the choices of $x_{\frak{p}_i}$'s  in the similar manner to 
the proof of Lemma 2.1.33 (1). By the definition of  $\delta_S^{x_S,x'_S}$, the formula follows from Proposition 2.1.18  by taking the sum over $\frak{p}_i \in S$. $\;\; \Box$\\
\\
We denote the cohomology class $[\lambda_S^{x_S}]$ by $[\lambda_S]$, which we call the {\em arithmetic Chern-Simons 1st cohomology class} for $\partial V_S$.

Let ${\cal L}_S^{x_S}$ be the product principal $\mathbb{Z}/N\mathbb{Z}$-bundle over ${\cal F}_S$:
$$ {\cal L}_S^{x_S} := {\cal F}_S \times \mathbb{Z}/N\mathbb{Z}, $$
on which $G$ acts from the right by
$$  {\cal L}_S^{x_S} \times G \rightarrow {\cal L}_S^{x_S}; \;\; ((\rho_S, m), g) \mapsto (\rho_S.g, m+\lambda_S^{x_S}(g,\rho_S)). $$
\\
{\bf Proposition 2.1.35.} {\em We have the following isomorphism of $G$-equivariant principal $\mathbb{Z}/N\mathbb{Z}$-bundles over ${\cal F}_S$}:
$$ \Phi_S^{x_S} :   {\cal L}_S \stackrel{\sim}{\longrightarrow} {\cal L}_S^{x_S}; \; [\alpha_S] = [(\alpha_{\frak{p}_1}, \dots , \alpha_{\frak{p}_r})] \mapsto 
(\varpi_S([\alpha_S]), \sum_{i=1}^r (\alpha_{\frak{p}_i} - x_{\frak{p}_i}(\varpi_{\frak{p}_i}(\alpha_{\frak{p}_i})) ).
$$
{\em For another section $x'_{S}$, we have the following isomorphism of $G$-equivariant $F$-line bundles over ${\cal F}_{S}$}
$$ \Phi_S^{x_{S}, x'_{S}} : {\cal L}_{S}^{x_{S}} \stackrel{\sim}{\longrightarrow} {\cal L}_{S}^{x'_{S}} : (\rho_{S},m) \mapsto (\rho_{S}, m + \delta_{S}^{x_{S}, x'_{S}}(\rho_{S})),$$
{\em where $\delta_{S}^{x_{S}, x'_{S}} : {\cal F}_{S} \rightarrow \mathbb{Z}/N\mathbb{Z}$ is the map in Proposition 2.1.34. For $x_{S}, x'_{S}, x''_{S} \in \Gamma({\cal F}_{S}, {\cal L}_{S})$ we have the equalities}
$$ \left\{ \begin{array}{l}
\Phi_S^{x_S,x'_S} \circ \Phi_S^{x_S} = \Phi_S^{x'_S},\\
\Phi_{S}^{x_{S}, x_{S} } = {\rm id}_{{\cal L}_S^{x_{S}}}, \;  \Phi_S^{x'_{S}, x_{S} } =  (\Phi_{S}^{x_{S}, x'_{S} })^{-1}, \;  \Phi_{S}^{x'_{S}, x''_{S} } \circ  \Phi_{S}^{x_{S}, x'_{S} } =  \Phi_{S}^{x_{S}, x''_{S} }.
\end{array} \right. $$
\\
{\it Proof.}   First, suppose $[(\alpha_{\frak{p}_1}, \dots , \alpha_{\frak{p}_r})] = [(\alpha'_{\frak{p}_1}, \dots , \alpha'_{\frak{p}_r})]$. Then $\varpi_{\frak{p}_i}(\alpha_{\frak{p}_i}) = \varpi_{\frak{p}_i}(\alpha'_{\frak{p}_i})$ and $\sum_{i=1}^r (\alpha'_{\frak{p}_i} - \alpha_{\frak{p}_i}) = 0$ by (2.1.28). So we have
$$ \begin{array}{ll} \displaystyle{ \sum_{i=1}^r (\alpha'_{\frak{p}_i} - x_{\frak{p}_i}( \varpi_{\frak{p}_i}(\alpha'_{\frak{p}_i}))) }& = \displaystyle{ \sum_{i=1}^r \left( (\alpha'_{\frak{p}_i} - \alpha_{\frak{p}_i}) + ( \alpha_{\frak{p}_i} -x_{\frak{p}_i}( \varpi_{\frak{p}_i}(\alpha'_{\frak{p}_i})))  \right) }\\
                                                       & = \displaystyle{ \sum_{i=1}^r (\alpha_{\frak{p}_i} - x_{\frak{p}_i}( \varpi_{\frak{p}_i}(\alpha_{\frak{p}_i}))).}\\
                                                       \end{array} $$
 The proofs of the assertions go well in the similar manner to those of Proposition 2.1.15 and Corollary 2.1.19, by taking  the sum over $\frak{p}_i \in S$. $\;\; \Box$
\\
 \\
Taking the quotient by the action of $G$, we obtain the principal $\mathbb{Z}/N\mathbb{Z}$-bundle $\overline{\varpi}_S^{x_S} : \overline{\cal L}_S^{x_S} \rightarrow {\cal M}_S$.  We call $\varpi_S^{x_S} : {\cal L}_S^{x_S} \rightarrow {\cal F}_S$ or $\overline{\varpi}_S^{x_S} : \overline{\cal L}_S^{x_S} \rightarrow {\cal M}_S$ the {\em arithmetic prequantization principal $\mathbb{Z}/N\mathbb{Z}$-bundle}  for $\partial V_S$ with respect to $x_S$. \\

Let $L_S$ be the $F$-line bundle associated to the principal $\mathbb{Z}/N\mathbb{Z}$-bundle ${\cal L}_S$ over ${\cal F}_S$ and the homomorphism 
$\mathbb{Z}/N\mathbb{Z} \rightarrow F^{\times}; m \mapsto \zeta_N^m$:
$$ \begin{array}{ll} L_{S} & := \displaystyle{ {\cal L}_{S} \times_{\mathbb{Z}/N\mathbb{Z}} F }\\
 & := ({\cal L}_{S} \times F)/([\alpha_S], z)\sim ([\alpha_S].m, \zeta_N^{-m}z) \;\; ([\alpha_S] \in {\cal L}_{S}, m \in \mathbb{Z}/N\mathbb{Z}, z \in F),
\end{array}  \leqno{(2.1.36)}$$ 
on which $G$ acts from the right by
$$ L_{S} \times G \longrightarrow L_{S}; \;\; ([([\alpha_S], z)],g) \mapsto  [([\alpha_S].g,z)]. \leqno{(2.1.37)}$$
The projection 
$$ \varpi_{S,F} : L_{S} \longrightarrow {\cal F}_{S}; \; [([\alpha_S], z)] \mapsto \varpi_{S}([\alpha_S]) $$
is a $G$-equivariant $F$-line bundle. We denote the fiber $\varpi_{S,F}^{-1}(\rho_{S})$ over $\rho_{S}$ by $L_{S}(\rho_{S})$, which is non-canonically bijective to $F$ by fixing $[\alpha_S] \in {\cal L}_S(\rho_S)$:
$$ L_{S}(\rho_{S}) := \{ [([\alpha_S],z)] \in L_{S} \, | \, \varpi_{S}([\alpha_S]) = \rho_{S} \} \stackrel{\sim}{\rightarrow} F; \; [([\alpha_S],z)] \mapsto z. \leqno{(2.1.38)}$$
Taking the quotient by the action of $G$, we obtain the $F$-line bundle $\overline{\varpi}_{S,F} : \overline{L}_S \rightarrow {\cal M}_S$. We call 
$\varpi_{S,F} : L_{S} \rightarrow {\cal F}_{S}$ or $\overline{\varpi}_{S,F} : \overline{L}_S \rightarrow {\cal M}_S$ the {\it arithmetic prequntization $F$-line bundle} for $\partial V_S$.

Let $L_S^{x_S}$  be the trivial $F$-line bundle over ${\cal F}_S$:
$$ L_S^{x_S} := {\cal F}_S \times F,$$
on which $G$ acts from the right by
$$  L_S^{x_S} \times G  \rightarrow L_S^{x_S}; \; ((\rho_S, z), g) \mapsto (\rho_S.g, z \zeta_N^{\lambda_S^{x_S}(g,\rho_S)}). $$
\\
{\bf Proposition 2.1.39.} {\em We have the following isomorphism of $G$-equivariant $F$-line bundles over ${\cal F}_S$}:
$$ \Phi_{S,F}^{x_S} :   L_S \stackrel{\sim}{\longrightarrow} L_S^{x_S}; \; [([\alpha_S], z)] \mapsto (\varpi_S([\alpha_S]), z\zeta_N^{\sum_{i=1}^r (\alpha_{\frak{p}_i} - x_{\frak{p}_i}(\varpi_{\frak{p}_i}(\alpha_{\frak{p}_i})))})$$
{\em For another section $x'_{S}$, we the following isomorphism of $G$-equivariant $F$-line bundles over ${\cal F}_{S}$}
$$ \Phi_{S, F}^{x_{S}, x'_{S}} : L_{S}^{x_{S}} \stackrel{\sim}{\longrightarrow} L_{S}^{x'_{S}} : [(\rho_{S},z)] \mapsto [(\rho_{S}, z\zeta_N^{\delta_{S}^{x_{S}, x'_{S}}(\rho_{S})})],$$
{\em where $\delta_{S}^{x_{S}, x'_{S}} : {\cal F}_{S} \rightarrow \mathbb{Z}/N\mathbb{Z}$ is the map in Proposition 2.1.34. For $x_{S}, x'_{S}, x''_{S} \in \Gamma({\cal F}_{S}, {\cal L}_{S})$,  we have the equalities}
$$ \left\{ \begin{array}{l} 
\Phi_{S,F}^{x_{S}, x'_{S} } \circ \Phi_{S,F}^{x_{S} } = \Phi_{S,F}^{x'_{S} },\\
 \Phi_{S,F}^{x_{S}, x_{S} } = {\rm id}_{{\cal L}_S^{x_{S}}}, \;  \Phi_{S,F}^{x'_{S}, x_{S} } =  (\Phi_{S,F}^{x_{S}, x'_{S} })^{-1}, \;  \Phi_{S,F}^{x'_{S}, x''_{S} } \circ  \Phi_{S,F}^{x_{S}, x'_{S} } =  \Phi_{S,F}^{x_{S}, x''_{S}}. \end{array} \right.$$
\\
\\
{\it Proof.} The assertions  can be proved in the similar manner to those of the assertions in Proposition 2.1.24, by taking the sum over $\frak{p}_i \in S$. $\;\; \Box$\\
\\
Taking the quotient by the action of $G$, we obtain the line $F$-bundle $\overline{\varpi}_{S,F}^{x_S} : \overline{L}_S^{x_S} \rightarrow {\cal M}_S$. We call $\varpi_{S,F}^{x_S} : L_S^{x_S} \rightarrow {\cal F}_S$ or $\overline{\varpi}_{S,F}^{x_S} : \overline{L}_S^{x_S} \rightarrow {\cal M}_S$ the {\em arithmetic prequantization  $F$-line bundle} for $\partial V_S$ with respect to $x_S$. 
\\

We may also give the description of $L_S$ in terms of the tensor product of $F$-line bundles. Let $p_i : {\cal F}_S \rightarrow {\cal F}_{\frak{p}_i}$ be the $i$-th projection. Let $p_i^*(L_{\frak{p}_i})$ be the $F$-line bundle over ${\cal F}_S$  induced from $L_{\frak{p}_i}$ by $p_i$:
$$ p_i^*(L_{\frak{p}_i}) := \{  (\rho_S, [(\alpha_{\frak{p}_i}, z_i)]) \in {\cal F}_S \times L_{\frak{p}_i} \, | \, p_i(\rho_S) = \varpi_{\frak{p}_i}(\alpha_{\frak{p}_i}) \},$$
and let 
$$ p_i^*(\varpi_{\frak{p}_i}) : p_i^*(L_{\frak{p}_i}) \longrightarrow {\cal F}_S; \; (\rho_S, [(\alpha_{\frak{p}_i}, z_i)]) \mapsto \rho_S$$
be the induced projection. The fiber over $\rho_S = (\rho_{\frak{p}_i})$ is given by
$$ \begin{array}{ll} p_i^*(\varpi_{\frak{p}_i})^{-1}(\rho_S) & = \{ \rho_S \} \times \{ [(\alpha_{\frak{p}_i}, z_i)] \in L_{\frak{p}_i} \, | \, \rho_{\frak{p}_i} = \varpi_{\frak{p}_i}(\alpha_{\frak{p}_i}), z_i \in F \}\\
& \simeq L_{\frak{p}_i}(\rho_{\frak{p}_i})\\
& \simeq F,
\end{array}$$
where $L_{\frak{p}_i}(\rho_{\frak{p}_i})$ is as in (2.1.22).  Let $L_{\frak{p}_1} \boxtimes \cdots \boxtimes L_{\frak{p}_r}$ be the tensor product of $p_i^*(L_{\frak{p}_i})$'s: 
$$ L_{\frak{p}_1} \boxtimes \cdots \boxtimes L_{\frak{p}_r} := p_1^*(L_{\frak{p}_1}) \otimes \cdots \otimes p_r^*(L_{\frak{p}_r}),$$
which is an $F$-line bundle over ${\cal F}_S$.  An element of  $L_{\frak{p}_1} \boxtimes \cdots \boxtimes L_{\frak{p}_r}$ is written by
$$ (\rho_S, [(\alpha_{\frak{p}_1},z_1)]\otimes \cdots \otimes [(\alpha_{\frak{p}_r},z_r)]),$$
where $\rho_S = (\rho_{\frak{p}_i}) \in {\cal F}_S$, $[(\alpha_{\frak{p}_i}, z_i)] \in L_{\frak{p}_i}(\rho_{\frak{p}_i})$. Let $\varpi_S^{\boxtimes} : L_{\frak{p}_1} \boxtimes \cdots \boxtimes L_{\frak{p}_r}  \rightarrow {\cal F}_S$ be the projection.  For fiber over $\rho_S$, we have
$$ (\varpi_S^{\boxtimes})^{-1}(\rho_S) \stackrel{\sim}{\rightarrow} F; (\rho_S, [(\alpha_{\frak{p}_1},z_1)]\otimes \cdots \otimes [(\alpha_{\frak{p}_r},z_r)]) \mapsto \prod_{i=1}^r z_i.  \leqno{(2.1.40)}$$ 
The right action of $G$ on $ L_{\frak{p}_1} \boxtimes \cdots \boxtimes L_{\frak{p}_r} $ is given by 
$$  \begin{array}{l} L_{\frak{p}_1} \boxtimes \cdots \boxtimes L_{\frak{p}_r} \times G \rightarrow L_{\frak{p}_1} \boxtimes \cdots \boxtimes L_{\frak{p}_r}; \\
((\rho_S, [(\alpha_{\frak{p}_1},z_1)]\otimes \cdots \otimes [(\alpha_{\frak{p}_r} ,z_r)]), g) \mapsto (\rho_S.g, [(\alpha_{\frak{p}_1}.g ,z_1)]\otimes \cdots \otimes [(\alpha_{\frak{p}_r}.g,z_r)] ).
\end{array} \leqno{(2.1.41)}$$
The projection $\varpi_S^{\boxtimes}$ is $G$-equivariant. \\
\\
{\bf Proposition 2.1.42.} {\it We have the following isomorphism of $G$-equivariant $F$-line bundles over ${\cal F}_S$}
$$  \begin{array}{l} \Phi_{S,F}^{\boxtimes} :  L_{\frak{p}_1} \boxtimes \cdots \boxtimes L_{\frak{p}_r} \stackrel{\sim}{\longrightarrow} L_S; \\
 (\rho_S, [(\alpha_{\frak{p}_1},z_1)]\otimes \cdots \otimes [(\alpha_{\frak{p}_r} ,z_r)]) \mapsto [([\alpha_S], \prod_{i=1}^r z_i)],
\end{array}$$
where $\rho_{S} = (\rho_{\frak{p}_i}) \in {\cal F}_S$, $[(\alpha_{\frak{p}_i}, z_i)] \in L_{\frak{p}_i}(\rho_{\frak{p}_i})$, and $\alpha_S = (\alpha_{\frak{p}_1},\dots , \alpha_{\frak{p}_r})$. \\
\\
{\em Proof.} If $(\alpha_{\frak{p}_i}, z_i)$ is changed to $(\alpha_{\frak{p}_i}.m_i, \zeta_N^{-m_i}z_i)$ for $m_i \in \mathbb{Z}/N\mathbb{Z}$, $(\alpha_S, \prod_{j=1}^r z_j)$ is changed to $([\alpha_S].m_i, \zeta_N^{-m_i}\prod_{j=1}^r z_j) \sim ([\alpha_S], \prod_{j=1}^r z_j)$. So, by (2.1.20) and (2.1.36),  $\Phi_{S,F}^{\boxtimes}$ is well-defined. \\
(i) It is easy to see that $\varpi_{S,F} \circ \Phi_{S,F}^{\boxtimes} = \varpi_S^{\boxtimes}$.\\
(ii) By (2.1.40), $\Phi_{S,F}^{\boxtimes}$ restricted to a fiber over $\rho_S$  is $F$-linear. \\
(iii) By (2.1.30), (2.1.37) and (2.1.41), we see that $\Phi_{S,F}^{\boxtimes}$ is $G$-equivariant. \\
Therefore $\Phi_{S,F}^{\boxtimes}$ is a morphism of $G$-equivariant $F$-line bundles over ${\cal F}_S$. The inverse is given by
$$ \begin{array}{l} (\Phi_{S,F}^{\boxtimes})^{-1} :  L_S \stackrel{\sim}{\longrightarrow} L_{\frak{p}_1} \boxtimes \cdots \boxtimes L_{\frak{p}_r}; \\
 ([\alpha_S],z) \mapsto  (\varpi_S([\alpha_S]), [(\alpha_{\frak{p}_1}, z)]\otimes [(\alpha_{\frak{p}_2}, 1)] \otimes \cdots \otimes [(\alpha_{\frak{p}_r},1)]),
\end{array}$$
Hence $\Phi_{S,F}^{\boxtimes}$ is a $G$-equivariant isomorphism. $\;\; \Box$\\
\\
{\bf 2.2. Arithmetic Chern-Simons functionals.}  Let ${\cal O}_k$ be the ring of integers of $k$. Let  $X_k := {\rm Spec}({\cal O}_k)$ and let $X_k^{\infty}$ denote the set of infinite primes of $k$. We set $\overline{X}_k := X_k \sqcup X_k^{\infty}$.  
 Let $S = \{ \frak{p}_1, \dots , \frak{p}_r \}$ be a finite set of finite primes of $k$.  Let $\overline{X}_S := \overline{X}_k \setminus S$. We denote by $\Pi_S$ the modified \'{e}tale fundamental group of $\overline{X}_S$ with geometric base point ${\rm Spec}(\overline{k})$ ($\overline{k}$ being a fixed algebraic closure of $k$), which is the Galois group of the maximal subextension $k_S$ of $\overline{k}$ over $k$, unramified outside $S$ (cf. [H; 2.1]). We assume that all maximal ideals of ${\cal O}_k$ dividing $N$ are contained in $S$ (in particular, $S$ is non-empty).

Let ${\cal F}_{\overline{X}_S}$ denote the set of continuous representations of $\Pi_S$ to $G$:
$$ {\cal F}_{\overline{X}_S} := {\rm Hom}_{\rm cont}(\Pi_S, G),$$
on which $G$ acts from the right by
$$ {\cal F}_{\overline{X}_S} \times G \rightarrow {\cal F}_{\overline{X}_S}; \;\; (\rho, g) \mapsto \rho.g := g^{-1} \rho g,  \leqno{(2.2.1)}$$
 and let ${\cal M}_{\overline{X}_S}$ denote the quotient set by this action:
$$ {\cal M}_{\overline{X}_S} := {\cal F}_{\overline{X}_S}/G.$$
Let ${\rm Map}({\cal F}_{\overline{X}_S}, \mathbb{Z}/N\mathbb{Z})$ be the additive group of maps from ${\cal F}_{\overline{X}_S}$ to $\mathbb{Z}/N\mathbb{Z}$, on which $G$ acts from the left by
$$ (g.\psi)(\rho) := \psi(\rho.g)      \leqno{(2.2.2)}$$
for $g \in G, \psi \in {\rm Map}({\cal F}_{\overline{X}_S}, \mathbb{Z}/N\mathbb{Z})$ and $\rho \in {\cal F}_{\overline{X}_S}$. 

We fix an embedding $\overline{k} \hookrightarrow \overline{k}_{\frak{p}_i}$, which induces the continuous homomorphism for each $1\leq i \leq r$
$$ \iota_{\frak{p}_i} : \Pi_{\frak{p}_i} \longrightarrow \Pi_{S}. $$
Let ${\rm res}_{\frak{p}_i}$ and ${\rm res}_S$ denote the restriction maps (the pull-backs by $\iota_{\frak{p}_i}$) defined by 
$$ \begin{array}{l}  {\rm res}_{\frak{p}_i} : {\cal F}_{\overline{X}_S} \longrightarrow {\cal F}_{\frak{p}_i}; \;\; \rho \mapsto \rho \circ \iota_{\frak{p}_i},\\
 {\rm res}_S := ({\rm res}_{\frak{p}_i}) : {\cal F}_{\overline{X}_S} \longrightarrow {\cal F}_{S}; \;\; \rho \mapsto (\rho \circ \iota_{\frak{p}_i}),
\end{array}  \leqno{(2.2.3)}$$
which are $G$-equivariant by (2.1.1), (2.1.25) and (2.2.1). We denote by  ${\rm Res}_{\frak{p}_i}$ and
${\rm Res}_S$ the homomorphisms on cochains defined by 
$$ \begin{array}{l}
{\rm Res}_{\frak{p}_i} : C^n(\Pi_S,\mathbb{Z}/N\mathbb{Z}) \longrightarrow C^n(\Pi_{\frak{p}_i}, \mathbb{Z}/N\mathbb{Z}); \;\; \alpha \mapsto \alpha \circ \iota_{\frak{p}_i},\\
{\rm Res}_S := ({\rm Res}_{\frak{p}_i}) : C^n(\Pi_S,\mathbb{Z}/N\mathbb{Z}) \longrightarrow \prod_{i=1}^r C^n(\Pi_{\frak{p}_i}, \mathbb{Z}/N\mathbb{Z}); \;\; \alpha \mapsto (\alpha \circ \iota_{\frak{p}_i}).\\
\end{array} \leqno{(2.2.4)}
$$
Firstly, we note the following\\
\\
{\bf Lemma 2.2.5.} {\em We have}
$$ H^3(\Pi_S, \mathbb{Z}/N\mathbb{Z}) = 0.$$
{\em Proof.} It suffices to show that the $p$-primary part $H^3(\Pi_S, \mathbb{Z}/N\mathbb{Z})(p) = 0$ for any prime number $p$. Since $H^3(\Pi_S, \mathbb{Z}/N\mathbb{Z})(p) = 0$ for $p \nmid N$, we may assume that $p \mid N$.\\
Case that $N>2$. Then $k$ is totally imaginary and so $\Pi_S = \Pi_{S \cup X_k^{\infty}}$ ($\Pi_{S \cup X_k^{\infty}} := \pi_1^{\mbox{\'{e}t}}({\rm Spec}({\cal O}_k \setminus S)$ being the Galois group of the maximal extension of $k$ unramified outside $S \cup X_k^{\infty}$). By our assumption on $S$,  all primes over $p$ are contained in $S$. So  the cohomological $p$-dimension ${\rm cd}_{p}(\Pi_S) \leq 2$ by  [NSW; Proposition 8.3.18]. Hence $H^3(\Pi_S, \mathbb{Z}/N\mathbb{Z})(p) = 0$.\\
Case that $N=2$ and so $p=2$. Since $S$ does not contain any real primes of $k$, the cohomological $2$-dimension ${\rm cd}_2(\Pi_S) \leq 2$ by [NSW;  Theorem 10.6.7]. Hence $H^3(\Pi_S, \mathbb{Z}/2\mathbb{Z})(2) = 0$. $\;\; \Box$
\\

Let $\rho \in {\cal F}_{\overline{X}_S}$ and so  $c \circ \rho \in Z^3(\Pi_S,\mathbb{Z}/N\mathbb{Z})$. By Lemma 2.2.5, there is $\beta_{\rho} \in C^2(\Pi_S,\mathbb{Z}/N\mathbb{Z})/B^2(\Pi_S,\mathbb{Z}/N\mathbb{Z})$ such that
$$ c \circ \rho = d\beta_{\rho},  \leqno{(2.2.6)}$$
where $d : C^2(\Pi_S,\mathbb{Z}/N\mathbb{Z}) \rightarrow C^3(\Pi_S,\mathbb{Z}/N\mathbb{Z})$ is the coboundary homomorphism. 
By (2.2.3), (2.2.4) and (2.2.6), we see that 
$$  c \circ {\rm res}_{\frak{p}_i}(\rho) = d\, {\rm Res}_{\frak{p}_i}(\beta_{\rho}) \leqno{(2.2.7)} $$
for  $1\leq i \leq r$. By (2.1.4) , (2.1.27) and (2.2.7),  we have 
$$ [{\rm Res}_{S}(\beta_{\rho})]  \in {\cal L}_{S}({\rm res}_S(\rho)).   \leqno{(2.2.8)}$$

Let ${\rm res}_S^{*}({\cal L}_S)$ be the $G$-equivariant principal $\mathbb{Z}/N\mathbb{Z}$-bundle over ${\cal F}_{\overline{X}_S}$ induced from ${\cal L}_S$ by ${\rm res}_S$:
$$  {\rm res}_S^{*}({\cal L}_S) := \{ (\rho, \alpha_S) \in {\cal F}_{\overline{X}_S} \times {\cal L}_S \, | \, {\rm res}_S(\rho) = \varpi_S(\alpha_S) \}. \leqno{(2.2.9)}$$
 and let ${\rm res}_S^{*}(\varpi_S)$ be the projection $ {\rm res}_S^{*}({\cal L}_S) \rightarrow {\cal F}_{\overline{X}_S}$.  The quotient by the action of $G$ is the principal $\mathbb{Z}/N\mathbb{Z}$-bundle ${\rm res}^*(\overline{\cal L}_S)$ over ${\cal M}_{\overline{X}_S}$ induced from $\overline{\cal L}_S$ by ${\rm res}_S$. By (2.2.9), a section of ${\rm res}_S^{*}(\varpi_S)$ is naturally identified with a map $y_S : {\cal F}_{\overline{X}_S} \rightarrow {\cal L}_S$ satisfying $\varpi_S \circ y_S = {\rm res}_S$:
$$ \Gamma({\cal F}_{\overline{X}_S}, {\rm res}_S^{*}({\cal L}_S)) = \{ y_S : {\cal F}_{\overline{X}_S} \rightarrow {\cal L}_S \, | \, \varpi_S \circ y_S = {\rm res}_S
\},  \leqno{(2.2.10)}$$
on which $G$ acts by $(g.y_S)(\rho) := y_S(\rho.g)$ for $\rho \in {\cal F}_{\overline{X}_S}, g \in G$. Let $\Gamma_G({\cal F}_{\overline{X}_S}, {\rm res}_S^{*}({\cal L}_S))$ denote the set of $G$-equivariant sections of ${\rm res}_S^{*}(\varpi_S)$. 
We define the (mod $N$)  {\em arithmetic Chern-Simons functional} $CS_{\overline{X}_S} : {\cal F}_{\overline{X}_S} \rightarrow {\cal L}_S$ by
$$ CS_{\overline{X}_S}(\rho) := [ {\rm Res}_S(\beta_{\rho})]   \leqno{(2.2.11)} $$
for $\rho \in {\cal F}_{\overline{X}_S}$. The value $CS_{\overline{X}_S}(\rho) \in {\cal L}_S$ is called the {\it arithmetic Chern-Simons invariant} of $\rho$.\\
\\
{\bf Lemma 2.2.12.}  (1) {\it $CS_{\overline{X}_S}(\rho)$ is independent of the choice of $\beta_{\rho}$.}\\
(2) $CS_{\overline{X}_S}$ is a $G$-equivariant section of ${\rm res}_S^{*}(\varpi_S)$:
$$ CS_{\overline{X}_S} \in \Gamma_G({\cal F}_{\overline{X}_S}, {\rm res}_S^{*}({\cal L}_S)) = \Gamma({\cal M}_{\overline{X}_S}, {\rm res}_S^{*}(\overline{{\cal L}}_S)).$$
\\
{\it Proof.} (1) Let $\beta'_{\rho} \in C^2(\Pi_S,\mathbb{Z}/N\mathbb{Z})/B^2(\Pi_S,\mathbb{Z}/N\mathbb{Z})$ be another choice  satisfying $c \circ \rho = d\beta'_{\rho}.$ Then we have $\beta'_{\rho} = \beta_{\rho} + z$ for some $z \in H^2(\Pi_S, \mathbb{Z}/N\mathbb{Z})$ and so
$$ {\rm Res}_{\frak{p}_i}(\beta'_{\rho}) - {\rm Res}_{\frak{p}_i}(\beta_{\rho}) = {\rm inv}_{\frak{p}_i}({\rm Res}_{\frak{p}_i}(z)) \;\; (1\leq i \leq r).$$
Noting that any primes dividing $N$ is contained in $S$, Tate-Poitou exact sequence ([NSW; 8.6.10]) implies that the composite of the following maps
$$ H^2(\Pi_S, \mathbb{Z}/N\mathbb{Z}) \stackrel{\prod_{ \frak{p} \in \overline{S}} {\rm Res}_{\frak{p}}}{\longrightarrow} \prod_{\frak{p} \in \overline{S}} H^2(\Pi_{\frak{p}}, \mathbb{Z}/N\mathbb{Z}) \stackrel{\sum_{\frak{p} \in \overline{S}} {\rm inv}_{\frak{p}}}{\longrightarrow} \mathbb{Z}/N\mathbb{Z}$$
is the zero map, where $\overline{S} = S \cup X_k^{\infty}$. For any infinite prime $v \in X_k^{\infty}$, the restriction map $\Pi_v := {\rm Gal}(\overline{k}_v/k_v) \rightarrow \Pi_S = {\rm Gal}(k_S/k)$ is the trivial homomorphism, because any infinite prime is unramified in $k_S/k$. So ${\rm Res}_v : H^2(\Pi_S, \mathbb{Z}/N\mathbb{Z}) \rightarrow H^2(\Pi_v, \mathbb{Z}/N\mathbb{Z})$ is the zero map. Hence we have
$$ \sum_{i=1}^r {\rm inv}_{\frak{p}_i}({\rm Res}_{\frak{p}_i}(z)) = 0.$$
By (2.1.28), we obtain 
$$[{\rm Res}_S(\beta'_{\rho})] =  [{\rm Res}_S(\beta_{\rho})]. $$
(2) By (2.2.8), (2.2.10) and (2.2.11), we have
$$ CS_{\overline{X}_S} \in \Gamma({\cal F}_{\overline{X}_S}, {\rm res}_S^*({\cal L}_S)).$$
So it suffices to show that $CS_{\overline{X}_S}$ is $G$-equivariant. By (2.1.5) and (2.2.6), we have
$$ d \beta_{\rho.g}  = c \circ (\rho.g) = (g.c) \circ \rho \\
                                                     = (c + d h_g)\circ \rho \\
                                                      = d(\beta_{\rho} + h_g \circ \rho).
 $$
for $g \in G$ and $\rho \in {\cal F}_{\overline{X}_S}$. Therefore there is $z \in H^2(\Pi_{S}, \mathbb{Z}/\mathbb{Z})$ such that $\beta_{\rho.g} = \beta_{\rho} + h_g \circ \rho + z$ and so
$$  \begin{array}{ll}  {\rm Res}_{S}(\beta_{\rho.g}) & = {\rm Res}_{S}(\beta_{\rho}) + h_g \circ {\rm res}_{S}(\rho) + {\rm Res}_{S}(z)\\
& =   {\rm Res}_S(\beta_{\rho}).g + {\rm Res}_S(z).
\end{array}$$    
By the same argument as in (1) above, we obtain
$$ CS_{\overline{X}_S}(\rho.g) = [{\rm Res}_S(\beta_{\rho.g})] = [{\rm Res}_S(\beta_{\rho})].g = CS_{\overline{X}_S}(\rho).g. \;\; \Box$$ 
\\                          

Let $x_S = [ (x_{\frak{p}_1}, \dots , x_{\frak{p}_r})] \in \Gamma({\cal F}_S, {\cal L}_S)$ be a section and  let
${\cal L}_S^{x_S}$ be the arithmetic prequantization principal $\mathbb{Z}/N\mathbb{Z}$-bundle over ${\cal F}_S$ with respect to $x_S$. Let ${\rm res}_S^{*}({\cal L}_S^{x_S})$ be the $G$-equivariant principal $\mathbb{Z}/N\mathbb{Z}$-bundle over ${\cal F}_{\overline{X}_S}$ induced from ${\cal L}_S^{x_S}$ by ${\rm res}_S$:
$$ \begin{array}{ll}
{\rm res}_S^{*}({\cal L}_S^{x_S}) & = \{ (\rho, (\rho_S,m)) \in {\cal F}_{\overline{X}_S} \times {\cal L}_S^{x_S} \, | \, {\rm res}_S(\rho) = \rho_S \}\\
                     & = {\cal F}_{\overline{X}_S} \times \mathbb{Z}/N\mathbb{Z}\\
                     \end{array}$$
                     by identifying $(\rho, (\rho_S,m))$ with $(\rho,m)$. So a section of ${\rm res}_S^{*}({\cal L}_S^{x_S})$ over ${\cal F}_{\overline{X}_S}$ is identified with a map ${\cal F}_{\overline{X}_S} \rightarrow \mathbb{Z}/N\mathbb{Z}$:
                     $$ \Gamma({\cal F}_{\overline{X}_S}, {\rm res}_S^{*}({\cal L}_S^{x_S})) = {\rm Map}({\cal F}_{\overline{X}_S},\mathbb{Z}/N\mathbb{Z}),$$
                     on which $G$ acts by (2.2.2). 
Therefore, letting  ${\rm Map}_G({\cal F}_{\overline{X}_S},\mathbb{Z}/N\mathbb{Z})$ denote the set of $G$-equivariant maps ${\cal F}_{\overline{X}_S} \rightarrow \mathbb{Z}/N\mathbb{Z}$,
we have the identification
$$ \begin{array}{ll} \Gamma_G({\cal F}_{\overline{X}_S}, {\rm res}_S^{*}({\cal L}_S^{x_S})) & = {\rm Map}_G({\cal F}_{\overline{X}_S},\mathbb{Z}/N\mathbb{Z}) \\
& = \{ \psi : {\cal F}_{\overline{X}_S} \rightarrow \mathbb{Z}/N\mathbb{Z}\, | \, \psi(\rho.g) = \psi(\rho) + \lambda_S^{x_S}(g, {\rm res}_S(\rho))\\
&  \;\;\;\;\;\;\;\;\;\;\;\;\;\;\;\;\;\;\;\;\;\;\;\;\;\;\;\;  \;\;\;\;\;\;\;\;\;\;\;\;\;\;\;\;\;\;\;\;\;\;\;\;\;\;\;\;  \mbox{for}\; \rho \in {\cal F}_{\overline{X}_S}, g \in G \}. \end{array}$$
The isomorphism $\Phi_S^{x_S} : {\cal L}_S \stackrel{\sim}{\rightarrow} {\cal L}_S^{x_S}$ in Proposition 2.1.35 induces the isomorphism
$$ \begin{array}{rcl} \Psi^{x_S} : \Gamma_G({\cal F}_{\overline{X}_S}, {\rm res}_S^{*}({\cal L}_S)) &  \stackrel{\sim}{\longrightarrow} & \Gamma_G({\cal F}_{\overline{X}_S}, {\rm res}_S^{*}({\cal L}_S^{x_S})) = {\rm Map}_G({\cal F}_{\overline{X}_S},\mathbb{Z}/N\mathbb{Z})\\
 y_S & \mapsto & \Phi_S^{x_S} \circ y_S.
 \end{array}
$$
We then define the  {\em arithmetic Chern-Simons functional} $CS_{\overline{X}_S}^{x_S} : {\cal F}_{\overline{X}_S} \rightarrow \mathbb{Z}/N\mathbb{Z}$ with respect to $x_S$ by the image of $ CS_{\overline{X}_S}$ under $\Psi^{x_S}$:
$$ CS_{\overline{X}_S}^{x_S} := \Psi^{x_S}(CS_{\overline{X}_S}). \leqno{(2.2.13)}$$
\\
{\bf Theorem 2.2.14.} (1) {\em For $\rho \in {\cal F}_{\overline{X}_S}$, we have}
$$ CS_{\overline{X}_S}^{x_S}(\rho) = \sum_{i=1}^r ({\rm Res}_{\frak{p}_i}(\beta_{\rho})  - x_{\frak{p}_i}({\rm res}_{\frak{p}_i}(\rho))),$$
which is independent of the choice of $\beta_{\rho}$.\\
(2) {\em We have the following equality in } $ C^1(G, {\rm Map}({\cal F}_{\overline{X}_S}, \mathbb{Z}/N\mathbb{Z}))$
$$ d CS_{\overline{X}_S}^{x_S} = {\rm res}^*(\lambda_S^{x_S}). $$
\\
{\em Proof.} (1) This follows from the definition of $\Phi_S^{x_S}$ in Proposition 2.1.35 and (2.2.13).\\
(2) Since $CS_{\overline{X}_S}^{x_S} \in {\rm Map}_G({\cal F}_{\overline{X}_S}, \mathbb{Z}/N\mathbb{Z})$, we have
$$ CS_{\overline{X}_S}^{x_S}(\rho.g) = CS_{\overline{X}_S}^{x_S}(\rho) + \lambda_S^{x_S}(g, {\rm res}_S(\rho))$$
for $g \in G$ and $\rho \in {\cal F}_{\overline{X}_S}$, which means the assertion.  $\;\; \Box$\\
\\
{\bf Proposition 2.2.15.} {\em Let $x'_S \in \Gamma({\cal F}_S, {\cal L}_S)$ be another section which yields $CS_{\overline{X}_S}^{x'_S}$, and let $\delta_S^{x_S, x'_S} : {\cal F}_S \rightarrow \mathbb{Z}/N\mathbb{Z}$ be the map in Proposition 2.1.34.
 Then we have}
$$ CS_{\overline{X}_S}^{x'_S}(\rho) - CS_{\overline{X}_S}^{x_S}(\rho) = \delta_S^{x_S,x'_S}({\rm res}_S(\rho)).$$
\\
{\it Proof.} By Proposition 2.2.14 (1) and Lemma 1.1.4 (1), we have
$$ \begin{array}{ll} CS_{\overline{X}_S}^{x'_S}(\rho) - CS_{\overline{X}_S}^{x_S}(\rho) & = \displaystyle{ \sum_{i=1}^r ({\rm Res}_{\frak{p}_i}(\beta_{\rho}) - x'_{\frak{p}_i}({\rm res}_{\frak{p}_i}(\rho))) - \sum_{i=1}^r ({\rm Res}_{\frak{p}_i}(\beta_{\rho}) - x_{\frak{p}_i}({\rm res}_{\frak{p}_i}(\rho))) }\\
& = \displaystyle{ \sum_{i=1}^r (x_{\frak{p}_i}({\rm res}_{\frak{p}_i}(\rho)) - x'_{\frak{p}_i}({\rm res}_{\frak{p}_i}(\rho))) }\\
& = \delta_S^{x_S, x'_S}({\rm res}_S(\rho)). \;\;\;\;\;\;\;\; \Box
\end{array}$$
\\

For $x_S, x'_S \in \Gamma({\cal F}_S, {\cal L}_S)$, the $G$-equivariant isomorphism $\Phi_S^{x_S,x'_S} : {\cal L}_S^{x_S} \stackrel{\sim}{\rightarrow} {\cal L}_S^{x'_S}$ induces the isomorphism
$$ \Psi^{x_S,x'_S} : \Gamma_G({\cal F}_{\overline{X}_S}, {\rm res}_S^{*}({\cal L}_S^{x_S})) \stackrel{\sim}{\longrightarrow} \Gamma_G({\cal F}_{\overline{X}_S}, {\rm res}_S^{*}({\cal L}_S^{x_S})); \; \psi^{x_S} \mapsto \Phi_S^{x_S, x'_S} \circ \psi^{x_S}.$$
By Proposition 2.1.35,  we have
$$ \left\{ \begin{array}{l} \Psi^{x_S,x'_S} \circ \Psi^{x_S} = \Psi^{x'_S}\\
\Psi^{x_S,x_S} ={\rm id}_{\Gamma_G({\cal F}_{\overline{X}_S}, {\rm res}_S^{*}({\cal L}_S^{x_S}))}, \Psi^{x'_S,x_S} = (\Psi^{x_S,x'_S})^{-1},  \Psi^{x'_S,x''_S} \circ \Psi^{x_S,x'_S} =\Psi^{x_S,x''_S}.
\end{array} \right.$$
So we can define the equivalence relation $\sim$ on the disjoint union of $\Gamma_G({\cal F}_{\overline{X}_S}, {\rm res}_S^{*}({\cal L}_S^{x_S})) $ over $x_S \in \Gamma({\cal F}_S, {\cal L}_S)$ by
$$ \psi^{x_S} \sim \psi^{x'_S}  \Longleftrightarrow \Psi^{x_S, x'_S}(\psi^{x_S}) = \psi^{x'_S} $$
for $\psi^{x_S} \in \Gamma_G({\cal F}_{\overline{X}_S}, {\rm res}_S^{*}({\cal L}_S^{x_S}))$ and $\psi^{x'_S} \in \Gamma_G({\cal F}_{\overline{X}_S}, {\rm res}_S^{*}({\cal L}_S^{x_S}))$.Since $\Phi_S^{x'_S} = \Phi_S^{x_S,x'_S} \circ \Phi_S^{x_S}$, $CS_{\overline{X}_S}^{x_S} \simeq CS_{\overline{X}_S}^{x'_S}$. Thus we have the following identification:
$$ \begin{array}{ccc} \Gamma_G({\cal F}_{\overline{X}_S}, {\rm res}_S^{*}({\cal L}_S)) & = &  \bigsqcup_{x_S \in \Gamma({\cal F}_S, {\cal L}_S)} \Gamma_G({\cal F}_{\overline{X}_S}, {\rm res}_S^{*}({\cal L}_S^{x_S})) / \sim; \\
\psi & \mapsto & [\Psi^{x_S}(\psi)] 
\end{array}
\leqno{(2.2.16)}
$$
where $CS_{\overline{X}_S}$ and $[CS_{\overline{X}_S}^{x_S}]$ are identified.
\\
\\
\begin{center}
{\bf 3. Quantum theory}
\end{center}

In this section, we construct the arithmetic quantum space and the arithmetic Dijkgraaf-Witten invariant over the moduli space of Galois representations. These constructions correspond to the quantum theory of topological Dijkgraaf-Witten TQFT. We keep the same notations and assumptions as in Section 2. We assume that  $F$ is a subfield of $\mathbb{C}$ such that $\zeta_N$ is contained in $F$ and $\overline{F} = F$ ($\overline{F}$ being the complex conjugate).  \\
\\
{\bf 3.1. Arithmetic quantum spaces.} Following the construction of the quantum Hilbert space, we define the {\em arithmetic quantum space} ${\cal H}_S$ for $\partial V_S$ by the space of $G$-equivariant sections of the arithmetic prequantization $F$-line bundle $\varpi_{S,F} : L_S \rightarrow {\cal F}_S$:
$$ {\cal H}_S := \Gamma_G({\cal F}_S, L_S) = \Gamma({\cal M}_S, \overline{L}_S).$$
It is a finite dimensional $F$-vector space. 

Let $x_S = [ (x_{\frak{p}_1}, \dots , x_{\frak{p}_r})] \in \Gamma({\cal F}_S, {\cal L}_S)$ be a section and  let
$L_S^{x_S}$ be the arithmetic prequantization $F$-line bundle over ${\cal F}_S$ with respect to $x_S$ and let
$$ \begin{array}{ll} {\cal H}_S^{x_S} & := \Gamma_G({\cal F}_S, L_S^{x_S}) = \Gamma({\cal M}_S, \overline{L}_S^{x_S})\\
                                                       & = \{ \theta : {\cal F}_S \rightarrow F \, | \, \theta(\rho_S.g) = \zeta_{N}^{\lambda_S^{x_S}(g, \rho_S)} \theta(\rho_S) \; \mbox{for}\; \rho_S \in {\cal F}_S, g \in G \},
\end{array} 
\leqno{(3.1.1)}$$
 which we call the {\em arithmetic quantum space} for $\partial V_S$ with respect to $x_S$.  The isomorphism $\Phi_{S,F}^{x_S} : L_S \stackrel{\sim}{\rightarrow} L_S^{x_S}$ in Proposition 2.1.39 induces the isomorphism
$$ \Theta^{x_S} : {\cal H}_S \stackrel{\sim}{\longrightarrow} {\cal H}_S^{x_S}; \; \theta  \mapsto  \Phi_{S,F}^{x_S} \circ \theta. \leqno{(3.1.2)}$$
We call an element of ${\cal H}_S$ or ${\cal H}_S^{x_S}$ an {\em arithmetic theta function} (cf. Remark 3.2.4 below). 

For $x_S, x'_S \in \Gamma({\cal F}_S,{\cal L}_S)$, the isomorphism $\Phi_{S,F}^{x_S,x'_S} : L_S^{x_S} \stackrel{\sim}{\rightarrow} L_S^{x'_S}$ induces the isomorphism of $F$-vector spaces:
$$ \Theta^{x_S, x'_S} : {\cal H}_S^{x_S} \stackrel{\sim}{\longrightarrow} {\cal H}_S^{x'_S}; \; \theta^{x_S} \mapsto \Phi_{S,F}^{x_S,x'_S} \circ \theta^{x_S}$$
and, by Proposition 2.1.39, we have
$$ \left\{ \begin{array}{l} \Theta^{x_S,x'_S} \circ \Theta^{x_S} = \Theta^{x'_S}\\
\Theta^{x_S,x_S} = {\rm id}_{{\cal H}_S^{x_s}}, \Theta^{x'_S,x_S} = (\Theta^{x_S,x'_S})^{-1}, \Theta^{x'_S,x''_S} \circ \Theta^{x_S,x'_S} = \Theta^{x_S,x''_S}.
\end{array} \right.
$$
So the equivalence relation $\sim$ is defined on the disjoint union of all ${\cal H}_S^{x_S}$ running over $x_S \in \Gamma({\cal F}_S, {\cal L}_S)$ by
$$ \theta^{x_S} \sim \theta^{x'_S}  \Longleftrightarrow \Theta^{x_S, x'_S}(\theta^{x_S}) = \theta^{x'_S} $$
for $\theta^{x_S} \in {\cal H}_S^{x_S}$ and $\theta^{x'_S} \in {\cal H}_S^{x'_s}$. Then we have the following identification:
$$ {\cal H}_S = \bigsqcup_{x_S \in \Gamma({\cal F}_S, {\cal L}_S)} {\cal H}_S^{x_S} / \sim. \leqno{(3.1.3)}$$
\\
{\bf Remark 3.1.4.} The arithmetic quantum space ${\cal H}_S$ is an arithmetic analog of the quantum Hilbert space ${\cal H}_{\Sigma}$ for a surface $\Sigma$ in (2+1)-dimensional Chern-Simons TQFT. We recall that ${\cal H}_{\Sigma}$ is known to coincides with the space of conformal blocks  ([BL]) and its dimension formula was shown by Verlinde ([V]). It would also be an interesting question in number theory to describe the dimension and a canonical basis of ${\cal H}_S$ in  comparison of Verlinde's formulas.\\
\\
{\bf 3.2.  Arithmetic Dijkgraaf-Witten partition functions.}  For $\rho_S \in {\cal F}_S$, we define the subset ${\cal F}_{\overline{X}_S}(\rho_S)$ of ${\cal F}_{\overline{X}_S}$ by
$$ {\cal F}_{\overline{X}_S}(\rho_S) := \{ \rho \in {\cal F}_{\overline{X}_S} \, | \, {\rm res}_S(\rho) = \rho_S \}.$$
 We then define the {\em arithmetic Dijkgraaf-Witten invariant} $Z_{\overline{X}_S}^{x_S}(\rho_S)$ of $\rho_S$ with respect to $x_S$ by
$$ Z_{\overline{X}_S}^{x_S}(\rho_S) := \displaystyle{ \frac{1}{\# G} \sum_{\rho \in {\cal F}_{\overline{X}_S}(\rho_S)} \zeta_N^{CS_{\overline{X}_S}^{x_S}(\rho)}. } \leqno{(3.2.1)}$$
\\
{\bf Theorem 3.2.2.} (1) {\em $Z_{\overline{X}_S}^{x_S}(\rho_S)$ is independent of the choice of $\beta_{\rho}$}.\\
(2) {\em We have}
$$ Z_{\overline{X}_S}^{x_S} \in {\cal H}_S^{x_S}.$$
\\
{\em Proof.} (1) This follows from Lemma 2.2.12 (1).\\
(2) This follows from Theorem 2.2.14 (2) and (3.2.1). \\
\\
We call $Z_{\overline{X}_S}^{x_S} \in {\cal H}_S^{x_S}$ the {\em arithmetic Dijkgraaf-Witten partition function} for $\overline{X}_S$ with respect to $x_S$. \\
\\
 The following proposition tells us how they are changed when we change $x_S$.\\
\\
{\bf Proposition 3.2.3.} {\em For sections $x_S, x'_S \in \Gamma({\cal F}_S, {\cal L}_S)$, we have}
$$ \Theta^{x_S,x'_S}(Z_{\overline{X}_S}^{x_S}) = Z_{\overline{X}_S}^{x'_S}.$$
\\
{\em Proof.} We have
$$ \begin{array}{ll} \Theta^{x_S,x'_S}(Z_{\overline{X}_S}^{x_S})(\rho_S) & = (\Phi_{S,F}^{x_S,x'_S} \circ Z_{\overline{X}_S})(\rho_S)\\
                                                                             & = Z_{\overline{X}_S}(\rho_S)\zeta_N^{\delta_S^{x_S,x'_S}(\rho_S)} \;\; \mbox{by Proposition 2.1.39} \\
                                                                             & = \displaystyle{ \frac{1}{\# G} \sum_{\rho \in {\cal F}_{\overline{X}_S}(\rho_S)} \zeta_N^{CS_{\overline{X}_S}^{x_S
                                                                           }(\rho) + \delta_S^{x_S,x'_S}(\rho_S)} } \;\; \mbox{by (3.2.1)}\\
                                                                           & 
= \displaystyle{ \frac{1}{\# G} \sum_{\rho \in {\cal F}_{\overline{X}_S}(\rho_S)} \zeta_N^{CS_{\overline{X}_S}^{x'_S}(\rho)} }  \;\; \mbox{by Proposition 2.2.15} \\
                                                                           & = Z_{\overline{X}_S}^{x'_S}(\rho_S)
                                                    \end{array}
                                                    $$
for $\rho_S \in {\cal F}_S$.   So we obtain the assertion.                                                $\; \; \Box$
\\
\\
By the identification (3.1.3), $Z_{\overline{X}_S}^{x_S}$ defines the element $Z_{\overline{X}_S}$ of ${\cal H}_S$ which is independent of the choice of $x_S$. We call it the {\em arithmetic Dijkgraaf-Witten partition function} for $\overline{X}_S$.\\
\\
{\bf Remark 3.2.4.}  In (2+1)-dimensional Chern-Simons TQFT, an element of ${\cal H}_{\Sigma}$ for a surface $\Sigma$ may be regarded as a (non-abelian) generalization of the classical theta function on the Jacobian manifold of $\Sigma$ (cf. [BL]. It goes back to Weli's paper [We]. See [Mo1] for an arithmetic analog.) In this respect, it may be interesting to observe that the Dijkgraaf-Witten partition function in (3.2.1) may look like a variant of (non-abelian) Gaussian sums. 
\\
\\
\begin{center}
{\bf 4. Some basic and functorial properties}
\end{center}

In this section, we study some basic and functorial properties of the objects constructed in Sections 2 and 3. We keep the same notations as in Sections 2 and 3. 
\\
\\
{\bf 4.1. Change of the 3-cocycle $c$.}  The theory given in Sections 2 and 3 depends on a chosen 3-cocycle $c$. We shall see in the following that when $c$ is changed in the cohomology class $[c]$, objects are changed to isomorphic ones, and hence the theory depends essentially on the cohomology class $[c]$.  Let $c' \in Z^3(G, \mathbb{Z}/N\mathbb{Z})$ be another 3-cocycle representing $[c]$. The objects constructed by using $c'$ will be denoted by using $'$, for example, by ${\cal L}'_{\frak{p}}, L'_{\frak{p}}, \dots$ etc.  

There is $b \in C^2(G, \mathbb{Z}/N\mathbb{Z})$ such that $c' -c = db$. Then we have the isomorphism of $\mathbb{Z}/N\mathbb{Z}$-torsors for $\rho_{\frak{p}} \in {\cal F}_{\frak{p}}$:
$$ {\cal L}_{\frak{p}}(\rho_{\frak{p}}) \stackrel{\sim}{\longrightarrow} {\cal L}'_{\frak{p}}(\rho_{\frak{p}}); \; \alpha_{\frak{p}} \mapsto \alpha_{\frak{p}} + b \circ \rho_{\frak{p}},$$
which induces the following isomorphisms of arithmetic quantization bundles:
$$  \begin{array}{ll}
\xi_{\frak{p}} : {\cal L}_{\frak{p}} \stackrel{\sim}{\longrightarrow} {\cal L}'_{\frak{p}}, \; &  \; \xi_{\frak{p},F} : L_{\frak{p}} \stackrel{\sim}{\longrightarrow} L'_{\frak{p}},\\
\xi_S : {\cal L}_S \stackrel{\sim}{\longrightarrow} {\cal L}'_{S},\;  & \; \xi_{S,F} :  L_S \stackrel{\sim}{\longrightarrow} L'_{S}.
\end{array} \leqno{(4.1.1)}
$$
Let $x_{\frak{p}} \in \Gamma({\cal F}_{\frak{p}}, {\cal L}_{\frak{p}})$ and $x_S = [(x_{\frak{p}_1},\dots , x_{\frak{p}_r})] \in \Gamma({\cal F}_S,{\cal L}_S)$, and let $x'_{\frak{p}} \in \Gamma({\cal F}'_{\frak{p}}, {\cal L}'_{\frak{p}})$ and $x'_S \in \Gamma({\cal F}'_F, {\cal L}'_S)$. 
Denote by $\lambda'_{\frak{p}}$ and $\lambda'_S$ the arithmetic Chern-Simons 1-cocycles for $\partial V_{\frak{p}}$ and $\partial V_S$ with respect to $x'_{\frak{p}}$ and $x'_S$, respectively. We define $\kappa_{\frak{p}} : {\cal F}_{\frak{p}} \rightarrow \mathbb{Z}/N\mathbb{Z}$ and $\kappa_S : {\cal F}_S \rightarrow \mathbb{Z}/N\mathbb{Z}$ by 
$$ \kappa_{\frak{p}}(\rho_{\frak{p}}) := (\xi_{\frak{p}}\circ x_{\frak{p}})(\rho_{\frak{p}}) - x'_{\frak{p}}(\rho_{\frak{p}}), \;\; \kappa_S(\rho_S) := \sum_{i=1}^r \kappa_{\frak{p}_i}(\rho_{\frak{p}_i}) $$
for $\rho_{\frak{p}} \in {\cal F}_{\frak{p}}$ and $\rho_S = (\rho_{\frak{p}_1},\dots , \rho_{\frak{p}_r}) \in {\cal F}_S$, respectively. Then we have
$$ \lambda'_{\frak{p}}(g) - \lambda_{\frak{p}}(g) = g.\kappa_{\frak{p}} - \kappa_{\frak{p}}, \;\; \lambda'_S(g) - \lambda_S(g) = g.\kappa_S - \kappa_S.$$
We note that if we take $x'_{\frak{p}} := \xi_{\frak{p}} \circ x_{\frak{p}}$ and $x'_S := \xi_S \circ x_S$, $\kappa_{\frak{p}} = 0$ and so $\kappa_S = 0$. As in Corollary 2.1.19, Propositions 2.1.24, 2.1.35 and 2.1.39, using $\kappa_{\frak{p}}$ and $\kappa_S$, we have the isomorphisms
$$  \begin{array}{ll}
 {\cal L}_{\frak{p}}^{x_{\frak{p}}} \stackrel{\sim}{\longrightarrow} {{\cal L}'_{\frak{p}}}^{x'_{\frak{p}}}, \; &  \; L_{\frak{p}}^{x_{\frak{p}}} \stackrel{\sim}{\longrightarrow} {L'_{\frak{p}}}^{x'_{\frak{p}}},\\
{\cal L}_S^{x_S} \stackrel{\sim}{\longrightarrow} {{\cal L}'_S}^{x'_S},\;  & \;  L_S^{x_S} \stackrel{\sim}{\longrightarrow} {L'_{S}}^{x'_S}.
\end{array}
$$
which are compatible with the isomorphisms in (4.1.1) via the isomorphisms ${\cal L}_{\frak{p}} \simeq {\cal L}_{\frak{p}}^{x_{\frak{p}}}, L_{\frak{p}} \simeq L_{\frak{p}}^{x_{\frak{p}}}, {\cal L}_S \simeq {\cal L}_S^{x_S}$ and $L_S \simeq L_S^{x_S}$ in Propositions 2.1.15, 2.1.24, 2.1.35 and 2.1.39. 

The isomorphism $\xi_S : {\cal L}_S \stackrel{\sim}{\rightarrow} {\cal L}'_S$ induces the isomorphism
$$ \Gamma_G({\cal F}_{\overline{X}_S}, {\rm res}_S^{*}({\cal L}_S)) \stackrel{\sim}{\longrightarrow} \Gamma_G({\cal F}_{\overline{X}_S}, {\rm res}_S^{*}({\cal L}'_S)) $$
which sends $CS_{\overline{X}_S}$ to $CS'_{\overline{X}_S}$, and the isomorphism $\xi_{S,F} : L_S \stackrel{\sim}{\rightarrow} L'_S$ induces the isomorphisms
$$  {\cal H}_S \stackrel{\sim}{\longrightarrow} {\cal H}'_S, \;\; {\cal H}_S^{x_S} \stackrel{\sim}{\longrightarrow} {{\cal H}'_S}^{x'_S}, $$
which sends $Z_{\overline{X}_S}$ to $Z'_{\overline{X}_S}$.
\\
\\
{\bf Remark 4.1.2.}  A cochain $\alpha \in C^n(G,A)$ is called {\em normalized} if  $\alpha(g_1, \dots , g_n) = 0$ whenever one of $g_i$'s is 1. It is known that any cocyle is cohomologous to a normalized one, namely, any cohomology class of $H^n(G,A)$ is represented by a normalized cocycle ([NSW; Chapter I, $\S 2$, Exercise 4], [EM; Lemma 6.1]). Therefore, by the above argument, we may assume that we can take the fixed cocycle $c \in Z^3(G, \mathbb{Z}/N\mathbb{Z})$ in our theory to be normalized.
\\
\\
{\bf 4.2. Change of number fields.} Let $k'$ be an another number field contains a primitive $N$-th root of unity and let $S'=\{ \frak{p}'_1, \dots, \frak{p}'_{r'} \}$ be a finite set of finite primes of $k'$ such that any finite prime dividing $N$ is contained in $S'$. The objects constructed by using $k'$ and $S'$ will be denoted by, for example,  ${\cal L}_{\frak{p}'}, L_{\frak{p}'}, {\cal L}_{S'}, L_{S'}, \dots$ etc, for simplicity of notations
. Assume that $r=r'$ and there are isomorphisms $\xi_i : k_{{\frak{p}}_i} \stackrel{\sim}{\rightarrow} k'_{\frak{p}'_i}$ for $1 \leq i \leq r$.  Then $\xi_i$'s induces the following isomorphisms of arithmetic quantization bundles:
$$\xi_{\frak{p}_i} : {\cal L}_{\frak{p}_i} \stackrel{\sim}{\longrightarrow} {\cal L}_{\frak{p}'_i} , \ \ \ \  \xi_{\frak{p}_i , F} :  L_{\frak{p}_i} \stackrel{\sim}{\longrightarrow}   L_{\frak{p}'_i}$$
$$\xi_S : {\cal L}_S \stackrel{\sim}{\longrightarrow} {\cal L}_{S'} , \ \ \ \  \xi_{S , F} :  L_S \stackrel{\sim}{\longrightarrow}  L_{S'} .$$
Let $x_{\frak{p}_i} \in \Gamma({\cal F}_{\frak{p}_i}, {\cal L}_{\frak{p}_i})$ and $x_S=[(x_{\frak{p}_1}, \dots , x_{\frak{p}_r})] \in \Gamma({\cal F}_S, {\cal L}_S)$, and let  $x_{\frak{p}'_i} \in \Gamma({\cal F}_{\frak{p}'_i}, {\cal L}_{\frak{p}'_i})$ and $x_{S'}=[(x_{\frak{p}'_1}, \dots , x_{\frak{p}'_r})] \in \Gamma({\cal F}_{S'}, {\cal L}_{S'})$. Then we have the isomorphisms of arithmetic prequantization bundles with respect to sections
$$ {\cal L}^{x_{\frak{p}_i}}_{\frak{p}_i} \stackrel{\sim}{\longrightarrow} {\cal L}_{\frak{p}'_i}^{x_{{\frak{p}}'_i} } , \ \ \ \    L^{x_{\frak{p}_i}}_{\frak{p}_i} \stackrel{\sim}{\longrightarrow}  L_{\frak{p}'_i}^{ x_{\frak{p}'_i} }$$
$$ {\cal L}^{x_S}_S \stackrel{\sim}{\longrightarrow} {\cal L}_{S'}^{x_{S'}} , \ \ \ \    L^{x_S}_S \stackrel{\sim}{\longrightarrow}  L_{S'}^{x_{S'}}.$$
Suppose further that there is an isomorphism $\tau : k \stackrel{\sim}{\rightarrow} k'$ of number fields which sends $\frak{p}_i$ to $\frak{p}'_i$ for $1\leq i \leq r$. so that we have the isomorphism 
$$\xi : \overline{X}_S := \overline{X}_k \setminus S \stackrel{\sim}{\longrightarrow} {\overline{X}_{k'}} \backslash S' =: {\overline{X}}_{S'}.$$
For example, let $k := \mathbb{Q}(\sqrt[3]{2}), k' := \mathbb{Q}(\sqrt[3]{2}\omega)$, $\omega := \exp(\frac{2\pi\sqrt{-1}}{3})$ and so  $N = 2$. Let $\xi$ be the isomorphism $k \stackrel{\sim}{\rightarrow} k'$ defined by $\xi(\sqrt[3]{2}) := \sqrt[3]{2}\omega$.  Noting $2{\cal O}_k = (\sqrt[3]{2})^2, X^3 -2 = (X-4)(X-7)(X-20)$ mod $31$, let $S := \{ \frak{p}_1 := (\sqrt[3]{2}), \frak{p}_2 := (31, \sqrt[3]{2} -4), \frak{p}_2 := (31, \sqrt[3]{2} - 7), \frak{p}_4 := (31, \sqrt[3]{2} - 20)\}, S' := \xi(S) = \{ \frak{p}_1' := (\sqrt[3]{2}\omega), \frak{p}_2' := (31, \sqrt[3]{2}\omega -4), \frak{p}_3' := (31, \sqrt[3]{2}\omega - 7), \frak{p}_4' := (31, \sqrt[3]{2}\omega - 20)\}$, so that we have $k_{\frak{p}_1} = k'_{\frak{p}'_1} = \mathbb{Q}_2$ and $k_{\frak{p}_i} = k'_{\frak{p}'_i} = \mathbb{Q}_{31}$ ($2\leq i \leq 4$). So this example satisfies the above conditions.

The isomorphism $\xi : \overline{X}_S \stackrel{\sim}{\rightarrow}  {\overline{X}}_{S'}$  induces the bijection $ \xi^* : {\cal F}_{\overline{X}_{S'}}  \stackrel{\sim}{\longrightarrow} {\cal F}_{\overline{X}_S}$. By the constructions in the subsection 2.2 and the section 3, we have the following\\
\\
{\bf Proposition 4.2.1.} {\em The isomorphism  $\xi_S : {\cal L}_S \stackrel{\sim}{\rightarrow} {\cal L}_{S'}$ induces the bijection
$$\Gamma_G({\cal F}_{\overline{X}_S}, {\rm res}^{*}_{S}({\cal L}_S)) \stackrel{\sim}{\longrightarrow}  \Gamma_G({\cal F}_{\overline{X}_{S'}}, {\rm res}^{*}_{S'}({\cal L}_{S'}))  $$
which sends $CS_{\overline{X}_S}$ to $CS_{\overline{X}_{S'}}$. The isomorphism $\xi_{S , F} :  L_S \stackrel{\sim}{\rightarrow}  L_{S'}$ induces the isomorphism
$${\cal H}_S \stackrel{\sim}{\longrightarrow} {\cal H}_{S'},$$
which sends $Z_{\overline{X}_S}$ to $Z_{\overline{X}_{S'}}.$ }
\\
\\
{\bf Remark 4.2.2.} Proposition 4.2.1 may be regarded as an arithmetic analogue of the axiom in $(2+1)$-diemnsional TQFT, which asserts that an orientation homeomorphism $f : \Sigma \stackrel{\approx}{\rightarrow} \Sigma'$ between closed surfaces induces an isomorphism ${\cal H}_{\Sigma} \stackrel{\sim}{\rightarrow} {\cal H}_{\Sigma'}$ of quantum Hilbert spaces and if $f$ extends to an orientation preserving homeomorphism $M \stackrel{\approx}{\rightarrow} M'$, with $\partial M = \Sigma, \partial M' = \Sigma'$, $Z_{M}$ is sent to $Z_{M'}$ under the induced isomorphism ${\cal H}_{\partial M} \stackrel{\sim}{\rightarrow} {\cal H}_{\partial M'}$.
\\
\\
{\bf 4.3. The case that $S$ is empty.} In the theory in Sections 2 and 3, we can include the case that $S$ is the empty set $\emptyset$ as follows. 

We define ${\cal F}_{\emptyset}$ to be the space of a single point, ${\cal F}_{\emptyset} := \{ * \}$. We define the arithmetic prequantization principal $\mathbb{Z}/N\mathbb{Z}$-bundle ${\cal L}_{\emptyset}$ to be $\mathbb{Z}/N\mathbb{Z}$, on which $G$ acts trivially, so that the map $\varpi_{\emptyset} : {\cal L}_{\emptyset} \rightarrow {\cal F}_{\emptyset}$ is $G$-equivariant. So the arithmetic prequantization $F$-line bundle $L_{\emptyset}$ is defined by $\mathbb{Z}/N\mathbb{Z} \times_{\mathbb{Z}/N\mathbb{Z}} F = F$. The arithmetic Chern-Simons $1$-cocycle $\lambda_{\emptyset}$ is defined to be $0$. 

Let $\tilde{\Pi}_k$ be the modified \'{e}tale fundamental group of $\overline{X}_k$ defined by considering the Artin-Verdier topology on $\overline{X}_k$, which takes the real primes into account (cf. [H; 2.1], [AC], [Bi], [Z]). It is the Galois group of the maximal extension of $k$ unramified at all finite and infinite primes.  We set
$$ {\cal F}_{\overline{X}_k} := {\rm Hom}_{\rm cont}(\tilde{\Pi}_k, G).$$
Following [H], we define the mod $N$ {\em arithmetic Chern-Simons invariant} $CS_{\overline{X}_k}(\rho)$ of $\rho \in {\cal F}_{\overline{X}_k}$ by the image of $c$ under the composition
$$ H^3(G, \mathbb{Z}/N\mathbb{Z}) \stackrel{\rho^{*}}{\rightarrow} H^3(\tilde{\Pi}_k, \mathbb{Z}/N\mathbb{Z}) \rightarrow H^3(\overline{X}_k, \mathbb{Z}/N\mathbb{Z}) \simeq \mathbb{Z}/N\mathbb{Z},$$
where the cohomology group of $\overline{X}_k$ is the modified \'{e}tale cohomology defined in the Artin-Verdier topology. Thus we have the arithmetic Chern-Simons functional $CS_{\overline{X}_k} : {\cal F}_{\overline{X}_k} \rightarrow \mathbb{Z}/N\mathbb{Z}$ and so we see that 
$$CS_{\overline{X}_k} \in \Gamma_G({\cal F}_{\overline{X}_k}, {\rm res}_{\emptyset}^*({\cal L}_{\emptyset})) = {\rm Map}({\cal M}_{\overline{X}_k}, \mathbb{Z}/N\mathbb{Z}),$$
where ${\rm res}_{\emptyset}$ is the (unique) restriction map ${\cal F}_{\overline{X}_k} \rightarrow {\cal F}_{\emptyset}$. Then we have
$$ dCS_{\overline{X}_k} = 0 = {\rm res}_{\emptyset}^{*}(\lambda_{\emptyset}).$$
 
The arithmetic quantum space ${\cal H}_{\emptyset}$ is defined by $\Gamma_G({\cal F}_{\emptyset}, L_{\emptyset}) = F$. Following [H], we define the arithmetic Dijkgraaf-Witten invariant $Z(\overline{X}_k)$ of $\overline{X}_k$  by
$$ \displaystyle{  Z(\overline{X}_k) := \frac{1}{\# G} \sum_{\rho \in {\cal F}_{\overline{X}_k}} \zeta_N^{CS_{\overline{X}_k}(\rho)}   }$$
and the arithmetic Dijkgraaf-Witten partition function by $Z_{\overline{X}_k} : {\cal F}_{\emptyset} \rightarrow F$ by $Z_{\overline{X}_k}(*) := Z(\overline{X}_k)$ for $* \in {\cal F}_{\emptyset}$. So we have
$$ Z_{\overline{X}_k} \in {\cal H}_{\emptyset}.$$
We note that when $[c]$ is trivial,  $Z({\overline{X}_k})$ coincides with the (averaged) number of continuous homomorphism from $\tilde{\Pi}_k$  to $G$:
$$ Z(\overline{X}_k) = \frac{ \# {\rm Hom}_{\rm cont}(\tilde{\Pi}_k,G) }{\# G}, $$
which is the classical invariant for the number field $k$.\\
\\
{\bf 4.4. Disjoint union of finite sets of primes and reversing the orientation of $\partial V_S$.} Let $S_1 = \{ \frak{p}_1, \dots, \frak{p}_{r_1} \}$ and $S_2 = \{ \frak{p}_{r_1+1}, \dots , \frak{p}_r \}$ be disjoint sets of finite primes of $k$ and let $S = S_1 \sqcup S_2$. We include the case where $S_1$ is empty, but $S_2$ is non-empty.  (For the case where $S_1$ and $S_2$ are both empty, the following arguments are trivial.) Then we have
$$ {\cal F}_S = {\cal F}_{S_1} \times {\cal F}_{S_2}.$$ 

For the arithmetic quantization principal $\mathbb{Z}/N\mathbb{Z}$-bundles, we define the map
$$ \boxplus  :  {\cal L}_{S_1} \times {\cal L}_{S_2} \longrightarrow {\cal L}_S,$$
 as follows. For the case that $S_1 = \emptyset$ (and so $S_2 = S$), we set
$$ m \boxplus [\alpha_{S_2}] := [\alpha_{S_2}].m   \leqno{(4.4.1)}$$
for $(m, [\alpha_S]) \in {\cal L}_{\emptyset} \times {\cal L}_{S_2}$. For the case that $S_1 \neq \emptyset$, we set
$$ [\alpha_{S_1}] \boxplus [\alpha_{S_2}] := [(\alpha_{S_1}, \alpha_{S_2})]   \leqno{(4.4.2)}$$
for $([\alpha_{S_1}],[\alpha_{S_2}]) \in {\cal L}_{S_1} \times {\cal L}_{S_2}$. 

For the arithmetic quantization $F$-line bundles, we let $p_i^{*}(L_{S_i})$ be the $G$-equivariant $F$-line bundle over ${\cal F}_S$ induced from $L_{S_i}$ by the projection $ p_i : {\cal F}_S \rightarrow {\cal F}_{S_i}$ for $i = 1, 2$:
$$ p_i^{*}(L_{S_i}) := \{ ( \rho_S, [([\alpha_{S_i}], z_i)]) \in {\cal F}_S \times L_{S_i} \, | \, \rho_{S_i} = \varpi_{S_i}([\alpha_{S_i}]) \; \}$$
for $\rho_S = (\rho_{S_1}, \rho_{S_2})$.  When $S_1 = \emptyset$, we think of $p_i^{*}(L_{\emptyset}) = F$ simply over ${\cal  F}_{\emptyset} = \{*\}$. Let
$$ p_i^{*}(\varpi_{S_i}) : p_i^{*}(L_{S_i}) \longrightarrow {\cal F}_S $$
be the projection. The fiber over $\rho_S = (\rho_{S_1}, \rho_{S_2})$ is given by
$$ \begin{array}{ll} p_i^{*}(\varpi_{S_i})^{-1}(\rho_S) & = \{ \rho_S \} \times \{ [ ([\alpha_{S_i}],z_i)] \in L_{S_i} \; |\; \rho_{S_i} = \varpi_{S_i}([\alpha_{S_i}]), z_i \in F \}\\
                 & = L_{S_i}(\rho_{S_i})\\
                 & \simeq F,
                 \end{array}$$
where $L_{S_i}(\rho_{S_i})$ is as in (2.1.38). We set
$$ L_{S_1} \boxtimes L_{S_2} := p_1^{*}(L_{S_1}) \otimes p_2^{*}(L_{S_2}),$$
which is the $F$-line bundle over ${\cal F}_S$ and whose element is written by
$$ (\rho_S, [ ([\alpha_{S_1}],z_1)] \otimes [ ([\alpha_{S_2}],z_2)]),$$
where $\rho_S = (\rho_{S_1}, \rho_{S_2}) \in {\cal F}_S$, $[([\alpha_{S_i}],z_i)] \in L_{S_i}(\rho_{S_i})$. The right action  on $L_{S_1} \boxtimes L_{S_2}$ is defined by
$$ (\rho_S, [ ([\alpha_{S_1}],z_1)] \otimes [ ([\alpha_{S_2}],z_2)]).g := (\rho_S.g, [ ([\alpha_{S_1}].g, z_1)] \otimes [ ([\alpha_{S_2}].g, z_2)])$$
so that the projection $L_{S_1} \boxtimes L_{S_2} \rightarrow {\cal F}_S$ is $G$-equivariant. 
Then, as in Proposition 2.1.42, we have the isomorphism of $G$-equivariant $F$-line bundles over ${\cal F}_S$:
$$ L_{S_1} \boxtimes L_{S_2} \stackrel{\sim}{\longrightarrow} L_S; \; (\rho_S, [ ([\alpha_{S_1}],z_1)] \otimes [ ([\alpha_{S_2}],z_2)]) \mapsto [([\alpha_S],z_1z_2)],$$
where $\alpha_S = (\alpha_{S_1},\alpha_{S_2})$. Choose $x_{S_i} \in \Gamma({\cal F}_{S_i}, {\cal L}_{S_i})$ and let $x_S := [(x_{S_1},x_{S_2})] \in \Gamma({\cal F}_S,{\cal L}_S)$. Then we see that
$$ \lambda_{S_1}^{x_{S_1}}(g, \rho_{S_1}) + \lambda_{S_2}^{x_{S_2}}(g, \rho_{S_2}) = \lambda_S^{x_S}(g,\rho_S)$$
for $g \in G, \rho_S = (\rho_{S_1},\rho_{S_2})$ and,  as in the case that $L_{S}$, we have the isomorphism
$$ L_{S_1}^{x_{S_1}} \boxtimes L_{S^2}^{x_{S_2}} := p_1^*(L_{S_1}^{x_{S_1}}) \otimes p_2^*(L_{S_2}^{x_{S_2}}) \stackrel{\sim}{\longrightarrow} L_S^{x_S}; \;\; ((\rho_{S_1},\rho_{S_2}), z_1\otimes z_2) \mapsto (\rho_S, z_1z_2)$$
for $\rho_S = (\rho_{S_1},\rho_{S_2})$, which is compatible with $L_{S_1} \boxtimes L_{S_2} \simeq L_S$ via Proposition 2.1.39.\\
\\
{\bf Proposition 4.4.3.} {\em For $\theta_i \in {\cal H}_{S_i}^{x_{S_i}}$ ($i=1,2$), we define $\theta_1\cdot \theta_2 \in {\cal H}_S^{x_S}$ by
$$ (\theta_1 \cdot\theta_2)(\rho_S) := \theta_1(\rho_{S_1}) \theta_2(\rho_{S_2})$$
for $\rho_S = (\rho_{S_1}, \rho_{S_2})$. Then we have the following isomorphism of $F$-vector spaces
$$ {\cal H}_{S_1}^{x_{S_1}} \otimes {\cal H}_{S_2}^{x_{S_2}} \stackrel{\sim}{\longrightarrow} {\cal H}_{S}^{x_S};\; \theta_1 \otimes \theta_2 \mapsto \theta_1 \cdot \theta_2.$$ 
For $\theta_i \in {\cal H}_{S_i}$ ($i=1,2$), we define $\theta_1 \boxtimes \theta_2 \in {\cal H}_S$ by
$$ (\theta_1 \boxtimes \theta_2)(\rho_S) := p_1^{*}(\theta_1(\rho_{S_1}))\otimes p_2^{*}(\theta_2(\rho_{S_2}))$$
for $\rho_S = (\rho_{S_1}, \rho_{S_2})$. Here $p_1^{*}(\theta_1(\rho_{S_1}))\otimes p_2^{*}(\theta_2(\rho_{S_2}))$ denotes $[([\alpha_S], z_1z_2)]$ when $\theta_i(\rho_{S_i}) = [([\alpha_{S_i}], z_i)], \alpha_S = (\alpha_{S_1}, \alpha_{S_2})$. Then we have the following isomorphism of $F$-vector spaces
$$ {\cal H}_{S_1} \otimes {\cal H}_{S_2} \stackrel{\sim}{\longrightarrow} {\cal H}_S; \; (\theta_1, \theta_2) \mapsto \theta_1 \boxtimes \theta_2.$$
The above isomorphisms are compatible via the isomorphisms $\Theta^{x_{S_i}} : {\cal H}_{S_i} \simeq {\cal H}_{S_i}^{x_{S_i}}$ $(i=1,2)$ and $\Theta^{x_S} : {\cal H}_S \simeq {\cal H}_S^{x_S}$ in (3.1.2). }\\
\\
{\em Proof.} We may assume by Remark 4.1.2 that the cocycle $c$ is normalized.  For $\theta \in {\cal H}_S^{x_S}$, set $\theta_1(\rho_{S_1}) := \theta(\rho_{S_1},1)$ and $\theta_2(\rho_{S_2}) := \theta(1,\rho_{S_2})$. Since $c$ is normalized, by (2.1.7) and (2.1.10), we have $\lambda_{\frak{p}}(g,1) = 0$ for $g \in G$ and $\frak{p} \in S_i$. From this, we have $\theta_i \in {\cal H}_{S_i}^{x_{S_i}}$.
Then the map ${\cal H}_S^{x_S} \rightarrow {\cal H}_{S_1}^{x_{S_1}}\otimes {\cal H}_{S_2}^{x_{S_2}}$; $\theta \mapsto \theta_1 \otimes \theta_2$, gives the inverse of the former map. By the definitions, the second map is compatible with the first one via $\Theta_{S_i}^{x_{S_i}} : {\cal H}_{S_i} \simeq {\cal H}_{S_i}^{x_{S_i}}$ $(i=1,2)$ and $\Theta^{x_S} : {\cal H}_S \simeq {\cal H}_S^{x_S}$  and so we have the following commutative diagram
$$ \begin{array}{rcc}
{\cal H}_{S_1} \otimes {\cal H}_{S_2} & \longrightarrow & {\cal H}_S \\
\Theta^{x_{S_1}}\otimes \Theta^{x_{S_2}}\; \wr \downarrow \;\;\;\; \;\;\;\; & & \;\;\;\;\;\;\;\; \downarrow \wr \; \Theta^{x_S}\\
 {\cal H}_{S_1}^{x_{S_1}} \otimes {\cal H}_{S_2}^{x_{S_2}} & \stackrel{\sim}{\longrightarrow} & {\cal H}_S^{x_S}, \\
\end{array}
$$
from which the second isomorphism follows. $\;\;\Box$
\\
\\
{\bf Remark 4.4.4.} Proposition 4.4.3 may be regarded as an arithmetic analog of the multiplicative property that ${\cal H}_{\Sigma_1 \sqcup \Sigma_2} = {\cal H}_{\Sigma_1} \otimes {\cal H}_{\Sigma_2}$ for disjoint surfaces $\Sigma_1$ and $\Sigma_2$ which is one of the axioms required in $(2+1)$-dimensional TQFT ([A1]).
\\
\\
For a finite prime $\frak{p}$ of $k$, the canonical isomorphism
$$ {\rm inv}_{\frak{p}} :  H^2_{\mbox{\'{e}t}}(\partial V_{\frak{p}}, \mathbb{Z}/N\mathbb{Z}) \stackrel{\sim}{\longrightarrow} \mathbb{Z}/N\mathbb{Z} $$
indicates  that  $\partial V_{\frak{p}}$ is ``orientable" and we choose (implicitly) the ``orientation" of $\partial V_{\frak{p}}$ corresponding $1 \in \mathbb{Z}/N\mathbb{Z}$. We let $\partial V^{*}_{\frak{p}} = \partial V_{\frak{p}}$ with the ``opposite orientation", namely, ${\rm inv}_{\frak{p}}([\partial V^*_{\frak{p}}]) = -1$. 

The arithmetic prequantization principal $\mathbb{Z}/N\mathbb{Z}$-bundle  for $\partial V_{\frak{p}}^*$, denoted by  ${\cal L}_{\frak{p}^*}$, is defined (formally) by ${\cal L}_{\frak{p}}$ with the opposite action of the structure group $\mathbb{Z}/N\mathbb{Z}$, $(\alpha_{\frak{p}},m) \mapsto \alpha_{\frak{p}}.(-m)$ for $\alpha_{\frak{p}} \in {\cal L}_{\frak{p}^*}$ and $m \in \mathbb{Z}/N\mathbb{Z}$. So the arithmetic prequantization $F$-line bundle $L_{\frak{p}^*}$ for $\partial V_{\frak{p}}^*$ is the dual bundle of $L_{\frak{p}}$,  $L_{\frak{p}^*} = L_{\frak{p}}^*$. Noting $\Gamma({\cal F}_{\frak{p}}, {\cal L}_{\frak{p}^*}) = \Gamma({\cal F}_{\frak{p}}, {\cal L}_{\frak{p}})$, the arithmetic Chern-Simons $1$-cocycle $\lambda_{\frak{p}^*}^{x_{\frak{p}}}$ for $\partial V_{\frak{p}}^*$ is given by  $- \lambda_{\frak{p}}^{x_{\frak{p}}}$ for $x_{\frak{p}} \in \Gamma({\cal F}_{\frak{p}}, {\cal L}_{\frak{p}^*})$. The actions of $G$ on ${\cal L}_{\frak{p}^*}^{x_{\frak{p}}} = {\cal F}_{\frak{p}} \times \mathbb{Z}/N\mathbb{Z}$ and $L_{\frak{p}^*}^{x_{\frak{p}}} = {\cal F}_{\frak{p}} \times F$ are changed to those via $\lambda_{\frak{p}^*}^{x_{\frak{p}}}$. 

For a finite set of finite primes $S = \{ \frak{p}_1, \dots , \frak{p}_r \}$, we set $\partial V_S^* := \partial V_{\frak{p}_1}^* \sqcup \cdots \sqcup \partial V_{\frak{p}_r}^*$. Then the arithmetic prequantization bundles ${\cal L}_{S^*}, L_{S^*}, {\cal L}_{S^*}^{x_S}$ and $L_{S^*}^{x_S}$ ($x_S \in \Gamma({\cal F}_{S}, {\cal L}_{S^*}) = \Gamma({\cal F}_{\frak{p}}, {\cal L}_{S})$) are defined in the similar manner. For the arithmetic Chern-Simons 1-cocycle, we have
$$ \lambda_{S^*}^{x_S} = - \lambda_S^{x_S}.$$

Let  ${\cal H}_{S^*}^{x_S}$ be the arithmetic quantum space for $\partial V_S^{*}$ with respect to $x_S$. Then we see that
$$ \begin{array}{ll} {\cal H}_{S^*}^{x_S} & = \{ \theta^* : {\cal F}_S \rightarrow F \; | \; \theta^*(\rho_S.g)  = \zeta_N^{\lambda_{S^*}^{x_S}(g,\rho_S)}\theta^*(\rho_S) \; \mbox{for}\; \rho_S \in {\cal F}_S, g \in G \}\\
& = \{ \theta^* : {\cal F}_S \rightarrow F \; | \; \theta^*(\rho_S.g)  = \zeta_N^{-\lambda_S(g,\rho_S)}\theta^*(\rho_S) \; \mbox{for}\; \rho_S \in {\cal F}_S, g \in G \}\\
&  = \overline{{\cal H}}_S^{x_S}, 
\end{array}
$$
where $\overline{{\cal H}}_S^{x_S}$ is the complex conjugate of ${\cal H}_S^{x_S}$. Since the pairing 
$$ {\cal H}_{S^*}^{x_S} \times {\cal H}_S^{x_S} \longrightarrow F; \; (\theta^*, \theta) \mapsto \displaystyle{  \sum_{\rho_S \in {\cal F}_S} \theta^*(\rho_S) \theta(\rho_S) }$$
is a (Hermitian) perfect pairing, together with (3.1.2), we have the following \\
\\
{\bf Proposition 4.4.5.} {\em ${\cal H}_{S^*}^{x_S}$ and ${\cal H}_{S^*}$ are the dual spaces of ${\cal H}_S^{x_S}$ and ${\cal H}_S$, respectively}:
$$ {\cal H}_{S^*}^{x_S} = ({\cal H}_S^{x_S})^*,\;\; {\cal H}_{S^*} = ({\cal H}_S)^*.$$  
\\
{\bf Remark 4.4.6.} Proposition 4.4.5 may be regarded as an arithmetic analog of the involutory property that ${\cal H}_{\Sigma^*} = {\cal H}_{\Sigma}^*$, where  $\Sigma^* = \Sigma$ with the opposite orientation,  which is one of the axioms required in $(2+1)$-dimensional TQFT ([A1]).\\

In the subsection 2.2 and the section 3, we have chosen implicitly the orientation of $\overline{X}_S$ so that the boundary $\partial \overline{X}_S$ with induced orientation may be identified with $\partial V_S$. Let $\overline{X}_S^{*}$ denote $\overline{X}_S$ with the opposite orientation. Then,  the arithmetic Chern-Simons functional and the Dijkgraaf-Witten partition function for $\overline{X}_S^{*}$ are given as follows:
$$  CS_{\overline{X}_S^{*}}^{x_S} = - CS_{\overline{X}_S}^{x_S}, \;\; Z_{\overline{X}_S^{*}}^{x_S}(\rho_S) = \frac{1}{\# G} \sum_{\rho \in {\cal F}_{\overline{X}_S}} \zeta_N^{- CS_{\overline{X}_S}^{x_S}(\rho)}.	\leqno{(4.4.7)}$$
\\
\begin{center}
 {\bf 5. Decomposition and gluing formulas}
\end{center}

In this section, we show a decomposition formula for arithmetic Chern-Simons invariants and a gluing formula for arithmetic  Dijkgraaf-Witten partition functions, which generalize the decomposition formula in [CKKPY] in our framework. We keep the same notations and assumptions as in Sections 2, 3 and 4.\\
\\
{\bf 5.1. Arithmetic Chern-Simons functionals and arithmetic Dijkgraaf-Witten partition functions for} $V_S$. 
For a  finite prime $\frak{p}$ of $k$, let ${\cal O}_{\frak{p}}$ denote the ring of $\frak{p}$-adic integers and we let $V_{\frak{p}} := {\rm Spec}({\cal O}_{\frak{p}})$. For a non-empty finite set of finite primes $S = \{ \frak{p}_1, \cdots , \frak{p}_r \}$ of $k$, let $V_S := V_{\frak{p}_1} \sqcup \cdots \sqcup V_{\frak{p}_r}$, which plays a role analogous to a tubular neighborhood of a link, and so $\partial V_S$ plays a role of the boundary of $V_S$.   In this subsection, we introduce the arithmetic Chern-Simons functional and arithmetic Dijkgraaf-Witten partition function for $V_S$, which will be used for our gluing formula in the next section.

Let $\tilde{\Pi}_{\frak{p}}$ be the \'{e}tale fundamental group of $V_{\frak{p}}$, namely, the Galois group of the maximal unramified extension of $k_{\frak{p}}$ and we set 
$${\cal F}_{V_{\frak{p}}} := {\rm Hom}_{\rm cont}(\tilde{\Pi}_{\frak{p}}, G),\;\; {\cal F}_{V_{S}} := {\cal F}_{V_{\frak{p}_1}} \times \cdots \times  {\cal F}_{V_{\frak{p}_r}}.$$
Since $\tilde{\Pi}_{\frak{p}} \simeq \hat{\mathbb{Z}}$ (profinite infinite cyclic group), ${\cal F}_{V_{\frak{p}}} \simeq G$. $G$ acts on ${\cal F}_{V_{S}}$ from the right by
$$ {\cal F}_{V_S} \times G \rightarrow {\cal F}_{V_S}; \;\; ((\tilde{\rho}_{\frak{p}_i})_i, g) \mapsto \rho.g := (g^{-1} \tilde{\rho}_{\frak{p}_i} g)_i,  $$
and let ${\cal M}_{V_S}$ denote the quotient set by this action:
$$ {\cal M}_{V_S} := {\cal F}_{V_S}/G.$$
Let $\tilde{{\rm res}}_{\frak{p}_i} : {\cal F}_{V_{\frak{p}_i}} \rightarrow {\cal F}_\frak{p} $ and $\tilde{{\rm res}}_S := (\tilde{\rm res}_{\frak{p}_i}) : {\cal F}_{V_S} \rightarrow {\cal F}_{S} $ denote the restriction maps induced by the natural continuous homomorphisms $v_{\frak{p}_i} : \Pi_{\frak{p}_i} \rightarrow \tilde{\Pi}_{\frak{p}_i}$ $(1 \leq i \leq r)$,   which are $G$-equivariant. We denote by $\tilde{{\rm Res}}_{\frak{p}_i}$ and $\tilde{{\rm Res}}_S$  the homomorphisms on cochains given as the pull-back by $v_{\frak{p}_i}$:
$$ \begin{array}{l}
\tilde{{\rm Res}}_{\frak{p}_i} : C^n(\tilde{\Pi}_{\frak{p}_i},\mathbb{Z}/N\mathbb{Z}) \longrightarrow C^n(\Pi_{\frak{p}_i}, \mathbb{Z}/N\mathbb{Z});  \alpha_i \mapsto \alpha_i \circ v_{\frak{p}_i},\\
\tilde{{\rm Res}}_S := (\tilde{{\rm Res}}_{\frak{p}_i}) : \prod_{i=1}^r  C^n(\tilde{\Pi}_{\frak{p}_i},\mathbb{Z}/N\mathbb{Z}) \longrightarrow \prod_{i=1}^r C^n(\Pi_{\frak{p}_i}, \mathbb{Z}/N\mathbb{Z}); (\alpha_i) \mapsto (\alpha_i \circ v_{\frak{p}_i}).\\
\end{array}$$
For $ \tilde{\rho}=(\tilde{\rho}_{{\frak{p}_i}})_i \in {\cal F}_{V_S}$, $c \circ \tilde{\rho}_{{\frak{p}_i}} \in Z^3(\tilde{\Pi}_{{\frak{p}_i}}, \mathbb{Z}/N\mathbb{Z})$. Since $H^3(\tilde{\Pi}_{\frak{p}_i},\mathbb{Z}/N \mathbb{Z})=0$, there is $\tilde{\beta}_{{\frak{p}_i}} \in C^2(\tilde{\Pi}_{{\frak{p}_i}}, \mathbb{Z}/N\mathbb{Z})$ such that
$$ c \circ \tilde{\rho}_{{\frak{p}_i}} = d \tilde{\beta}_{{\frak{p}_i}}.$$
We see that 
$$c \circ  \tilde{{\rm res}}_{{\frak{p}_i}} (\tilde{\rho}_{{\frak{p}_i}}) = d \tilde{{\rm Res}}_{{\frak{p}_i}} (\tilde{\beta}_{{\frak{p}_i}}) $$
for  $1\leq i \leq r$ and we have
$$ [\tilde{{\rm Res}}_{S}((\tilde{\beta}_{{\frak{p}_i}})_i)]  \in {\cal L}_{S}(\tilde{{\rm res}}_S(\tilde{\rho})).$$
Let $\tilde{{\rm res}}_S^{*}({\cal L}_S)$ be the $G$-equivariant principal $\mathbb{Z}/N\mathbb{Z}$-bundle over ${\cal F}_{V_S}$ induced from ${\cal L}_S$ by $\tilde{{\rm res}}_S$:
$$  \tilde{{\rm res}}_S^{*}({\cal L}_S) := \{ (\tilde{\rho}, \alpha_S) \in {\cal F}_{V_S} \times {\cal L}_S \, | \, \tilde{{\rm res}}_S(\tilde{\rho}) = \varpi_S(\alpha_S) \}$$
and let $\tilde{{\rm res}}_S^{*}(\varpi_S)$ be the projection $ \tilde{{\rm res}}_S^{*}({\cal L}_S) \rightarrow {\cal F}_{V_S}$.
We define the {\em arithmetic Chern-Simons functional} $CS_{V_S} : {\cal F}_{V_S} \rightarrow {\cal L}_S$ by
$$ CS_{V_S}(\tilde{\rho}) := [ \tilde{{\rm Res}}_S((\tilde{\beta}_{{\frak{p}_i}})_i)]$$
for $\tilde{\rho} \in {\cal F}_{V_S}$. The value $CS_{V_S}(\tilde{\rho})$ is called the {\it arithmetic Chern-Simons invariant} of $\tilde{\rho}$.\\
\\
{\bf Lemma 5.1.1.}  (1) {\it $CS_{V_S}(\tilde{\rho})$ is independent of the choice of $\tilde{\beta}_{{\frak{p}_i}}$.}\\
(2) $CS_{V_S}$ is a $G$-equivariant section of $\tilde{{\rm res}}_S^{*}(\varpi_S)$:
$$ CS_{V_S} \in \Gamma_G({\cal F}_{V_S}, \tilde{{\rm res}}_S^{*}({\cal L}_S)) = \Gamma({\cal M}_{V_S}, \tilde{{\rm res}}_S^{*}(\overline{{\cal L}}_S)).$$
\\
{\it Proof.} (1) This follows from the fact that the cohomological dimension of $\tilde{\Pi}_{\frak{p}_i}$ is one.\\
(2) The proof of this lemma is almost same as Lemma 2.2.12. (2). $\;\; \Box$\\ 
\\
For a section $x_S = [ (x_{\frak{p}_1}, \dots , x_{\frak{p}_r})] \in \Gamma({\cal F}_S, {\cal L}_S)$, the isomorphism $\Phi_S^{x_S} : {\cal L}_S \stackrel{\sim}{\rightarrow} {\cal L}_S^{x_S}$ induces the isomorphism
$$ \begin{array}{rcl} \tilde{\Psi}^{x_S} : \Gamma_G({\cal F}_{V_S}, \tilde{{\rm res}}_S^{*}({\cal L}_S)) & \stackrel{\sim}{\longrightarrow} & \Gamma_G({\cal F}_{V_S}, \tilde{{\rm res}}_S^{*}({\cal L}_S^{x_S})) = {\rm Map}_G({\cal F}_{V_S},\mathbb{Z}/N\mathbb{Z}); \\
y_S  & \mapsto &  \Phi_S^{x_S} \circ y_S.\\
\end{array}$$ 
We  define the {\em arithmetic Chern-Simons functional} $CS_{V_S}^{x_S} : {\cal F}_{V_S} \rightarrow \mathbb{Z}/N\mathbb{Z}$ with respect to $x_S$ by the image of $ CS_{V_S}$ under $\tilde{\Psi}^{x_S}$. 
\\
\\
{\bf Proposition 5.1.2.} (1) {\em For $\rho \in {\cal F}_{V_S}$, we have}
$$ CS_{V_S}^{x_S}(\tilde{\rho}) = \sum_{i=1}^r (\tilde{{\rm Res}}_{S}(\tilde{\beta}_{{\frak{p}_i}})  - x_{\frak{p}_i}(\tilde{{\rm res}}_{\frak{p}_i}(\tilde{\rho}_{\frak{p}_i}))). $$\\
(2) {\em We have the following equality in } $ C^1(G, {\rm Map}({\cal F}_{V_S}, \mathbb{Z}/N\mathbb{Z}))$
$$ d CS_{V_S} = \tilde{{\rm res}}^*(\lambda_S^{x_S}). $$
\\
{\em Proof.} (1) This follows from the definition of $\tilde{\Psi}^{x_S}$.\\
(2) Since $CS_{V_S} \in {\rm Map}_G({\cal F}_{V_S}, \mathbb{Z}/N\mathbb{Z})$, we have
$$ CS_{V_S}^{x_S}(\tilde{\rho}.g) = CS_{V_S}^{x_S}(\tilde{\rho}) + \lambda_S^{x_S}(g, \tilde{{\rm res}}_S(\tilde{\rho}))$$
for $g \in G$ and $\tilde{\rho} \in {\cal F}_{V_S}$, which means the assertion.  $\;\; \Box$\\
\\
{\bf Proposition 5.1.3.} {\em Let $x'_S \in \Gamma({\cal F}_S, {\cal L}_S)$ be another section, which yields $CS_{V_S}^{x'_S}$ and let $\delta_S^{x_S, x'_S} : {\cal F}_S \rightarrow \mathbb{Z}/N\mathbb{Z}$ be the map in Proposition 2.1.34. 
 Then we have}
$$ CS_{V_S}^{x'_S}(\tilde{\rho}) - CS_{V_S}^{x_S}(\tilde{\rho}) = \delta_S^{x_S,x'_S}(\tilde{{\rm res}}_S(\tilde{\rho})).$$
\\
{\it Proof.} This follows from Proposition 5.1.2. (1) and Lemma 1.1.4. $\;\; \Box$\\
\\
For $\rho_S \in {\cal F}_S$, we define the subset ${\cal F}_{V_S}(\rho_S)$ of ${\cal F}_{V_S}$ by
$$ {\cal F}_{V_S}(\rho_S) := \{ \tilde{\rho} \in {\cal F}_{V_S} \, | \, \tilde{{\rm res}}_S(\tilde{\rho}) = \rho_S \}.$$
 We then define the {\em arithmetic Dijkgraaf-Witten invariant} $Z_{V_S}(\rho_S)$ of $\rho_S$ with respect to $x_S$ by
$$ Z_{V_S}^{x_S}(\rho_S) := \frac{1}{\# G} \sum_{\tilde{\rho} \in F_{V_{S} }(\rho_{S})} {\zeta_N}^{ CS^{x_{{S}}}_{ V_{S} }(\tilde{\rho}) }.$$
\\
{\bf Theorem 5.1.4.} (1) {\em $Z_{V_S}^{x_S}(\rho_S)$ is independent of the choice of $\tilde{\beta}_{\rho_{\frak{p}_i}}$}.\\
(2) {\em We have}
$$ Z_{V_S}^{x_S} \in {\cal H}_S^{x_S}.$$
\\
{\em Proof.} (1) This follows from Proposition 5.1.1 (1). \\
(2) This follows from Proposition 5.1.2. (2). $\;\; \Box$\\ \\
\\
We call $Z_{V_S}^{x_S}$ the {\em arithmetic Dijkgraaf-Witten partition function} for $V_S$ with respect to $x_S$.\\
\\
{\bf Proposition 5.1.5.} For sections $x_S, x'_S \in \Gamma({\cal F}_S, {\cal L}_S)$ we see that
$$ \Theta^{x_S,x'_S} (Z_{V_S}^{x_S})=Z_{V_S}^{x'_S}. $$
{\it Proof.} This follows from Proposition 5.1.3. $\;\; \Box$\\
\\
By the identification (3.1.3), $Z_{V_S}^{x_S}$ defines the element $Z_{V_S}$ of ${\cal H}_S$ which is independent of the choice of $x_S.$ We call it the {\em arithmetic Dijkgraaf-Witten partition function} for $V_S$.\\

In the above, the orientation of  $V_S$ is chosen so that it is compatible with that of $\partial V_S$ as explained in the subsection 4.4. Let $V_S^*$ denote $V_S$ with opposite orientation. Then, following (4.4.7), the arithmetic Chern-Simons functional and the arithmetic Dijkgraaf-Witten partition function are given by
$$  CS_{V_S^{*}}^{x_S} = - CS_{V_S}^{x_S}, \;\; Z_{V_S^{*}}^{x_S}(\rho_S) = \frac{1}{\# G} \sum_{\tilde{\rho} \in {\cal F}_{V_S}(\rho_S)} \zeta_N^{- CS_{\overline{X}_S}^{x_S}(\tilde{\rho})}.	  \leqno{(5.1.6)}$$
\\
{\bf 5.2. Gluing formulas for arithmetic Chern-Simons invariants and gluing formulas for arithmetic Dijkgraaf-Witten partition functions.}  Let $S_1$ and $S_2$ be disjoint sets of finite primes of $k$, where $S_1$ may be empty and $S_2$ is non-empty. We assume that any prime dividing $N$ is contained in $S_2$ if $S_1$ is empty and that any prime dividing $N$ is contained in $S_1$ if $S_1$ is non-empty. We let $S := S_1 \sqcup S_2$. We may think of $\overline{X}_{S_1}$ as the space obtained by gluing $\overline{X}_S$ and $V_{S_2}^{*}$ along $\partial V_{S_2}$. Let  $\eta_S : \Pi_S \rightarrow \Pi_{S_1}$, $\iota_{\frak{p}} : \Pi_{\frak{p}} \rightarrow \Pi_S$, $v_{\frak{p}} : \Pi_{\frak{p}} \rightarrow \tilde{\Pi}_{\frak{p}}$, and   $u_{\frak{p}} :  \tilde{\Pi}_{\frak{p}} \rightarrow \Pi_{S_1}$ be the natural homomorphisms, where $\frak{p} \in S_2$,  so that  we have  $\eta_S \circ \iota_{\frak{p}}= u_{\frak{p}} \circ v_{\frak{p}}$ for $\frak{p} \in S_2$. 
\\
\\
\setlength{\unitlength}{0.8mm}
\begin{picture}(120,28)(0,0)
\put(64,10){\vector(2, -1){16}}
\put(56,11){$\Pi_{\frak{p}}$ }
\put(64, 15){\vector(2,1){16}}
\put(82,0){$\tilde{\Pi}_{\frak{p}}$}
\put(90,2){\vector(2,1){16}}
\put(82,22){$\Pi_S$}
\put(90,24){\vector(2,-1){16}}
\put(108,11){$\Pi_{S_1}$}
\put(68,4){$v_{\frak{p}}$}
\put(68,20){$\iota_{\frak{p}}$}
\put(100,4){$u_{\frak{p}}$}
\put(98,22){$\eta_S$}
\end{picture}
\\
\\
Let  $\boxplus  : {\cal{L}}_{S_1} \times {\cal{L}}_{S_2} \rightarrow \cal{L}_{S}$ be the map defined as in (4.4.1) and (4.4.2).  Now we have the following decomposition formula.\\ 
\\
{\bf Theorem 5.2.1} ({\em Decomposition formula}). For $\rho \in {\rm Hom}_{\rm cont}(\Pi_{S_1}, G)$, we have
$$CS_{\overline{X}_{S_1}}(\rho) \boxplus CS_{V_{S_2}}((\rho \circ u_{\frak{p}})_{\frak{p} \in S_2} ) = CS_{\overline{X}_S}(\rho \circ \eta_S).  $$
\\
{\em Proof.} Case that $S_1= \emptyset$. Although this may be well known, we give a proof for the sake of readers. By the Artin--Verdier Duality for compact support \'etale cohomologies ([Mil; Chapter II. Theorem 3.1]) and modified \'etale cohomologies ([Bi; Theorem 5.1]), we have the following isomorphisms for a fixed $\zeta_N \in \mu_N$,
$$H^3_{\rm comp}(X_S,\mathbb{Z}/N \mathbb{Z}) \cong {\rm Hom}_{X_S}(\mathbb{Z}/N \mathbb{Z}, \mathbb{G}_{m, X_S} )^{\sim} \cong \mu_N(k)^{\sim} \cong \mathbb{Z}/N \mathbb{Z}, $$
$$H^3(\overline{X}_k,\mathbb{Z}/N \mathbb{Z}) \cong {\rm Hom}_{\overline{X}_k}(\mathbb{Z}/N \mathbb{Z}, \mathbb{G}_{m, \overline{X}_k} )^{\sim} \cong \mu_N(k)^{\sim} \cong \mathbb{Z}/N \mathbb{Z}, $$
where $\mathbb{G}_{m, X_S}$ (resp. $\mathbb{G}_{m, \overline{X}_k}$) is the  sheaf of units on $X_S$ (resp. $\overline{X}_k$) and $(-)^{\sim}$ is given by ${\rm Hom}(-,\mathbb{Q}/\mathbb{Z})$. We denote the isomorphisms above by ${\rm inv}' :H^3_{\rm comp}(X_S,\mathbb{Z}/N \mathbb{Z}) \rightarrow  \mathbb{Z}/N \mathbb{Z}$ and  ${\rm inv} :H^3(\overline{X}_k,\mathbb{Z}/N \mathbb{Z}) \rightarrow \mathbb{Z}/N \mathbb{Z}$. Now we recall the definition of $H^3_{\rm comp}(X_S,\mathbb{Z}/N \mathbb{Z})$ ([Mil; p.165]). We define the complex $C_{comp} (\Pi_S,\mathbb{Z}/N \mathbb{Z})$ by
$$C^n_{comp} (\Pi_S,\mathbb{Z}/N \mathbb{Z}) :=C^n(\Pi_S,\mathbb{Z}/N \mathbb{Z}) \times  \prod_{\frak{p}\in S} C^{n-1} (\Pi_{\frak{p}} ,\mathbb{Z}/N\mathbb{Z}), $$
$$ d(a ,(b_\frak{p}) ) := (da, ({\rm Res}_\frak{p} (a) - d b_\frak{p})), $$
where $a \in C^n(\Pi_S,\mathbb{Z}/N \mathbb{Z})$ and $ (b_\frak{p}) \in \prod_{\frak{p}\in S} C^{n-1} (\Pi_{\frak{p}} ,\mathbb{Z}/N\mathbb{Z}).$ $H^n_{\rm comp}(X_S,\mathbb{Z}/N \mathbb{Z})$ is defined by
$$ H^n_{\rm comp}(X_S,\mathbb{Z}/N \mathbb{Z}) := H^n(C^{*}_{comp} (\Pi_S,\mathbb{Z}/N \mathbb{Z})).$$ Then we can describe ${\rm inv}' :H^3_{\rm comp}(X_S,\mathbb{Z}/N \mathbb{Z}) \rightarrow  \mathbb{Z}/N \mathbb{Z}$ as follows. Let $[(a,(b_{\frak{p}}))] \in  H^3_{\rm comp}(X_S,\mathbb{Z}/N \mathbb{Z})$. Since $da=0$ and $H^3(\Pi_S,\mathbb{Z}/N \mathbb{Z})=0$, there is a cochain $b\in C^2 (\Pi_S,\mathbb{Z}/N \mathbb{Z})$ such that $db=a$. Then we have
$${\rm inv}'  ([(a,(b_\frak{p})]) = {\displaystyle \sum_{\frak{p} \in S} {\rm inv}_\frak{p} ([{\rm Res}_\frak{p} (b) -b_\frak{p} ])}, $$
where ${\rm inv}_\frak{p} : H^2(\Pi_\frak{p},\mathbb{Z}/N \mathbb{Z}) \rightarrow \mathbb{Z}/N \mathbb{Z}$ is the canonical isomorphism given by the theory of Brauer groups. We note that the right side of the equation above doesn't depend on the choice of $b$. Recall that $\tilde{\Pi}_k$ denotes the modified \'{e}tale fundamental group of $\overline{X}_k$. Let $j_3 : H^3(\tilde{\Pi}_k , \mathbb{Z}/ N \mathbb{Z}) \rightarrow H^3(\overline{X}_k, \mathbb{Z}/N \mathbb{Z})$ be the natural homomorphism induced by the modified Hochschild-Serre spectral sequence ([H; Corollary 2.2.8]). We describe the image of the cohomology class $[c \circ \rho ] \in H^3(\tilde{\Pi}_k,\mathbb{Z}/ N\mathbb{Z})$ by the composed map
$${\rm inv'}^{-1} \circ {\rm inv} \circ j_3 : H^3(\tilde{\Pi}_k , \mathbb{Z}/ N \mathbb{Z}) \rightarrow H^3_{\rm comp}(X_S,\mathbb{Z}/N \mathbb{Z}).$$
Since $c \circ (\rho \circ \eta_S) \in Z^3(\Pi_S,\mathbb{Z}/ N \mathbb{Z})$ and $H^3(\Pi_S,\mathbb{Z}/N \mathbb{Z})=0$, there exists a cochain $\beta_{\rho \circ \eta_S} \in C^2(\Pi_S,\mathbb{Z}/N \mathbb{Z})$ such that $d \beta_{\rho \circ \eta_S} =c \circ (\rho \circ \eta_S).$  We note that $d {\rm Res}_\frak{p} (\beta_{\rho \circ \eta_S})=d (\beta_{\rho \circ \eta_S} \circ \iota_\frak{p}) =c \circ \rho \circ u_\frak{p} \circ v_{\frak{p}}$.  Since $c \circ (\rho \circ u_\frak{p}) \in Z^3(\tilde{\Pi}_\frak{p},\mathbb{Z}/ N \mathbb{Z})$ and $H^3(\tilde{\Pi}_\frak{p},\mathbb{Z}/N \mathbb{Z})=H^2(\tilde{\Pi}_\frak{p},\mathbb{Z}/N \mathbb{Z})=0$, there  exists a cochain $\tilde{\beta}_{\rho \circ u_\frak{p}} \in C^2(\tilde{\Pi}_\frak{p},\mathbb{Z}/N \mathbb{Z})$ such that $d \tilde{\beta}_{\rho \circ u_\frak{p}} =c \circ (\rho \circ u_\frak{p})$. We set $\beta_{\rho \circ u_{\frak{p}} \circ v_{\frak{p}}} := \tilde{\beta}_{\rho \circ u_\frak{p}} \circ v_\frak{p} \in C^2(\Pi_\frak{p},\mathbb{Z}/N \mathbb{Z}).$ So we have $d \beta_{\rho \circ u_\frak{p} \circ v_{\frak{p}} } = c \circ (\rho \circ u_\frak{p} \circ v_{\frak{p}}).$ Then we obtain
$$ ({\rm inv'}^{-1} \circ {\rm inv} \circ j_3)([c \circ \rho])=[(c \circ (\rho \circ \eta_S),(\beta_{\rho \circ u_\frak{p}}))].$$
We see that $[{\rm Res}_\frak{p} (\beta_{\rho \circ \eta_S})], [\beta_{\rho \circ u_\frak{p} \circ v_{\frak{p}} }] \in {\cal L}_\frak{p} (\rho \circ u_\frak{p} \circ v_{\frak{p}} )$. Thus we obtain
\begin{eqnarray}
CS_{\overline{X}_k}(\rho)&=& ({\rm inv} \circ j_3)([c\circ \rho]) \nonumber \\
                               &=&  ({\rm inv'} \circ {\rm inv'}^{-1} \circ {\rm inv} \circ j_3)([c \circ \rho]) \nonumber  \\
                               &=& {\rm inv'} ([(c \circ (\rho \circ \eta_S),(\beta_{\rho \circ u_\frak{p} \circ v_{\frak{p}} }))]) \nonumber \\
                               &=& {\displaystyle \sum_{\frak{p} \in S} {\rm inv}_\frak{p} ([{\rm Res}_\frak{p} (\beta_{\rho \circ \eta_S}) -\beta_{\rho \circ u_\frak{p} \circ v_{\frak{p}} } ])} \nonumber \\
                               &=& CS_{\overline{X}_S}(\rho \circ \eta_S)-CS_{V_S}( {(\rho \circ u_\frak{p})}_{\frak{p} \in S}). \nonumber
\end{eqnarray}
Case that $S_1 \neq \emptyset$. Let $\beta_\rho \in C^2(\Pi_{S_1} ,\mathbb{Z} / N \mathbb{Z})$ be a cochain such that $d \beta_\rho=c \circ \rho$. We have $d (\beta_\rho \circ \eta_S)=c \circ (\rho \circ \eta_S)$ and $d (\beta_\rho \circ u_\frak{p})=c \circ (\rho \circ u_\frak{p})$ for $\frak{p} \in S_2$. So we obtain
\begin{eqnarray}
CS_{\overline{X}_{S_1}}(\rho) \boxplus CS_{V_{S_2}}( {(\rho \circ u_\frak{p})}_{\frak{p} \in S_2}) &=& [(\beta_\rho \circ \eta_S \circ \iota_\frak{p})_{ \frak{p} \in S_1} ] \boxplus [(\beta_\rho \circ u_\frak{p} \circ v_\frak{p})_{ \frak{p} \in S_2}] \nonumber \\
                                   &=& [(\beta_\rho \circ u_\frak{p} \circ v_{\frak{p}})_{ \frak{p} \in S}] \nonumber  \\
                                   &=& [(\beta_\rho \circ \eta_S \circ \iota_\frak{p})_{ \frak{p} \in S}] \nonumber \\
                                   &=& CS_{\overline{X}_S}({\rho \circ \eta_S}). \;\;\;\; \Box \nonumber
\end{eqnarray}
\\

 Let $x_{S_i} \in \Gamma({\cal F}_{S_i}, {\cal{L}}_{S_i}) \ (i=1,2))$ be any sections.  We define the section $x_S \in \Gamma({\cal F}_S, \cal{L}_S)$ by
$$x_S(\rho_{S_1}, \rho_{S_2}) := x_{S_1}(\rho_{S_1}) \boxplus x_{S_1}(\rho_{S_2}). $${\cal 
By the proof of Theorem 5.2.1, we have the following \\ 
\\
{\bf Corollary 5.2.2.} Notations being as above, we have the following equality in $\mathbb{Z}/ N \mathbb{Z}$.
$$CS^{x_{S_1}}_{\overline{X}_{S_1}}(\rho)+CS^{x_{S_2}}_{V_{S_2}}((\rho \circ u_{\frak{p}})_{\frak{p} \in S_2} ) = CS^{x_{S}}_{\overline{X}_S}(\rho \circ \eta_S).  $$
\\
We consider the situation that we obtain the space $\overline{X}_{S_1}$  by gluing $\overline{X}_S$ and $V_{S_2}^*$ along $\partial V_{S_2}$.  We define the pairing $<\ ,\ > : {\cal H}^{x_S}_{S} \times {\cal H}^{x_{S_2}}_{{S_2}^*} \rightarrow {\cal H}^{x_{S_1}}_{S_1} $ by
$$<\theta_S , \theta_{{S_2}^*}>(\rho_{S_1}) := \#G \sum_{\rho_{S_2} \in {\cal F}_{S_2}} \theta_S(\rho_{S_1},\rho_{S_2}) \theta_{{S_2}^*}(\rho_{S_2}) \leqno{(5.2.3)}$$
for $\theta_S \in {\cal H}_S^{x_S}, \theta_{S_2^*} \in {\cal H}_{S_2^*}^{x_{S_2}}$ and $\rho_{S_1} \in {\cal F}_{S_1}$.   This induces the pairing  $<\ ,\ > : {\cal H}_{S} \times {\cal H}_{{S_2}^*} \rightarrow {\cal H}_{S_1}$ by (3.1.2).  Now we prove the following gluing formula.\\
\\
{\bf Theorem 5.2.4} {\em (Gluing formula}). Notations being as above, We have the following equality
$$<Z_{{\overline{X}}_S}, Z_{V_{S_2}^*}  > \; = \; Z_{{\overline{X}}_{S_1}}. $$
\\
{\em Proof.} We show the equality
$$<Z^{x_{S}}_{\overline{X}_S}, Z^{x_{S_2}}_{V_{S_2}^*}  > \; = \; Z^{x_{S_1}}_{{\overline{X}}_{S_1}}$$
for any sections $x_{S_i}  \in \Gamma({\cal F}_{S_i}, {\cal{L}}_{S_i}) \ (i=1,2)$. Noting (5.1.6), we have
\begin{eqnarray}
<Z^{x_{S}}_{{\overline{X}}_S}, Z^{x_{S_1}}_{V_{S_2}^* } > (\rho_{S_1})  &=&  \#G \sum_{\rho_{S_2} \in {\cal F}_{S_2}} \Bigl( \frac{1}{\# G} \sum_{\rho' \in {\cal F}_{ {\overline{X}}_S }(\rho_{S_1}, \rho_{S_2})} {\zeta_N}^{ CS^{x_{S}}_{ \overline{X}_{S} }(\rho') } \Bigr)  \Bigl( \frac{1}{\# G} \sum_{\tilde{\rho} \in {\cal F}_{V_{S_2} }(\rho_{S_2})} {\zeta_N}^{ -CS^{x_{{S_2}}}_{ V_{S_2} }(\tilde{\rho}') } \Bigr) \nonumber \\
                                   &=& \sum_{\rho_{S_2} \in {\cal F}_{S_2}} \Bigl( \frac{1}{\# G}  \sum_{ (\rho', \tilde{\rho}) \in {\cal F}_{ {\overline{X}}_S }(\rho_{S_1}, \rho_{S_2}) \times {\cal F}_{V_{S_2} }(\rho_{S_2})} {\zeta_N}^{ CS^{x_{S}}_{ \overline{X}_{S} }(\rho')  -CS^{x_{{S_2}}}_{ V_{S_2} }(\tilde{\rho})} \Bigr)   \nonumber
\end{eqnarray}
for $\rho_{S_1} \in {\cal F}_{S_1}$. We define the map
$$ \chi (\rho_{S_1}) : {\cal F}_{ {\overline{X}}_{S_1} }(\rho_{S_1}) \rightarrow \underset{ \rho_{S_2} \in {\cal F}_{S_2} } {\bigsqcup} \Bigl( {\cal F}_{ {\overline{X}}_S }(\rho_{S_1}, \rho_{S_2}) \times {\cal F}_{V_{S_2} }(\rho_{S_2}) \Bigr)$$
by
$$ \chi (\rho_{S_1})(\rho_1) = (\rho_1 \circ \eta_S, (\rho_1 \circ u_\frak{p})_{\frak{p} \in S_2} )$$
for $\rho_1 \in {\cal F}_{\overline{X}_{S_1}}(\rho_{S_1})$.  In order to obtain the required statement by Corollary 5.2.2, it suffices to show that  $\chi (\rho_{S_1})$ is bijective. (Though this may be seen by noticing that $\Pi_{S_1}$ is the push-out of the maps $\iota_{\frak{p}}$ and $v_{\frak{p}}$ ($\Pi_{S_1}$ is the amalgamated product of $\Pi_S$ and $\tilde{\Pi}_k$ along $\Pi_{\frak{p}}$) for $S_2 = \{ \frak{p} \}$, we give here a straightforward proof.) \\
$\chi(\rho_{S_1})$ is injective: Suppose $\chi(\rho_{S_1})(\rho_1) = \chi(\rho_{S_1})(\rho_1')$ for $\rho_1, \rho_1' \in {\cal F}_{ {\overline{X}}_{S_1} }(\rho_{S_1})$. Then $\rho_1 \circ \eta_S = \rho_1' \circ \eta_S$. Since $\eta_S$ is surjective, $\rho_1 = \rho_1'$.\\
$\chi (\rho_{S_1})$ is surjective: Let $(\rho, (\tilde{\rho}_{\frak{p}})_{\frak{p} \in S_2}) \in {\cal F}_{ {\overline{X}}_S }(\rho_{S_1}, \rho_{S_2}) \times {\cal F}_{V_{S_2} }(\rho_{S_2})$. Then we have 
$${\rm res}_{S_1}(\rho) = \rho_{S_1},   {\rm res}_{S_2}(\rho) = \rho_{S_2}, \tilde{\rm res}_{S_2}( (\tilde{\rho}_{\frak{p}})_{\frak{p} \in S_2}) = \rho_{S_2}. $$ 
Since $\tilde{\rm res}_{\frak{p}}(\tilde{\rho}_{\frak{p}})$ is unramified representation of $\Pi_{\frak{p}}$ for $\frak{p} \in S_2$, $\rho$ is unramified over $S_2$. Therefore there is $\rho_1 \in {\cal F}_{\overline{X}_{S_1}}$ such that $\rho = \rho_1 \circ \eta_S$. Since we see that
$$ \rho_1 \circ u_{\frak{p}} \circ v_{\frak{p}} = \rho_1 \circ \eta_S \circ \iota_{\frak{p}} = \rho \circ \iota_{\frak{p}} = \tilde{\rho}_{\frak{p}} \circ v_{\frak{p}}$$
for $\frak{p} \in S_2$ and $v_{\frak{p}}$ is surjective, we have $\rho_1 \circ u_{\frak{p}} = \tilde{\rho}_{\frak{p}}$ for $\frak{p} \in S_2$. Hence $\chi(\rho_{S_1})(\rho_1) = (\rho, (\tilde{\rho}_{\frak{p}})_{\frak{p} \in S_2})$ and so $\chi(\rho_{S_1})$ is surjective. $\;\; \Box$\\
\\
\vspace{0.4cm}\\
{\bf References}
\begin{flushleft}
{[AC]} E. Ahlqvist, M. Carlson, The cohomology of the ring of integers of a number field. arXiv:1803.08437.\\
{[A1]} M. Atiyah, Topological quantum field theories. Inst. Hautes Etudes Sci. Publ. Math. No. {\bf 68} (1988), 175--186.\\
{[A2]} M. Atiyah, The Geometry and Physics of Knots. Cambridge University Press, 1990.\\
{[BL]} A. Beauville, Y. Laszlo, Conformal blocks and generalized theta functions. Commun. Math. Phys., {\bf 164} (1994), 385--419.\\
{[Bi]} M. Bienenfeld, An etale cohomology duality theorem for number fields with a real embedding. Trans. Amer. Math. Soc. {\bf 303} (1987), no. 1, 71--96. \\
{[BCGKPT]} F.M. Bleher, T. Chinburg, R. Greenberg, M. Kakde, G. Pappas, M. Taylor, Cup products in the \'{e}tale cohomology of number fields. New York J. Math. {\bf 24} (2018), 514--542.\\
{[CKKPPY]} H.-J. Chung, D. Kim, M. Kim, G. Pappas, J. Park, H. Yoo, Abelian arithmetic Chern-Simons theory and arithmetic linking numbers.  Int. Math. Res. Not.  2019, no. {\bf 18}, 5674--5702.\\
{[CKKPY]} H.-J. Chung, D. Kim, M. Kim, J. Park and H. Yoo, Arithmetic Chern-Simons theory II. to appear in Proceedings of the Simons Symposium on $p$-adic Hodge theory (2019). (available at arXiv:1609.03012.)\\
{[DW]} R. Dijkgraaf, E. Witten, Topological gauge theories and group coho-
mology. Commun. Math. Phys. {\bf 129} (1990), 393--429.\\
{[EM]} S. Eilenberg, S. MacLane, Cohomology theory in abstract groups I, Annals of Math., {\bf 48}, (1947), 51--78.\\
{[FQ]} D. Freed, F. Quinn, Chern-Simons theory with finite gauge group. Comm. Math. Phys. {\bf 156} (1993), no.3, 435--472.\\
{[G]} K. Gomi, Extended topological quantum field theory: a toy model. (Japanese) In: Report of the 4th Kinosaki Seminar, 2007, 18 pages.\\
{[H]} H. Hirano, On mod $2$ Dijkgraaf-Witten invariants for certain real quadratic number fileds. arXiv:1911.12964.\\
{[Ki]} M. Kim, Arithmetic Chern-Simons Theory I.  arXiv:1510.05818.\\
{[Ko]} T. Kohno, Conformal field theory and topology. Translations of Mathematical Monographs, {\bf 210}. Iwanami Series in Modern Mathematics. Amer. Math. Soc. Providence, RI, 2002. \\
{[LP]} J. Lee, J. Park, Arithmetic Chern-Simons theory with real places,
 arXiv:1905.13610.\\
{[My]} J. P. May, Simplicial objects in algebraic topology. Van Nostrand Mathematical Studies, No.{\bf 11} D. Van Nostrand Co., Inc., Princeton, N.J.-Toronto, Ont.-London 1967.\\
{[Mz]} B. Mazur, Notes on etale cohomology of number fields. Ann. Sci. Ecole Norm. Sup. (4) 6 (1973), 521--552.\\
{[Mh]} T. Mihara, Cohomological approach to class field theory in arithmetic topology. Canad. J. Math. {\bf 71} (2019), no. 4, 891--935.  \\
{[Mil]} J. S. Milne, Arithmetic Duality  Theorems, Perspectives in Math. Vol. {\bf 1}, Academic Press, 1986.\\
{[Mo1]} M. Morishita, Integral representations of unramified Galois groups and matrix divisors over number fields. Osaka J. Math. {\bf 32} (1995), no. 3, 565--576.\\
{[Mo2]} M. Morishita, Knots and Primes -- An introduction to Arithmetic Topology. Universitext. Springer, London, 2012. \\
{[NSW]} J. Neukirch, A. Schmidt, K. Wingberg, Cohomology of number fields. Second edition. Grundlehren der Mathematischen Wissenschaften, {\bf 323}, Springer-Verlag, Berlin, 2008.\\
{[NU]} H. Niibo, J. Ueki, Id\`{e}lic class field theory for $3$-manifolds and very admissible links, Transactions of the AMS, {\bf 371}, No.12, (2019), 8467--8488.\\
{[S1]} J.-P. Serre, Cohomologie Galoisienne. Lecture Note in Math. {\bf 5}, Springer, 1964.\\
{[S2]} J.-P. Serre, Corps Locaux. Hermann, 1968.\\
{[V]} E. Verlinde, Fusion rules and modular transformations in 2-D conformal field theory. Nucl. phys. {\bf B300[FS22]} (1988), 360-376.\\
{[Wa]} M. Wakui, On Dijkgraaf-Witten invariant for $3$-manifolds. Osaka J. Math. {\bf 29} (1992), 675--696.\\
{[We]} A. Weil, G\'{e}n\'{e}ralisation des fonctions ab\'{e}liennes. J. Math. Pure Appl. {\bf 17} (1938), no.9, 47--87.\\
{[Wi]} E. Witten, Quantum field theory and the Jones polynomial. Commun. Math. Phys. {\bf 121} (1989), 351-399.\\
{[Y]} D. Yetter, Topological quantum field theories associated to finite groups and crossed $G$-sets. J. Knot Theory and its Ramifications. {\bf 1} (1992), 1--20.\\
{[Z]} T. Zink, Etale cohomology and duality in number fields. Appendix 2 In: Haberland, Galois cohomology of algebraic number fields. VEB Deutscher Verlag der Wissenschaften, Berlin, 1978. \\
\end{flushleft}
\vspace{0.4cm}
{\small 
H. Hirano:\\
Graduate School of Mathematics, Kyushu University, \\
744, Motooka, Nishi-ku, Fukuoka  819-0395, Japan.\\
e-mail: h-hirano@math.kyushu-u.ac.jp\\
\\
J. Kim:\\
Graduate School of Mathematics, Kyushu University, \\
744, Motooka, Nishi-ku, Fukuoka  819-0395, Japan.\\
Current address:\\
3-18-3, Megurohoncho, Meguro-ku, Tokyo 152-0002, Japan\\
e-mail: res1235@gmail.com\\
\\
M. Morishita:\\
Graduate School of Mathematics, Kyushu University, \\
744, Motooka, Nishi-ku, Fukuoka  819-0395, Japan.\\
e-mail: morisita@math.kyushu-u.ac.jp \\
}
\end{document}